\documentclass[10pt,a4paper,reqno]{article}
\usepackage{amsfonts}
\usepackage{amsthm}
\usepackage{amsmath}
\usepackage{amscd}
\usepackage[latin2]{inputenc}
\usepackage{t1enc}
\usepackage[mathscr]{eucal}
\usepackage{graphicx}
\usepackage{graphics}
\usepackage{pict2e}
\usepackage{epic}
\usepackage{verbatim}
\usepackage{stmaryrd}
\usepackage{systeme}
\usepackage{amssymb}
\usepackage{mathtools}
\usepackage{subcaption}
\usepackage{textcomp}
\usepackage{changes}
\usepackage{indentfirst}

\newcommand{\mm}[1]{{\color{black}{#1}}}

\numberwithin{equation}{section}
\usepackage[margin=2.3cm]{geometry}

\newcommand{\tro}{\tilde{r}_1}
\newcommand{\trt}{\tilde{r}_2}

\title{A direct sampling method for \mm{simultaneously recovering inhomogeneous inclusions of different nature}}
\date{}
\author{Yat Tin Chow \thanks{Department of Mathematics, University of California, Riverside. The work of this author is supported by Omnibus Research and Travel Award 2019, University of California, Riverside. ({ytchow@ucr.edu}).} 
\and Fuqun Han \thanks{Department of Mathematics, The Chinese University of Hong Kong, Shatin, N.T., Hong Kong.  ({fqhan@math.cuhk.edu.hk}).}
\and Jun Zou \thanks{Department of Mathematics, The Chinese University of Hong Kong, Shatin, N.T., Hong Kong. The work of this author was substantially supported by Hong Kong RGC grant (Projects 14322516 and 14306718). ({zou@math.cuhk.edu.hk}).}}

\begin{document}
\maketitle
\begin{abstract}
\mm{In this work, we investigate a class of elliptic inverse problems and aim to simultaneously recover 
multiple inhomogeneous inclusions arising from two different physical parameters, using very limited boundary Cauchy data
collected only at one or two measurement event. We propose a new fast, stable and highly parallelable direct sampling method (DSM) 
for the simultaneous reconstruction process. Two groups of probing and index functions are constructed, 
and their desired properties are analyzed. 
In order to identify and decouple the multiple inhomogeneous inclusions of different physical nature, 
we introduce a new concept of mutually almost orthogonality property that generalizes 
the important concept of almost orthogonality property in classical DSMs for inhomogeneous inclusions of same 
physical nature 
in \cite{DOT,Chow_2014, chow2018a, Ito_2012, li2013a}. 
With the help of this new concept, we develop a reliable strategy 
to distinguish two different types of inhomogeneous inclusions with noisy data collected at one or two measurement event. 
We further improve the decoupling effect by choosing an appropriate boundary influx. 
Numerical experiments are presented to illustrate the robustness and efficiency of the proposed method.
}
\end{abstract}


\section{Introduction}

In this work, we propose a novel \mm{parameter reconstruction} method in which we decouple measurements from one (or at most two) pair(s) of Cauchy data and locate two different types of inhomogeneities in the model. Let us consider an open bounded domain $\Omega$ in $\mathbb{R}^d$ ($d = 2$, $3$) with \mm{a smooth boundary $\partial \Omega$}, and the following elliptic PDE: 
\begin{equation}
\label{eqn_prostso}
\begin{cases}
  -\nabla \cdot(\sigma \nabla u) + Vu = 0 \quad &\text{in} \quad \Omega\,, \\
  \frac{\partial u}{\partial \nu} = f \quad &\text{on} \quad \partial \Omega\,,
\end{cases}
\end{equation}
where \mm{the coefficients $\sigma$, $V \in L^{\infty}(\Omega)$ represent two unknown physical inclusions 
in the physical ranges} 
$ c < \sigma < C $ and $ -C < V < C $ for some $ c > 0$, $C < \infty$. Let $\sigma_0$ and $V_0$ be the respective coefficients describing the homogeneous background medium $u_0$. 
\mm{We assume two physical inclusions are in the interior of the domain, i.e.,}
$\text{supp}(\sigma - \sigma_0)$, $\text{supp}(V - V_0)$.
Our goal is \mm{to simultaneously identify and reconstruct these two inclusions,  
i.e., $\text{supp}(\sigma - \sigma_0)$ and $\text{supp}(V- V_0)$,  
using the data $u$ measured on the boundary 
corresponding to a boundary influx $f$. 
 } 
\mm{We like to point out that our proposed method can be appropriately generalized to handle other types of boundary conditions that may arise in real applications, e.g. the Robin boundary condition, 
although this work focuses only on a Neumann boundary condition (cf.\,\eqref{eqn_prostso}).}

\mm{Inverse problems of the elliptic system \eqref{eqn_prostso} may arise from a wide range of applications, 
such as medical imaging,  geophysical prospecting, nano-optics, and nondestructive testing; see, e.g., \cite{dot_review, novotny2006principles, schmerr2016fundamentals, geo} and the references therein.} 
\mm{The solution $u$ and two coefficients $\sigma$ and $V$ may represent different physical state and 
parameters in different applications.
For instance, in the diffusion-based optical tomography \cite{arridge1999optical}, 
$u$, $\sigma$ and $V$ represent the photon density, 
diffusion and absorption coefficients, respectively; 
}
\mm{Identification of locations of inhomogeneities of $\sigma$ and $V$ helps determine the distribution of different types of tissues.}
\mm{The model \eqref{eqn_prostso} can also represent the inverse electromagnetic scattering problem. 
Under the transverse electric symmetry, the three-dimensional full Maxwell equations 
may be reduced to \eqref{eqn_prostso}, where $\sigma$ and $-V$ stand for the permeability and permittivity of the media \cite{asym_exp}. 
The system \eqref{eqn_prostso} is also adopted in the ultrasound medical imaging, 
where $\sigma$ and $V$ represent the volumetric mass density and bulk modulus, respectively, 
while $u$ describes the acoustic pressure \cite{ammari2012direct}. 
For the convenience of descriptions, we shall often call $\sigma$ and $V$ as conductivity and potential 
throughout this work.
}

\mm{The uniqueness and simultaneous identifiability for the elliptic inverse problem \eqref{eqn_prostso}) 
have been widely investigated.
}
In particular, a negative result 
was proved in \cite{Arridge_1998}, 
\mm{that is, no uniqueness for the simultaneous reconstruction of $\sigma$ and $V$ 
when both coefficients are smooth. 
For piecewise constant $\sigma$ and piecewise analytic $V$, 
the uniqueness and simultaneous identifiability were established 
in \cite{Harrach_2009} for real-valued coefficients, as long as all possible Neumann-to-Dirichlet data are available. 
This uniqueness result also hints why there are many reasonable numerical results for simultaneous reconstructions, 
even though there is still no general uniqueness result. 
}

\mm{During the recent two decades, many efficient numerical methods were proposed for the inverse problem \eqref{eqn_prostso}. Minimizing a least-squared functional with appropriate regularizations 
is a very popular methodology 
in many applications, along with iterative methods; see, e.g., 
\cite{beilina2015optimization, Dorn_1998, dorn2006level, Kolehmainen_1999, Wang_2010}.
}
\mm{Usually, a locally convergent Newton-type method is employed.  
However, an iterative scheme may be trapped often in local optima, 
owing to high ill-posedness and high non-convexity of the objective functional. 
Moreover, the high dimension of the optimization problem also hinders the performance of this type of algorithms.  
Therefore there is a significant interest to develop some alternative numerical methods, 
that are fast, computationally cheap and robust against noisy data, to provide a reasonable initial guess 
for these iterative methods. On the other hand, some rough estimates of the inhomogeneous inclusions 
directly from the measurement data may be sufficient for many practical applications.
}

\mm{Motivated by these two important applications, 
many non-iterative schemes were developed for a large class of inverse problems for parameter identifications. 
Most of those methods are sampling-type, which rely on an appropriately designed functional that is expected to attain 
relatively large values inside the inhomogeneity. These include linear sampling method \cite{Colton_1996}, singular source method \cite{potthast2001point}, and factorization method \cite{kirsch1999factorization}, etc.
Recently, MUSIC-type method using the multistatic response matrix (MSR) \cite{ammari2005reconstruction, Kirsch_2002}, algorithms based on the topological derivative \cite{Bellis_2013}, and the reverse time migration \cite{chen2013reverse} 
were also developed for the purpose.
We refer to several recent monographs
\cite{cakoni2011linear,xudong_book_2018, kirsch2008factorization, Potthast_2006} 
for more developments in this direction. 
Nevertheless, to the best of our knowledge, 
there seems to exist little development of sampling type methods for simultaneously reconstructing 
two different types of inhomogeneities. 
}

\mm{In this work, we make the first effort to develop a new sampling type method, a direct sampling method (DSM), 
for simultaneously identifying and recovering multiple inhomogeneous inclusions corresponding to 
two different physical parameters. In particular, a specific attempt is made 
to ensure that the method can apply to the important scenarios where very limited data is available, 
e.g., only noisy data collected at one or two measurement event.
}
\mm{DSMs have been developed recently through a series of efforts,
e.g., \cite{DOT,Chow_2014, chow2018a, Ito_2012, li2013a, potthast2001point}, for recovering the inhomogeneous media, 
first for the wave type inverse problems, and then for the non-wave inverse problems. 
}
\mm{This family of direct sampling methods construct an index function that leverages upon an almost orthogonality property 
between the family of fundamental's functions of the forward problem and a particular family of 
probing functions under a properly selected Sobolev duality product.
}
\mm{All the existing DSMs were designed for the cases when there are only inhomogeneous inclusions of same physical nature. 
In this work, we make the first attempt to design DSMs for simultaneously recovering multiple inhomogeneous inclusions 
of two very different physical parameters. These inverse scattering problems are much more ill-posed and challenging than those 
associated with inclusions of same physical nature. A natural mathematical and technical issue 
is how to identify which inclusions come from one physical parameter, not from the other; and how to locate and separate 
the multiple inclusions corresponding to one parameter from those corresponding to the other. 
We shall make use of an important observation that the near field or scattered data satisfies a fundamental property
that it can be approximated as a combination of Green's functions and their gradients at a set of discrete points.  
With this observation, we shall develop two separate families of probing functions, 
namely, the monopole and dipole probing functions, which enable us to construct two separate index functions 
for decoupling the multiple inhomogeneous inclusions associated with one physical parameter 
from those associated with the other parameter. 
In order for this decoupling to function effectively, 
we introduce a new and key concept, the mutually almost orthogonality property, between the family of fundamental functions 
and their gradients, and two families of monopole and dipole probing functions. 
Furthermore, we take advantage of an additional parameter, namely the probing direction of the dipole probing function, 
and an appropriate boundary influx to decouple the multiple inhomogeneous inclusions of one parameter 
from those of the other parameter.
As we will see, the new method is computationally cheap and numerically stable, and works quite satisfactorily, 
as demonstrated in section\,\ref{sec_numerical} by several typical challenging numerical examples 
with very limited data available, e.g., only noisy data collected at one or two measurement event. 
The outputs generated by the new method can serve as reasonable approximations for many important 
applications where general rough locations and shapes of inhomogeneous inclusions are sufficient, 
or as a quick and stable initial guess of some expensive nonlinear optimization approaches 
when more accurate reconstructions are needed. 
}

\mm{The rest of our work is as follows. We address in 
section \ref{sec_principle} the general principles of DSMs, including 
the fundamental property and the new mutually almost orthogonality property.} 
We then show in section \ref{sec_fundamental} that the fundamental property holds for our inverse problem in many cases 
that we encounter in practice. 
\mm{We propose in section \ref{sec_probing} two index functions} for the reconstruction process and discuss 
their properties, including an alternative characterization. In section \ref{sec_explicit}, we derive some explicit representations of the probing and index functions in some special \mm{sampling} domains and 
discuss the mutually almost orthogonality properties in those cases. 
We will also address some appropriate boundary influxes to further decouple the monopole and the dipole effects 
in the measurement. 
Numerical experiments are conducted in section \ref{sec_numerical} to illustrate the effectiveness of the new method.

\section{Principles of DSMs with coupled measurement}

\label{sec_principle}
We briefly explain in this section some general observations that motivate our direct sampling method with coupled measurement.
The development of our DSM hinges on a basic fact that our measurement data can be approximated by a sum of Green's functions of the homogeneous equation and their gradients. 
\mm{With this in mind, along with an appropriate choice of the Sobolev duality product,
those Green's functions and their gradients located at different sampling points are respectively nearly orthogonal with two properly selected families of probing functions.  
These two families of probing functions are monopole-type and dipole-type functions, and 
couple well with the Green's function and its gradient respectively.  
This is a very important property to our new method, and will be called the mutually almost orthogonality property, 
namely, the Green's functions interact well only with monopole probing functions, 
while the gradient of Green's functions interact well solely with dipole probing functions.
This allows us to decouple the monopole and the dipole effects. Moreover, 
different types of boundary influxes and probing directions can be chosen to maximize the decoupling effect.
}

\mm{To be more precise, we aim to make use of the following two properties to develop an effective and robust 
direct sampling method:}
\begin{enumerate}
\item
(Fundamental property) 
\mm{The boundary data, i.e., $u-u_{0}$ on $\partial\Omega$,  of the model \eqref{eqn_prostso} 
can be represented approximately by a sum of Green's functions of the homogeneous medium and their gradients:
}
\[
(u-u_{0})(x) \approx \sum_{j=1}^{n} c_j \, G_{q_j}(x) + \sum_{i=1}^{m} a_i \, d_i\cdot \nabla G_{p_i}(x) \,, \quad 
x\in \partial\Omega
\]
for some choices of coefficients $\{c_j\}_{j=1}^{n} \in \mathbb{C} $, $\{(a_i, d_i)\}_{i=1}^{n} \in \mathbb{C} \times \mathbb{S}^{d-1} $, \mm{and the sets of discrete points} 
$\{q_j \}_{j=1}^n \in \text{supp}(V - V_0) $, $\{p_i\}_{i=1}^m \in \text{supp}(\sigma - \sigma_0) $.

\item
(Mutually almost orthogonality property) 
\mm{There are two sets of probing functions, namely $ \{ \zeta_{x} \}_{x \in \Omega}$ representing a family of monopole probing functions \mm{at sources} $x \in \Omega$, 
and $ \{ \eta_{x,d} \}_{x \in \Omega, d \in \mathbb{S}^{d-1}}$ representing a family of dipole probing functions at 
sources $x \in \Omega$ and dipole directions $d \in \mathbb{S}^{d-1}$, such that 
the following four kernels
}
\begin{eqnarray*}
  (x,z) &\mapsto& {K_1}(x,z) := \frac{ (\zeta_{x}, G_z )_{\text{mo}} }{C_{\text{mo}}(x) }\,, \\  
  (x,z,d_z) &\mapsto& {K_{2,d_z} }(x,z) := \frac{ (\zeta_{x}, d_z \cdot \nabla G_z )_{\text{mo}} }{C_{\text{mo}}(x) }\,, \\
  (x,z,d_x) &\mapsto& {K_{3,d_x}}(x,z) := \frac{ (\eta_{x,d_x}, G_z )_{\text{di}} }{C_{\text{di}}(x, d_x) }\,, \\
  (x,z,d_x,d_z) &\mapsto& {K_{4, d_x, d_z} }(x,z) := \frac{ (\eta_{x,d_x}, d_z \cdot \nabla G_z )_{\text{di}} }{C_{\text{di}}(x,d_x) }\,
\end{eqnarray*}
\mm{have the following properties, 
under two appropriate couplings $( \cdot , \cdot )_{\text{mo}}$, $( \cdot , \cdot )_{\text{di}}$ and weights 
$C_{\text{mo}} (x)$ for $x \in \Omega$ and $C_{\text{di}} (x,d)$ for $x \in \Omega$, $d \in \mathbb{S}^{d-1}$: 
}
\begin{eqnarray*}
  {K_1}(x,z) & &\text{ is of large magnitude if $x$ is close to $z$, and is small otherwise, } \\
  {K_{2,d_z} }(x,z) & &\text{ is relatively small,} \\
  {K_{3,d_x}}(x,z) & &\text{ is relatively small,} \\
  {K_{4, d_x, d_z} }(x,z) & & \text{ is of large magnitude if $x \approx z$ and $d_x \approx d_z$, and is small otherwise.} 
\end{eqnarray*}
\end{enumerate}
\mm{The above mutually almost orthogonality property 
means that the two families of probing functions, i.e., monopole and dipole probing functions, interact well with only the Green's functions and their gradients respectively. This is a very important property that allows us to decouple the monopole and
dipole effects in the measurement data.
}

\mm{With the above definitions and the fundamental property, we can define two index functions}
\begin{eqnarray}
I_{\text{mo}}(x) := \frac{ (\zeta_{x}, u- u_0 )_{\text{mo}} }{C_{\text{mo}}(x) } \,  \quad \text{and}\quad 
I_{\text{di}}(x,d_x) := \frac{ (\eta_{x,d_x}, u-u_0 )_{\text{di}} }{C_{\text{di}}(x, d_x) }\,, \label{eq:dsm}
\end{eqnarray}
\mm{which have the approximations} 
\begin{eqnarray*}
I_{\text{mo}}(x) &\approx & \sum_{j=1}^n c_j K_1(x,q_j) + \sum_{i=1}^m a_iK_{2,d_{i}}(x,p_i) \,, \\
I_{\text{di}}(x,d_x) &\approx& \sum_{j=1}^n c_j K_{3,d_x}(x,q_j) + \sum_{i=1}^m a_i K_{4,d_x, d_{i}} (x,p_i) \,. 
\end{eqnarray*}
\mm{From the above, we can see from the mutually almost orthogonality property that} the index function
$I_{\text{mo}}(x)$ has a large magnitude if $x$ is close to one of the points $\{q_j\}_{j=1}^m$ inside 
the potential inclusions,  i.e., \text{supp}$(V-V_0)$, 
and is small otherwise. Meanwhile, the index function
$I_{\text{di}}(x,d_x) $ has a large magnitude if $x$ is close to one of the points $\{p_i\}_{i=1}^n$ inside  
the conductivity inclusions, i.e., \text{supp}$(\sigma-\sigma_0)$, 
as well as $d_x \approx d_{i}$ for such $i$, and is small otherwise.
Therefore, this decouples the effect of Green's functions and their gradients \mm{in the near field or scattered 
data} with the help of monopole and dipole probing functions, thanks to the mutually almost orthogonality property.  
In order to maximize such a decoupling effect, different types of boundary influxes and probing directions 
\mm{are also analysed.}  
The above properties and strategies for decoupling will be addressed in further detail in the rest of the work.

\smallskip
Under the settings above, two index functions in \eqref{eq:dsm} give rise to our new Direct Sampling Method: 

\smallskip
Given the measurement data 
$u-u_{0}$ on $\partial\Omega$, and a set of discrete sampling points $x\in\Omega$, 

(i) evaluate $I_{\text{mo}}$ to recover the potential inclusions, i.e., \text{supp}$(V-V_0)$;  

(ii) evaluate $I_{\text{di}}$ to recover the conductivity inclusions, i.e., \text{supp}$(\sigma-\sigma_0)$.  

\section{Fundamental property}
\label{sec_fundamental}
In this section, \mm{we aim to verify the fundamental property introduced in section \ref{sec_principle} 
for some typical cases} that we encounter in real applications.
In particular, we intend to derive an approximation of the measurements as a combination of the Green's functions of the homogeneous medium and their gradients \mm{when $\sigma$ is either smooth or piecewise constant.}  

\mm{Associated with the model \eqref{eqn_prostso}, the incident field $u_0$ from the homogeneous background 
satisfies}
\begin{equation}
\label{eqn_prostsh}
\begin{cases}
  -\nabla \cdot(\sigma _0 \nabla u_0) + V_0u_0 = 0 \quad &\text{in} \quad \Omega\,, \\
  \frac{\partial u_0}{\partial \nu} = f \quad &\text{on} \quad \partial \Omega\,.
\end{cases}
\end{equation}
Combining \mm{the systems \eqref{eqn_prostso} and \eqref{eqn_prostsh}, we readily see}
\begin{equation}
  \label{eqn_diffrep}
  \begin{cases}
  -\Delta (u-u_0)+\frac{V_0}{\sigma_0}(u-u_0) = \frac{1}{\sigma_0}[\nabla\cdot((\sigma-\sigma _0)\nabla u) -(V-V_0)u] \quad &\text{in} \quad \Omega\,, \\
  \frac{\partial (u-u_0)}{\partial \nu} = 0  &\text{on} \quad \partial \Omega\,.
  \end{cases}
\end{equation}

\mm{If $V_0 \neq 0$, we consider the Green's function $G_x$ for $x\in \Omega$ satisfying} 
\begin{equation}
  \label{eqn_greeno}
    -\Delta G_x + \frac{V_0}{\sigma _0}G_x = \delta _x \quad \text{in} \quad \Omega\,, \qquad \frac{\partial G_x}{\partial \nu} = 0 \quad \text{on} \quad \partial \Omega \,.
\end{equation}
Then the difference $u-u_0$ 
can be represented by 
\begin{equation}
\label{eqn_greenrep}
  (u-u_0) (x) = \frac{1}{\sigma_0}\int _{\Omega} \bigg[\nabla _y\cdot((\sigma-\sigma _0)\nabla _y u) -(V-V_0)u\bigg]G_xdy\,.
\end{equation}

On the other hand, if $V_0 = 0$, we consider the following Green's function $G_x$ for $x\in \Omega$ instead:
\begin{equation}
\label{eqn_greenrepmu0}
-\Delta G_x = \delta_x \quad \text{in} \quad \Omega\,, \qquad \frac{\partial G_x}{\partial \nu} = -\frac{1}{|\partial \Omega|} \quad \text{on} \quad \partial \Omega \,, \qquad \int_{\partial \Omega}G_x ds = 0\,.
\end{equation}
\mm{Then we can obtain a similar representation to \eqref{eqn_greenrep}.}

From now on, we shall consider only the following two typical cases:  
\mm{either $\sigma$ is smooth or piecewise constant.}  
\mm{First for the case when $\sigma \in C^{1}(\Omega)$, 
by writing $D:= \text{supp}(\sigma - \sigma_0) \bigcup \text{supp}(V - V_0) \Subset \Omega$, 
we can readily derive from \eqref{eqn_greenrep} by the divergence theorem that 
}
\begin{equation}
  \label{eqn_greenrepsm}
  \begin{split}
  (u-u_0) (x) =& \,\frac{1}{\sigma_0} \bigg[\int_{\partial \Omega} (\sigma -\sigma _0)G_x \frac{\partial u}{\partial \nu} ds(y) - \int_{\Omega}(\sigma - \sigma_0)\nabla_y u \cdot \nabla _y G_x dy - \int_{\Omega}(V-V_0)uG_xdy \bigg]\\ 
   =& -\frac{1}{\sigma_0}\bigg[\int_{\Omega}(\sigma - \sigma_0)\nabla_y u \cdot \nabla _y G_x dy + \int_{\Omega}(V-V_0)uG_xdy \bigg]\,.
   \end{split}
\end{equation}

\mm{Next, we consider the case when $\sigma$ is piecewise constant.}
We assume that $D = \overline{\cup_{i=1}^m\Omega _i} $, where $\Omega_i$ are open subsets of $\Omega$ with smooth boundary such that $\Omega_i \bigcap \Omega_j = \emptyset$, and that $\sigma= \sigma_i$ in $\Omega _i$ 
for some constant $\sigma_i$. 
\mm{And we further write $\Omega_0 = \bar{ \Omega }\setminus D$ for simplicity. 
Then for all $\phi \in H^1({\Omega})$, we derive from \eqref{eqn_prostso} that 
}
\begin{equation}
\begin{split}
\label{eqn_variational}
    0 &= \sum_{i=0}^m \bigg(\int_{\Omega _i} \sigma_i\nabla u \cdot \nabla \phi dy\bigg) - \int_{\partial \Omega} \sigma_0 f \phi ds(y) + \int_{\Omega}V u \phi dy \\
    &= \sum_{i=0}^m \bigg[\int_{\Omega_i} \bigg(- \sigma_i \Delta u + V u \bigg)\phi dy\bigg] + \sum_{i=1}^m \bigg[\int_{\partial \Omega_i} \bigg(\sigma_i \frac{\partial u^-}{\partial \nu} - \sigma_0\frac{\partial u^+}{\partial \nu}\bigg)\phi ds(y)\bigg]\,.\\
  \end{split}
\end{equation}
\mm{Noticing that the normal derivative of $u$ has a jump across $\partial \Omega_i$, we get for $v := \sigma u$ 
from \eqref{eqn_variational} that 
}
\begin{equation*}
   - \Delta v + \frac{V_0}{\sigma_0} v = (\frac{\sigma}{\sigma_0}V_0 - V)u \quad \text{in} \quad\Omega\setminus(\cup_{i=1}^m\partial \Omega_i )\,; \quad \frac{\partial v^+}{\partial \nu}|_{\partial \Omega_i} = \frac{\partial v^-}{\partial \nu}|_{\partial \Omega_i}
  ~~\text{on} ~~\partial \Omega_i\,, 
  \quad \frac{\partial v}{\partial \nu} = \frac{f}{\sigma_0} ~~\text{on} ~~\partial \Omega\,, 
\end{equation*}
\mm{where we have chosen the normal vector to point towards $\Omega_0$ on each 
$\partial \Omega_i$, and will write the jump of any function $w$ across the boundary $\partial \Omega_i$ 
as $[w] := w^+-w^-$. 
The above equation readily implies the equation for $\gamma: = \sigma u - \sigma_0 u_0$
}
\begin{equation}
\begin{cases}
      -\Delta \gamma +\frac{V_0}{\sigma_0}\gamma = -(V-V_0)u + (\sigma -\sigma_0)\frac{V_0}{\sigma_0} u &\text{in} \quad\Omega\setminus(\cup_{i=1}^m\partial \Omega_i )\,, 
    \\ \frac{\partial \gamma}{\partial \nu} = 0 \quad \text{on} \quad \partial \Omega\,,
    \qquad \frac{\partial \gamma^+}{\partial \nu}|_{\partial \Omega_i} = \frac{\partial \gamma^-}{\partial \nu}|_{\partial \Omega_i} &\text{on} \quad \partial \Omega_i\,.
  \end{cases}
\end{equation}
\mm{For any $x\in \Omega_0$, we can easily write} 
\begin{equation}
\begin{split}
\label{pc1}
    &\sum_{i=1}^m \int_{\partial \Omega_i} [\gamma] \frac{\partial G_x}{\partial \nu} ds(y)\\ 
    =& \sum_{i=1}^m \int_{\partial \Omega_i}\bigg(\gamma^+\frac{\partial G_x}{\partial \nu} - \frac{\partial \gamma^+ }{\partial \nu} G_x \bigg) ds(y)- \sum_{i=1}^m\int_{\partial \Omega_i}\bigg(\gamma^-\frac{\partial G_x}{\partial \nu} - \frac{\partial \gamma^- }{\partial \nu} G_x \bigg)ds(y)\\
     =& \bigg[ \sum_{i=1}^m \int_{\partial \Omega_i}\bigg(\gamma^+ \frac{\partial G_x}{\partial \nu} - \frac{\partial \gamma^+ }{\partial \nu} G_x\bigg) ds(y)- \int_{\partial \Omega}\bigg(\gamma \frac{\partial G_x}{\partial \nu} - \frac{\partial \gamma }{\partial \nu} G_x \bigg) ds(y)\bigg]\\
  &- \bigg[ \sum_{i=1}^m \int_{\partial \Omega_i}\bigg(\gamma^- \frac{\partial G_x}{\partial \nu} - \frac{\partial \gamma^- }{\partial \nu} G_x \bigg)ds(y) \bigg]\,.
    \end{split}
\end{equation}
Applying the Green's formula in $\Omega_0$ to the first part of the above difference, we obtain
\begin{equation}
  \begin{split}
  \label{pc2}
    &\sum_{i=1}^m \bigg[\int_{\partial \Omega_i}\bigg(\gamma^+\frac{\partial G_x}{\partial \nu} - \frac{\partial \gamma^+ }{\partial \nu} G_x \bigg)ds(y)\bigg] - \int_{\partial \Omega}\bigg(\gamma \frac{\partial G_x}{\partial \nu} - \frac{\partial \gamma}{\partial \nu} G_x\bigg)ds(y)\\ =& \int_{\Omega_0} \bigg(\Delta \gamma G_x - \gamma \Delta G_x\bigg)dy 
    = \int_{\Omega_0} \bigg[(V-V_0)uG_x + \sigma_0(u-u_0) \bigg]dy\,.
  \end{split}
\end{equation}
Meanwhile, for the second part of the difference in \eqref{pc1}, we notice the following for each $\Omega_i$:
\begin{equation}
  \begin{split}
  \label{pc3}
    &\int_{\partial \Omega_i}\bigg(\gamma^-\frac{\partial G_x}{\partial \nu} - \frac{\partial \gamma^- }{\partial \nu} G_x\bigg) ds(y)
    = \int_{\Omega_i}\bigg[ -(V-V_0)uG_x + (\sigma_i - \sigma_0)u \Delta G_x \bigg]dy \\
    =& -\int_{\Omega_i}\bigg[(V-V_0)uG_x+(\sigma_i-\sigma_0)\nabla u\cdot \nabla G_x \bigg] dy + \int_{\partial \Omega_i}(\sigma_i-\sigma_0)u^-\frac{\partial G_x}{\partial \nu} ds(y)\,.
  \end{split}
\end{equation}
Combining \eqref{pc1}-\eqref{pc3}, we come to the difference of the potentials
\begin{align}
 \label{eqn_greenreppc}
    &(u-u_0)(x) \\
    =&-\frac{1}{\sigma_0} \bigg\{ \int_{\Omega}(V-V_0)uG_xdy
    + \sum_{i=1}^m\bigg[\int_{\Omega_i}(\sigma_i - \sigma_0)\nabla u \cdot \nabla G_xdy + \int_{\partial \Omega_i}([\gamma] - (\sigma_i - \sigma_0)u^-) \frac{\partial G_x}{\partial \nu} ds(y)\bigg]\bigg\}\,.\notag
\end{align}
\mm{Now using some appropriate numerical quadrature rule, 
we can easily see from the expressions \eqref{eqn_greenrepsm} and \eqref{eqn_greenreppc} that 
the boundary data or the scattered field 
can be approximated by 
}
\begin{equation}
\label{eqn_solreps}
(u-u_{0})(x) \approx \sum_{j=1}^{n} c_jG_{q_j}(x) + \sum_{i=1}^{m} a_i d_i\cdot \nabla G_{p_i}(x)\,, 
\quad x\in \partial\Omega
\end{equation}
for some coefficients $a_i\in\mathbb{C}$, $c_j\in \mathbb{C}$, 
$d_i \in \mathbb{S}^{d-1}$, 
\mm{and some quadrature points} $p_i \in \text{supp}(\sigma - \sigma_0)$ and $q_j \in \text{supp}(V-V_0)$.
We have therefore verified the fundamental property introduced in section \ref{sec_principle}. 

\section{Probing and index functions}
\label{sec_probing}
\subsection{Monopole and dipole probing functions}

In order to accurately locate the respective medium inhomogeneities $ \text{supp}(\sigma - \sigma_0)$ and 
$\text{supp}(V-V_0)$, 
\mm{we are expected} to decouple the effects of the Green's function $G_x$ and $\nabla G_x$ in \eqref{eqn_solreps}. 
For this purpose, we define two groups of probing functions, $ \{ \zeta_{x} \}_{x \in \Omega}$ representing a family of monopole probing functions from sources $x \in \Omega$, and $ \{ \eta_{x,d} \}_{x \in \Omega, d \in \mathbb{S}^{d-1}}$ representing a family of dipole probing functions from sources $x \in \Omega$ and dipole directions $d \in \mathbb{S}^{d-1}$.

We first introduce the family of monopole probing functions $ \{ \zeta_{x} \}_{x \in \Omega}$. 
\mm{For a point $x\in \Omega$, we consider a monopole potential $v_x$ satisfying}
\begin{equation}
  \begin{cases}
    -\Delta v_{x}+\frac{V_0}{\sigma_0}v_{x} = \delta _x &\text{in} \quad \Omega\,, \\
    v_{x} = 0 &\text{on} \quad \partial \Omega\,.
  \end{cases} 
\end{equation}
We then define $\zeta_x$ as the boundary flux of $v_x$
\begin{equation}
  \label{eqn_defzeta}
 \zeta_{x} := -\frac{\partial v_x}{\partial \nu} \quad \text{on} \quad \partial \Omega\,.
\end{equation}
\mm{To avoid the approximation of a delta measure in computing $\zeta_{x}$, we may evaluate $v_x$ 
using its equivalent expression $v_x = v^{(1)}_x - v^{(2)}_x$, where $v^{(1)}_{x}$ is the fundamental solution 
in the whole space $\mathbb{R}^d$ with any appropriate boundary condition, namely 
}
\begin{equation}
  \label{eqn_defmonoknown1}
  -\Delta v^{(1)}_{x}+\frac{V_0}{\sigma_0}v^{(1)}_x = \delta_x \quad \text{in} \quad\mathbb{R}^d\,, 
\end{equation}
\mm{while $v^{(2)}_x$ solves} 
\begin{equation}
\label{eqn_defmonoknown2}   
-\Delta v^{(2)}_x + \frac{V_0}{\sigma_0}v^{(2)}_x = 0 \quad \text{in} \quad\Omega \,, \quad v^{(2)}_x = v^{(1)}_x \quad \text{on} \quad\partial \Omega\,.
\end{equation}

Next we define another family of dipole probing functions $ \{ \eta_{x,d} \}_{x \in \Omega, d \in \mathbb{S}^{d-1}}$.  Given $x\in \Omega$ and $d \in \mathbb{S}^{d-1}$, \mm{we consider the dipole potential $w_{x,d}$ satisfying}
\begin{equation}
\label{eqn_defwxd}
  \begin{cases}
    -\Delta w_{x,d}+\frac{V_0}{\sigma_0}w_{x,d} = -d\cdot\nabla \delta _x &\text{in} \quad \Omega\,, \\
    w_{x,d} = 0 &\text{on} \quad \partial \Omega\,, 
  \end{cases}
\end{equation}
then we define $\eta_{x,d}$ as the boundary flux
\begin{equation}
  \label{eqn_defeta}
  \eta_{x,d} := -\frac{\partial w_{x,d}}{\partial \nu}\quad \text{on} \quad \partial \Omega\,.
\end{equation}
\mm{Similarly, to avoid the approximation of a delta measure in computing $\eta_{x,d}$, we may evaluate $w_{x,d}$ 
using its equivalent expression $w_{x,d} = w_{x,d}^{(1)} - w_{x,d}^{(2)}$, where $w_{x,d}^{(1)}$ is defined 
as \eqref{eqn_defmonoknown1} with the right-hand side replaced by $-d\cdot \nabla \delta_x$ while 
$w_{x,d}^{(2)}$ is defined as \eqref{eqn_defmonoknown2}. 
}

\subsection{Monopole and dipole index functions}
\mm{We are now ready to define two critical index functions that give rise to our new direct sampling method.  
For this purpose, 
for a given $\gamma \geq 0$ and an auxiliary choice of $l \geq 0$, 
we introduce a Sobolev duality product
}
\begin{equation}
\label{def_innerpr}
  \langle f,g \rangle_{H^{\gamma}(\partial \Omega)}:= \int_{\partial \Omega}(-\Delta _{\partial \Omega})^{\gamma} \overline{f(x)} g(x) d s(x) \quad \forall \,f\in H^{2\gamma + l }(\partial \Omega), \quad g\in H^{-l}(\partial \Omega)\,.
\end{equation}
\mm{We notice that for $f$, $g\in H^\gamma(\partial \Omega)$,  
the above duality product is the standard definition of a $\gamma$-semi-inner product on $H^\gamma(\partial \Omega)$. 
However, 
the argument $g$ in \eqref{def_innerpr} will play the role of the noisy measurement from the forward problem, which exists generally only in $H^{-l}(\partial \Omega)$ for some $l\geq 0$.  
For simplicity, we will often write $\langle \cdot , \cdot \rangle_{H^{\gamma}}$ instead of $\langle \cdot, \cdot \rangle_{H^{\gamma}(\partial \Omega)}$, and 
use $|\cdot|_{H^{\gamma}(\partial \Omega)}$ as the $H^\gamma$ semi-norm induced by the duality product in \eqref{def_innerpr}.
$\gamma$ is often called a Sobolev scale.
}

\mm{We are now ready to introduce our two index functions.  
First, for any $x \in \Omega$, $d \in \mathbb{S}^{d-1}$, we know 
$ \zeta_x$, $\eta_{x,d} \in H^{2 \gamma - l}(\partial \Omega)$ for any $\gamma$, $l \geq 0$. 
Then 
corresponding to the monopole probing functions in \eqref{eqn_defzeta} and the dipole probing functions \eqref{eqn_defeta}, 
we define the index functions as follows:
}  
\begin{eqnarray}
  I_{\text{mo}}(x) &:=& \frac{\langle\zeta_x, u_s \rangle_{H^{\gamma_{\text{mo}}}(\partial \Omega)}}{|\zeta_x|_{H^{\gamma_{\text{mo}}}(\partial \Omega)}^{n_1}\cdot|G_x|_{H^{\gamma_{\text{mo}}}(\partial \Omega)}^{n_2}}  \,,  \label{eqn_defindexmo} \\
  I_{\text{di}}(x,d_x) &:=& \frac{\langle\eta_{x,d_x}, u_s \rangle_{H^{\gamma_{\text{di}}}(\partial \Omega)}}{|\eta_{x,d_x}|_{H^{\gamma_{\text{di}}}(\partial \Omega)}^{m_1}\cdot| d_x \cdot \nabla G_x|_{H^{\gamma_{\text{di}}}(\partial \Omega)}^{m_2}} \,,   \label{eqn_defindexdi} 
\end{eqnarray} 
\mm{under appropriate choices of two Sobelov scales $\gamma_{\text{mo}}$ and $\gamma_{\text{di}}$ and 
the coefficients $n_i$ and $m_i$. 
}

\mm{Using \eqref{eqn_solreps}, we have the approximations}
\begin{eqnarray*}
I_{\text{mo}}(x) &\approx & \sum_{j=1}^n c_j \frac{\langle\zeta_{x}, G_{q_j}\rangle_{H^{\gamma_{\text{mo}}}}}{|\zeta_{x}|_{H^{\gamma_{\text{mo}}}}^{n_1}\cdot|G_x|_{H^{\gamma_{\text{mo}}}}^{n_2}} + \sum_{i=1}^m a_i \frac{\langle\zeta_{x}, d_i\cdot \nabla G_{p_i}\rangle_{H^{\gamma_{\text{mo}}}}}{|\zeta_{x}|_{H^{\gamma_{\text{mo}}}}^{n_1}\cdot|G_x|_{H^{\gamma_{\text{mo}}}}^{n_2}} \\
&=&  \sum_{j=1}^n c_j \, K_1(x,q_j) + \sum_{i=1}^m a_i \, K_{2, d_i } (x,p_i) \,, \\
I_{\text{di}}(x,d_x) &\approx& \sum_{j=1}^n c_j \frac{\langle\eta_{x,d_x}, G_{q_j}\rangle_{H^{\gamma_{\text{di}}}}}{|\eta_{x,d_x}|_{H^{\gamma_{\text{di}}}}^{m_1}\cdot| d_x \cdot \nabla G_x|_{H^{\gamma_{\text{di}}}}^{m_2}} 
+
 \sum_{i=1}^m a_i \frac{\langle\eta_{x,d_x}, d_i\cdot \nabla G_{p_i}\rangle_{H^{\gamma_{\text{di}}}}}{|\eta_{x,d_x}|_{H^{\gamma_{\text{di}}}}^{m_1}\cdot| d_x \cdot \nabla G_x|_{H^{\gamma_{\text{di}}}}^{m_2}} \\
&=& \sum_{j=1}^n c_jK_{3,d_x}(x,q_j) + \sum_{i=1}^m a_i K_{4,d_x, d_i}(x,p_i) \,,
\end{eqnarray*}
\mm{where the kernels $K_s$ for $s = 1,2,3,4$ are now respectively given by}
\begin{align}
 K_1(x,z) &= \frac{\langle\zeta_{x}, G_z\rangle_{H^{\gamma_{\text{mo}}}}}{|\zeta_{x}|_{H^{\gamma_{\text{mo}}}}^{n_1}\cdot|G_x|_{H^{\gamma_{\text{mo}}}}^{n_2}} \,, 
 \quad &K_{2,d_z}(x,z) &= \frac{\langle\zeta_{x}, d_z\cdot \nabla G_z\rangle_{H^{\gamma_{\text{mo}}}}}{|\zeta_{x}|_{H^{\gamma_{\text{mo}}}}^{n_1}\cdot|G_x|_{H^{\gamma_{\text{mo}}}}^{n_2}} \, ; \label{kernal_mo} 
  \\ K_{3,d_x}(x,z) &= \frac{\langle\eta_{x,d_x}, G_z\rangle_{H^{\gamma_{\text{di}}}}}{|\eta_{x,d_x}|_{H^{\gamma_{\text{di}}}}^{m_1}\cdot| d_x\cdot \nabla G_x|_{H^{\gamma_{\text{di}}}}^{m_2}} \,,
  \quad &K_{4,d_x,d_z}(x,z) &=\frac{\langle\eta_{x,d_x}, d_z\cdot \nabla G_z\rangle_{H^{\gamma_{\text{di}}}}}{|\eta_{x,d_x}|_{H^{\gamma_{\text{di}}}}^{m_1}\cdot| d_x\cdot \nabla G_x|_{H^{\gamma_{\text{di}}}}^{m_2}} \,.   \label{kernal_di}  
\end{align}
Therefore, if we have the mutually almost orthogonality property between the two families of probing functions and the fundamental solution with its gradient respectively under the aforementioned duality product, 
we shall be able to decouple the effects coming from monopoles and dipoles, and 
reconstruct inhomogeneous inclusions as well as recognize their types with one or two pair(s) of Cauchy data. 
In section \ref{sec_explicit}, we will verify these desired properties of probing functions under our special choice of the duality product in some typical sampling domains. 

\mm{We end this subsection with} two \mm{helpful remarks:}
\begin{enumerate}
\item
\mm{In order to numerically evaluate our index functions efficiently from the measurement data, we need only to compute 
the Sobolev duality product approximately after discretization. 
The approximations of the $H^{\gamma}$ norm and pointwise values of probing functions can be all computed off-line. 
The entire algorithm does not involve any iterative procedure or matrix inversion.  
}

\item
\mm{We would like to comment on the intuition of what the surface Laplacian in \eqref{def_innerpr} does.  
Considering the fact that when $x$ approaches the boundary, one may represent the Laplacian in terms of the surface Laplacian operator (up to the boundary)
}
\begin{equation}
\label{boundary_laplacian}
   \Delta_{\partial \Omega} u(x) = - \Delta u(x) + \frac{V_0}{\sigma _0} u(x) + \Delta_{\partial \Omega} u(x) = - \frac{\partial^2 u}{\partial \nu^2}(x) - (d-1)H(x)\frac{\partial u }{\partial \nu} (x) + \frac{V_0}{\sigma _0} u(x) \, ,
\end{equation}
\mm{where $H(x)$ represents the mean curvature of $\partial \Omega$ embedded in $\mathbb{R}^d$ at the point $x$ and the normal derivative is taken outward from the inside. Therefore, we may expect that, by choosing a larger value of Sobolev scale 
$\gamma$, 
we are essentially taking a higher order normal derivative of the boundary measurement in the distributional sense, 
i.e., a higher order flux of the measurement at the boundary.
Hence, taking a bigger $\gamma$ in the duality product amounts to comparing the higher order details of probing functions along the boundary (either in the tangential or normal direction) with that of monopole/dipole functions in the measurement. This can improve 
the reconstruction results; see our numerical studies in 
Example 1 of section \ref{sec_numerical}. 
}


\end{enumerate}

\subsection{Alternative characterization of index functions}
\label{subsec_alternative}
\mm{In order to simplify the computation and obtain a better understanding of the index functions 
\eqref{eqn_defindexmo} and \eqref{eqn_defindexdi}, as well as to make an optimal choice of the probing direction $d_x$ 
there 
we now present an alternative characterization of the index functions. 
For this purpose, let us consider $\phi$ to be an auxiliary function that solves 
}
\begin{equation}
  \label{eqn_defchara}
  \begin{cases}
    -\Delta \phi + \frac{V_0}{\sigma_0}\phi = 0 &\text{in} \quad \Omega\,,
    \\ \phi = (- \Delta_{\partial \Omega})^{\gamma}(u-u_0) & \text{on} \quad\partial \Omega\,.
  \end{cases}
\end{equation}
\mm{where the boundary condition is understood in the distributional sense. 
Using the definitions \eqref{eqn_defeta} and \eqref{def_innerpr}, we can easily observe that 
}
\begin{equation}
  \label{eqn_charadiG}
  \begin{split}
\langle\eta_{x,d},u_s\rangle_{H^{\gamma}(\partial \Omega)} =& \int_{\partial \Omega} (-\Delta_{\partial \Omega})^{\gamma} (u-u_0) \, \overline{\eta_{x,d}} dy 
  = -\int_{\Omega}\bigg( \phi \overline{\Delta w_{x,d}} + \nabla \phi \cdot \overline{\nabla w_{x,d}} \bigg)dy\\
  =& \int_{\Omega} \bigg(\frac{V_0}{\sigma _0}\phi \overline{w_{x,d}}-\phi \overline{\Delta w_{x,d}} \bigg)dy
  =\,\, d\cdot \nabla \phi(x)\,.
    \end{split}
\end{equation}
Similarly, from definitions \eqref{eqn_defzeta} and \eqref{def_innerpr}, we readily obtain
\begin{equation}
  \label{eqn_charamoG}
  \begin{split}
  \langle\zeta_{x,d},u_s\rangle_{H^{\gamma}(\partial \Omega)} =& \int_{\partial \Omega} (-\Delta_{\partial \Omega})^{\gamma} (u-u_0) \, \overline{ \zeta_{x}} dy   = \phi(x)\,.
      \end{split}
\end{equation}
With the help of the above expressions, we can therefore rewrite \eqref{eqn_defindexmo} and \eqref{eqn_defindexdi} as
\begin{equation}
  \label{eqn_moalternatedef}
  I_{\text{mo}}(x) = \frac{\phi(x)}{|\zeta_{x} |_{H^{\gamma}}^{n_1}\cdot|G_x|_{H^{\gamma}}^{n_2}}\,, \qquad   I_{\text{di}}(x,d_x) = \frac{d_x\cdot \nabla \phi(x)}{|\eta_{x,d_x} |_{H^{\gamma}}^{m_1}\cdot|d_x\cdot \nabla G_x|_{H^{\gamma}}^{m_2}}\,.
\end{equation}
The above understanding of the index functions helps in two folds:
\begin{enumerate}
\item
First, this provides us another way to quickly compute index functions. In particular, given that $\partial \Omega$ is smooth enough, we could quickly evaluate the surface Laplacian. It then remains to numerically solve a Dirichlet boundary value problem for $\phi$ by 
any appropriate numerical method.
\item
This expression helps us obtain an optimal choice of the probing direction $d_x$ at each point $x\in \Omega$. In fact, based on the expression \eqref{eqn_moalternatedef}, we can see that the magnitude of $I_{\text{di}}(x,d_x) $ can be maximized by choosing $d_x$ parallel to $\nabla \phi(x)$, and minimized when we choose a $d_x$ that is orthogonal to $\nabla \phi(x)$. Therefore, in order to locate $\text{supp}(\sigma - \sigma_0)$, we may therefore maximize $I_{\text{di}}(x,d_x) $ by choosing 
\begin{equation}
\label{optimal_optimal}
d_x = \frac{{\nabla \phi(x)}}{|\nabla \phi(x)|} \,.
\end{equation}
\mm{This serves as a guide for an optimal probing direction.}
\end{enumerate}

\section{Explicit expressions of probing functions and index functions in some special domains}
\label{sec_explicit}
In this section, we aim at obtaining some explicit expressions of our choices of probing functions in some special domains for more efficient numerical computation.  With the same technique, We can also obtain explicit expressions of kernels $K_i$ introduced in \eqref{kernal_mo} and \eqref{kernal_di} in those cases, which help us understand more precisely the behaviour of those kernels, and verify the mutually almost orthogonality properties.

\mm{For the notational sake, we shall write $k^2 := {V_0}/{\sigma_0}$ from now on.
The Poincare-Steklov operator plays an essential role in our subsequent analysis.  
We define the Neumann-to-Dirichlet map (NtD) as $\Lambda f = g$, where $f$ and $g$ satisfy the equations
}
\begin{equation}
\label{eqn_defNtD}
\begin{cases}
    -\Delta \Phi + k^2 \Phi = 0 \quad &\text{in} \quad  \Omega \,,\\
    \frac{\partial}{\partial \nu} \Phi = {f} \quad &\text{on} \quad\partial \Omega\,, \\
    \Phi = g \quad &\text{on} \quad \partial \Omega \,.
\end{cases}
\end{equation}
We recall that $\Lambda: H^{-\frac{1}{2}} (\partial \Omega) \rightarrow H^{\frac{1}{2}} (\partial \Omega)$ is a compact self-adjoint operator when we restrict ourselves to $L^2(\partial \Omega)$. Therefore, there exists a complete orthonormal basis consists of eigenfunctions of $\Lambda$.
We notice that, in some special cases, this set of eigenfunctions coincides with the set of eigenfunctions of the surface Laplacian $\Delta_{\partial \Omega}$.  This helps us to write both probing functions and the $H^{\gamma}(\partial \Omega)$ semi-inner product defined in \eqref{def_innerpr} explicitly via Fourier coefficients with respect to the same orthonormal basis. In this section, we will focus on one such case, that is when $\partial \Omega = R \mathbb{S}^{d-1}$ for some $R > 0$ and $d \geq 2$, 
\mm{which is a typical geometric shape used in many applications.}  

\mm{We like to point out that,} although the two sets of eigenfunctions differ in general, they are comparable to each other based on the following observation: if we denote  $\| \xi\|^2_{g(x)} := \langle \xi, \,g^{-1}(x) \xi\rangle$, the dual norm of $\xi$ under the metric $g(x)$ on the surface, then the principle symbol of $\Delta_{\partial \Omega}$ is 
$\| \xi \|_{g(x)}^2$, while that of $\Lambda$ is $\| \xi \|_{g(x)}^{-1}$ (Proposition 8.53, \cite{DtN_Gunther}). 
With this, via an application of the generalized Weyl's law, we can obtain a precise comparison of the pointwise asymptotic average squared density between the two sets of eigenfunctions.  In fact, one readily checks that the volume of the variety coming from the two Hamiltonians $ \{ \xi :  \| \xi \|_{g(x)}^2 = 1 \} $ and $ \{ \xi:  \| \xi \|_{g(x)}^{-1} = 1 \} $ are in fact the same, and the generalized Weyl's law will therefore render us that the two sets of eigenfunctions have the same pointwise asymptotic average squared density 
\mm{in some sense mathematically.
We skip} the details of this argument for the sake of exposition, and focus only 
on the case $\partial \Omega = R \mathbb{S}^{d-1}$ for some $R > 0$, when the two sets of eigenfunctions coincide.

\subsection{Circular domains}
\label{subsec_circular}
Now let us consider the special case when the domain $\Omega = B_R \subset \mm{\mathbb {R}^2}$ is a disk 
with radius $R > 0$ centered at the origin. We consider the following Poincare-Steklov eigenvalue problem:
\begin{equation}
\label{eqn_defNtDcir}
\begin{cases}
    -\Delta \varphi_n + k^2 \varphi_n = 0 \quad &\text{in} \quad B_R\,, \\
    \frac{\partial}{\partial \nu} \varphi_n = \frac{1}{\lambda_n}f_n \quad &\text{on} \quad \partial B_R\,, \\
        \varphi_n = f_n \quad &\text{on} \quad \partial B_R\,.
\end{cases}
\end{equation}
Writing $I_n$ as the modified Bessel function of the first kind of order $n$, we readily obtain, via a separation of variables, that eigenfunctions of $\Lambda$ and their associated eigenvalues are given by 
\begin{equation}
    \label{def_modbess}
    \varphi_n = \begin{cases}
        \frac{I_n(kr)}{I_n(kR)}e^{in\theta}\,, \quad & k^2\neq 0\,;  \\
                \frac{r^{|n|}}{R^{|n|}}e^{in\theta}\,, \quad & k^2 = 0\,; 
    \end{cases}
    \\
    \qquad \lambda_n = \begin{cases}
        \frac{I_n(kR)}{kI_n'(kR)}\,, \quad & k^2\neq 0\,;  \\
        \frac{R}{|n|}\,, \quad   & k^2 = 0\quad (n\neq 0)\,.
\\
    \end{cases}
\end{equation}
\mm{From these explicit expressions, one can readily find for $k^2 \neq 0$ and $k=0$ that
}
\begin{eqnarray}\label{eq:In}
    \nabla \varphi_n 
    =&\frac{e^{in\theta}}{I_n(kR)}\begin{pmatrix}
        \cos(\theta) & -\sin(\theta) \\
        \sin(\theta) & \cos(\theta)
    \end{pmatrix}  
    \begin{pmatrix}kI_n'(kr) \\\frac{inI_n(kr)}{r} \end{pmatrix}\,
    \quad &\mbox{for} ~~k\ne 0\,, \\
\label{eq:Iin}
    \nabla\varphi_n =&\frac{r^{|n|-1}}{R^{|n|}}e^{in\theta}\begin{pmatrix}
        \cos(\theta) & -\sin(\theta) \\
        \sin(\theta) & \cos(\theta)
    \end{pmatrix} \begin{pmatrix} |n| \\ in \end{pmatrix} \quad &\mbox{for} ~~k=0\,. 
\end{eqnarray}
\mm{Recalling the definition of the dipole probing function in \eqref{eqn_defeta}, 
we obtain their Fourier coefficients}
\begin{equation}
\label{eqn_cirFouco1}
\begin{split}    
    R\int_{\partial B_R}e^{in\theta_y}\eta_{x,d}d\theta_y &= -\int_{\partial B_R}\varphi _n \eta_{x,d} \, ds(y)
    = \int_{\partial B_R}\bigg( w_{x,d}\frac{\partial \varphi _n}{\partial \nu}-\varphi _n\frac{\partial w_{x,d}}{\partial \nu}  \bigg)ds(y)  \\ 
    &= \int_{B_R}\bigg(  k^2  w_{x,d} \varphi_n -\Delta w_{x,d} \varphi_n\bigg) dy
    =\,\,d\cdot \nabla \varphi_n(x)\,.
    \end{split}
\end{equation}
Similarly, from the definition of the monopole probing function in \eqref{eqn_defzeta}, we derive 
\begin{equation}
\label{eqn_cirFouco2}
    \begin{split}    
    R\int_{\partial B_R}e^{in\theta_y}\zeta_{x}d\theta_y 
= \int_{\partial B_R}\bigg( v_{x}\frac{\partial \varphi_n}{\partial \nu} -\varphi_n\frac{\partial v_{x}}{\partial \nu}\bigg)ds(y)  
    =  \varphi_n(x)\,.
    \end{split}
\end{equation}
On the other hand, we can deduce from definitions \eqref{eqn_greeno} and \eqref{eqn_defNtDcir} that
\begin{equation}
\label{eqn_cirFouco3}
\begin{split}
    R\int_{\partial B_R}e^{in\theta_y}G_{x}d\theta_y 
    = \lambda_n \int_{\partial B_R}\bigg( G_x \frac{\partial \varphi_n}{\partial \nu}  - \varphi_n \frac{\partial{G_x}}{\partial \nu} \bigg)ds(y) 
    = \lambda_n \, \varphi_n (x)\,.
\end{split}
\end{equation}
Differentiating \eqref{eqn_cirFouco3} with respect to $x$, and considering the symmetry of the Green's function $G_x$ in \eqref{eqn_greeno}, i.e., $\nabla G_x = \nabla_x G_x$, we obtain 
\begin{equation}
\label{eqn_cirFouco4}
\begin{split}
 R\int_{\partial B_R}e^{in\theta_y} d\cdot\nabla G_{x}d\theta_y   = R\int_{\partial B_R}e^{in\theta_y} d\cdot\nabla_x G_{x}d\theta_y   = \lambda_n \, d\cdot\nabla \varphi_n (x)\,.
\end{split}
\end{equation}    

Now let us recall the definition of the duality product in \eqref{def_innerpr}.   
\mm{When $\Omega = B_R$, with 
$\hat{f}(n): = \int_{\partial B_R} {f(\theta)}e^{-in\theta}d\theta$, 
one may readily check that $\Delta_{\partial \Omega} e^{i n \theta} = n^2 e^{i n \theta} $, and therefore 
}
\begin{equation}
    \label{eqn_fourexpancir}
    \begin{split}
    \langle f,g\rangle _{H^{\gamma}(\partial B_R)} = \sum_{n=-\infty}^{\infty}\frac{R |n|^{2\gamma}}{2\pi} \overline{\hat{f}(n)}\hat{g}(n)\, ,
        \end{split}
    \end{equation}
\mm{Using \eqref{eqn_cirFouco1}-\eqref{eqn_fourexpancir}, 
we can obtain the explicit expressions of the duality products and $H^{\gamma}$ semi-norms
}
\begin{eqnarray}
\label{eqn_ipcir1}
    \langle\eta_{x_1,d_1}, d_2 \cdot \nabla G_{x_2} \rangle_{H^{\gamma}(\partial B_R)} &=& \sum_{n=-\infty}^{\infty}\bigg\{\frac{|n|^{2\gamma}}{2\pi R} \overline{(d_1\cdot \nabla_x \varphi_n (x_1) )} {(\lambda_n d_2\cdot\nabla_x \varphi_n(x_2))}\bigg\}\,,\\
\label{eqn_ipcir2}
    \langle\eta_{x_1,d_1}, G_{x_2}\rangle_{H^{\gamma}(\partial B_R)} &=& \sum_{n=-\infty}^{\infty}\bigg\{\frac{ |n|^{2\gamma}}{2\pi R} \overline{(d_1\cdot \nabla_x \varphi_n (x_1) )} {(\lambda_n \varphi_n(x_2) )}\bigg\}\,,\\
\label{eqn_ipcir3}
    \langle\zeta_{x_1}, d_2\cdot \nabla G_{x_2} \rangle_{H^{\gamma}(\partial B_R)} &=& \sum_{n=-\infty}^{\infty}\bigg\{\frac{|n|^{2\gamma}}{2\pi R} \overline{(\varphi_n (x_1) )} {(\lambda_n d_2\cdot\nabla_x \varphi_n(x_2))}\bigg\}\,,\\
\label{eqn_ipcir4}
    \langle\zeta_{x_1}, G_{x_2}\rangle_{H^{\gamma}(\partial B_R)} &=& \sum_{n=-\infty}^{\infty}\bigg\{\frac{ |n|^{2\gamma}}{2\pi R} \overline{(\varphi_n (x_1) )} {(\lambda_n \varphi_n(x_2))}\bigg\}\,;
\end{eqnarray}
\begin{align}
\label{eqn_seminor1}
    |\eta_{x,d}|_{H^{\gamma}}^2 &= \sum_{n=-\infty}^{\infty}\frac{ |n|^{2\gamma}}{2\pi R} |d \cdot \nabla \varphi_n(x)|^2 \,,
\,\,   &|\zeta_{x}|_{H^{\gamma}}^2& = \sum_{n=-\infty}^{\infty}\frac{ |n|^{2\gamma}}{2\pi R} |\varphi_n(x)|^2\,;\\
\label{eqn_seminor3}
    |d\cdot \nabla G_{x}|_{H^{\gamma}}^2 &= \sum_{n=-\infty}^{\infty}\frac{ |n|^{2\gamma}}{2\pi R} |\lambda_n d \cdot \nabla \varphi_n(x)|^2\,,
\,\,  &|G_{x}|_{H^{\gamma}}^2& = \sum_{n=-\infty}^{\infty}\frac{ |n|^{2\gamma}}{2\pi R} |\lambda_n \varphi_n(x)|^2\,.
\end{align}

\subsubsection{\mm{More about} the mutually almost orthogonality property} 
\label{subsec_more_verification}
\mm{We shall focus only on the case of Sobolev scale $\gamma = 1$, and the cases of 
other $\gamma \geq 0$ follow similarly.
}

\noindent \textbf{Case 1: $V_0 = 0$.}
\mm{For given $|c|<1$, one may quickly obtain}
\begin{equation}
    \sum_{n=1}^{\infty} n c^n = \frac{c}{(1-c)^2}, \quad \sum_{n=1}^{\infty} n^2 c^n = \frac{c(1+c)}{(1-c)^3},\quad \sum_{n=1}^{\infty}n^3 c^n = \frac{c(c^2+4c+1)}{(1-c)^4},
    \quad \sum_{n=1}^{\infty}n^4c^n = \frac{c^4+11c^3+11c^2+c}{(1-c)^5} \, .
\end{equation}

We first consider $K_{4,d_1,d_2}(x_1,x_2)$. 
\mm{For convenience, we write 
$d_i = (-\sin(\alpha_i), \cos(\alpha_i))$,} $x_i = (r_i, \theta_i)$ in the polar coordinates 
and $\tilde{r}_i = {r_i}/{R}$. Using the fact that $\tilde{r}_i<1$, \eqref{eqn_ipcir1} can be simplified as 
    \begin{equation}
    \label{eqn_pwex1}
    \begin{split}
       |  \langle\eta_{x_1,d_1}, d_2 \cdot \nabla G_{x_2} \rangle_{H^{1}(\partial B_1)}  |=&\,\, \left| \sum_{n=1}^{\infty}\frac{n^{3}}{\pi (r_1\,r_2)}(\tro\,\trt)^{n}\text{cos}((n-1)( \theta_1 - \theta_2) + \alpha_1 - \alpha_2) \right|\\
        \leq &\,\, \frac{|(\tro^2 \,\trt^2\,e^{2i(\theta_1-\theta_2)}+4\,\tro\,\trt \,e^{i(\theta_1-\theta_2)}+1)|}{\pi R^2|(1-\tro\,\trt\, e^{i(\theta_1-\theta_2)})^4|} \leq \, \frac{\tro^2\, \trt^2+4\,\tro\,\trt+1}{\pi R^2(1-\tro\,\trt)^4}\,.
        \end{split}
    \end{equation}
We may notice that the above inequalities 
become equalities if $\alpha_1 - \alpha_2 = n \pi$ (i.e. $d_1 = \pm d_2 $) and $\theta_1 = \theta_2$, that is, when the maximum is attained for fixed $r_1$ and $r_2$. 
\mm{Applying a similar trick, we further obtain 
from \eqref{eqn_seminor1} and \eqref{eqn_seminor3} that }
\begin{align}
\label{eqn_pwex3}
    &|\eta_{x_1,d_1}|_{H^{1}}^2 = \sum_{n=1}^{\infty}\frac{n^4R}{\pi}\,\tro^{2n-2} = \frac{(\tro^6+11\tro^4+11\tro^2+1)}{\pi R(1-\tro^2)^5}\,, 
     &|\zeta_{x_1}|_{H^1}^2& = \sum_{n=1}^{\infty}\frac{n^{2}}{\pi R} \,\tro^{2n} = \frac{\tro^2(1+\tro^2)}{\pi R(1-\tro^2)^3}\,;\\
\label{eqn_pwex4}
    &|d_1\cdot \nabla G_{x_1}|_{H^1}^2 = \sum_{n=1}^{\infty}\frac{n^2 R}{\pi}\,\tro^{2n-2} = \frac{R(1+\tro^2)}{\pi(1-\tro^2)^3}\,,
    &\,\,|G_{x_1}|_{H^{1}}^2& = \sum_{n=1}^{\infty}\frac{R}{\pi}\,\tro^{2n}= \frac{R\,\tro^2}{\pi(1-\tro^2)}\,.
\end{align} 

\begin{figure}
\centering
    \includegraphics[scale = 0.5]{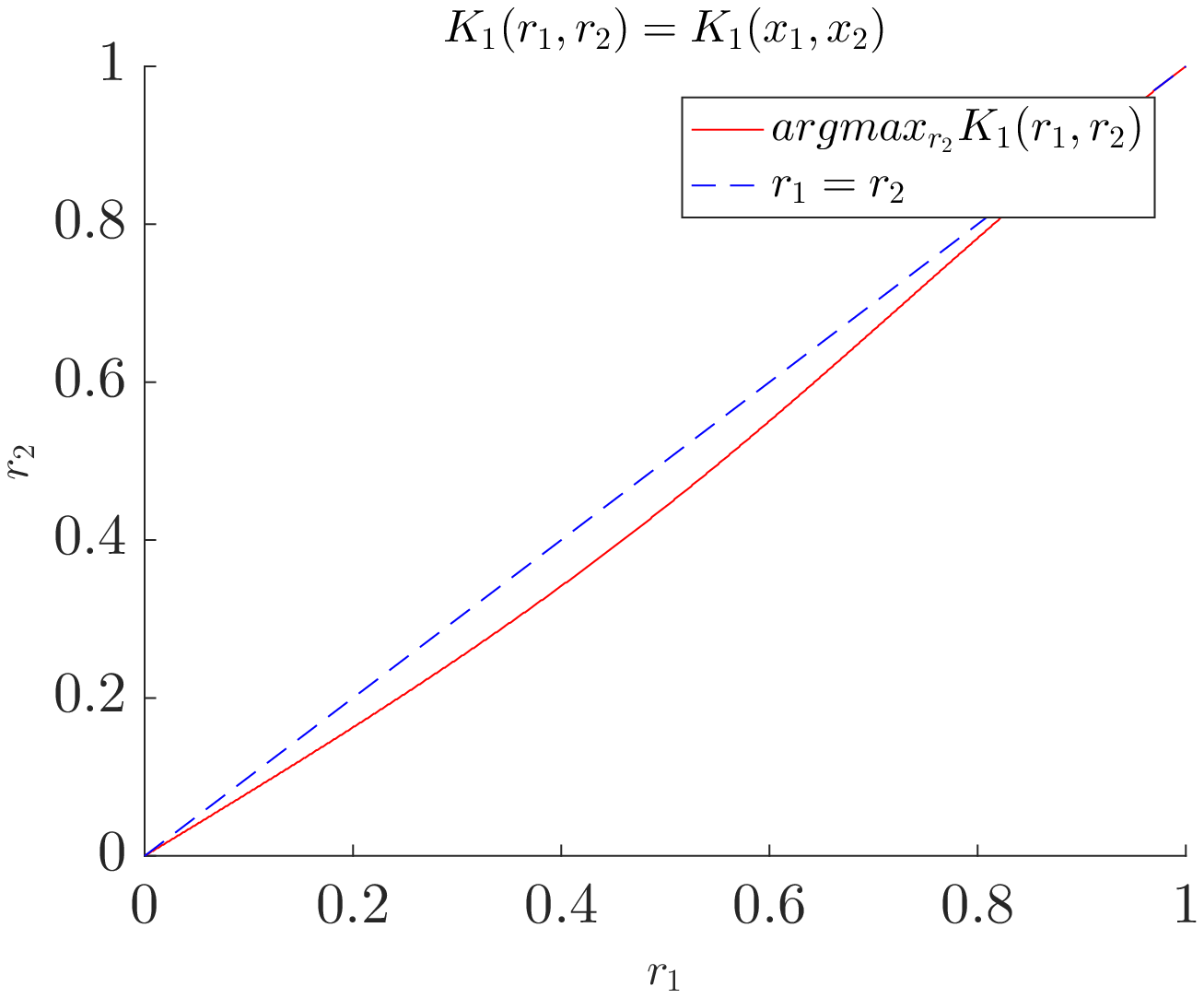}
       \includegraphics[scale = 0.5]{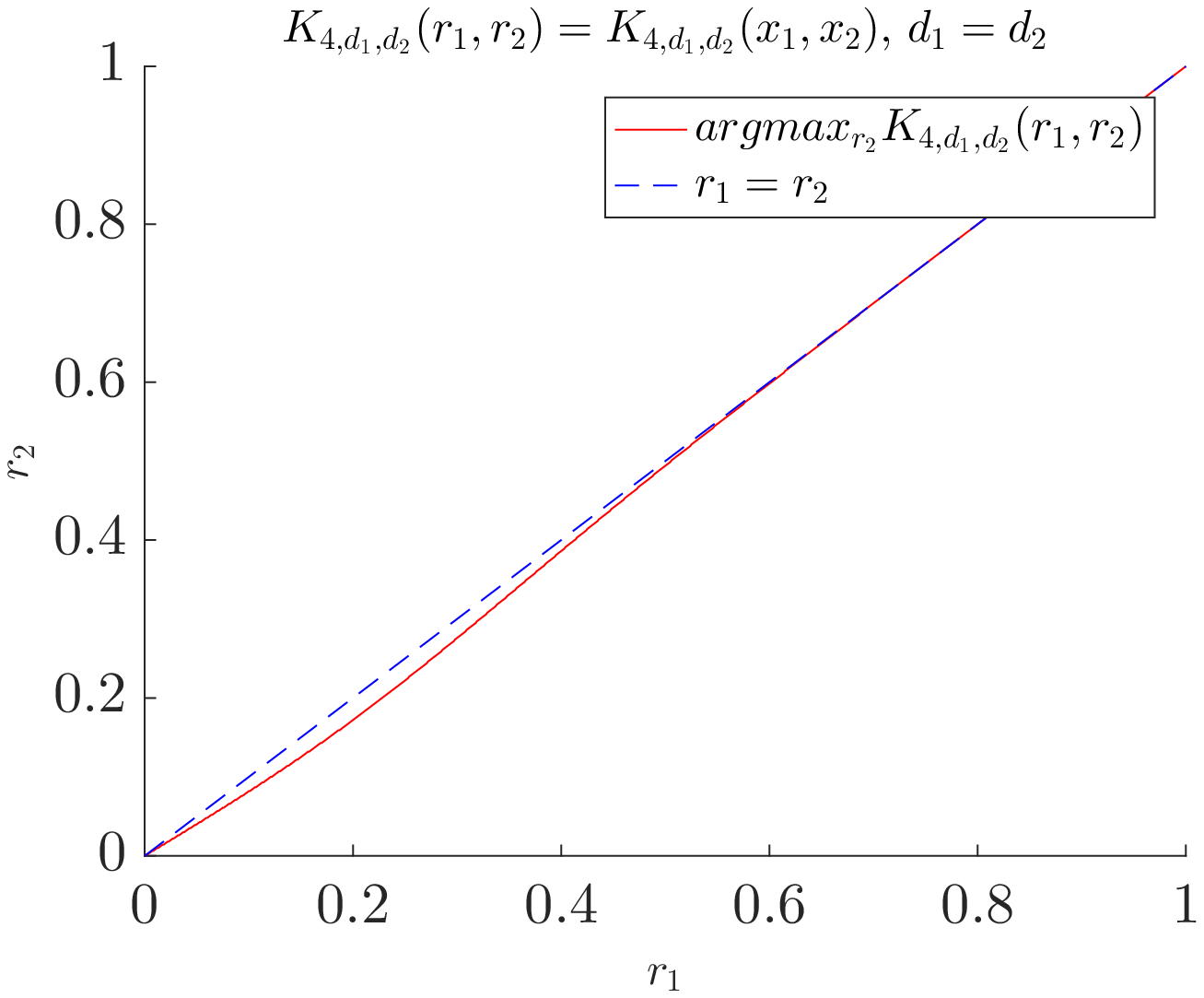}
    \caption{The location of the maximum value of kernels $K_1(x_1, x_2)$ and $K_{4,d_1,d_2}(x_1,x_2)$ defined in \eqref{kernal_mo} and \eqref{kernal_di} when $V_0 = 0$, under $\gamma = 1$, $m_i = n_i = 1/2$ ($i=1$, $2$), $d_1 = d_2$, and $\theta_1 = \theta_2$, where $x_i = (r_i,\,\theta_i)$.}
    \label{ortho_V0_vanish}
\end{figure}

To better understand the behaviour of the kernel $K_{4, d_1, d_2} (x_1, x_2)$, let us fix $\theta_1 = \theta_2$ and 
$r_1$ in \eqref{eqn_pwex1} for the time being.  
Then we like to check if the maximum of $K_{4,d_1,d_2}$, which is now a rational function of $r_2$, 
is attained when $r_2 \approx r_1$.  
While the explicit optimum is hard to find analytically, we can obtain them by solving \mm{the KKT optimality system 
via numerical approximations.} The second plot in 
Fig.\,\ref{ortho_V0_vanish} shows the value of $r_2$ that maximizes $K_{4,d_1,d_2}(x_1,\,x_2)$ with $m_1 = m_2 = 1/2$, $d_1 = d_2$, and $\theta_1 = \theta_2$. We may observe that, the function $\text{argmax}_{r_2}K_{4,d_1,d_2}(x_1,x_2)$ is very close to the linear function $r_1 = r_2$. 
For instance, we may check that when $r_1 = 0.4$, the maximum value is attained when $r_2 \approx 0.386$; and when $r_1 = 0.6$, the maximum value is attained when $r_2 \approx 0.598$. 
Therefore, we can verify the almost orthogonality property numerically in the most part of the domain 
$\Omega$ for $K_{4,d_x,d_z}$. 

We next study $K_1 (x_1,x_2)$ defined as in \eqref{kernal_mo}.  We can similarly deduce the explicit expression of the numerator of $K_1$ when $\gamma = 1$ as
\begin{equation}
    \label{K1_exp_V0}
    \begin{split}
        | \langle\zeta_{x_1}, G_{x_2}\rangle_{H^{1}(\partial B_R)}  | = &\,\, \left| \sum_{n=1}^{\infty}\frac{n}{\pi}\,(\tro\,\trt)^{n}\cos(n\theta_1-n\theta_2) \right| 
        = \,\, \left|  \text{Re}\left\{ \frac{\tro\,\trt e^{i(\theta_1-\theta_2)}}{\pi(1-\tro\,\trt\,e^{i(\theta_1-\theta_2)})^2 }\right\} \right| \\
        \leq &\,\, \frac{\tro\,\trt}{\pi}\frac{|e^{i(\theta_1-\theta_2)}|}{|1-\tro\,\trt\,e^{i(\theta_1-\theta_2)}|^2}\leq\, \frac{\tro\,\trt}{\pi(1-\tro\,\trt)^2}\,.
        \end{split}
    \end{equation}
We can see that the equalities hold when $\theta_1 = \theta_2$ in \eqref{K1_exp_V0}, that is, when the maximum is achieved for fixed ${r}_1$ and ${r}_2$. Let us now fix $\theta_1 = \theta_2$ and $r_1$ in \eqref{K1_exp_V0}, 
we like to check again if the maximum of $K_1$, which is a rational function of $r_2$, is attained when $r_2 \approx r_1$.  
\mm{Similarly, 
we may approximate them by solving the KKT optimality system via numerical approximations.}
The first plot in Fig.\,\ref{ortho_V0_vanish} describes the value of $r_2$ that maximizes $K_{1}(x_1,x_2)$ with $n_1 = n_2 = 1/2$. 
We may observe that, the function $\text{argmax}_{r_2}K_{1}(x_1,x_2)$ is very close to the linear function $r_1 = r_2$. For instance, we may check that when $r_1 =0.4 $, the maximum occurs at $r_2\approx 0.342$; and when $r_1 = 0.7$, the maximum happens at $r_2\approx 0.666$. \mm{Therefore we have verified numerically that the maximum of $K_{1}(x_1,x_2)$ occurs} 
when $x_1$ is very close to $x_2$, which is the desired almost orthogonality property.

Now we consider the decoupling effect, i.e., to check the full version of the mutually almost orthogonality property.  
For this purpose, we like to compare behaviours of
$K_{2, d_2} (x_1,x_2)$ and $K_{3, d_1} (x_1,x_2)$ with $K_1(x_1,x_2)$ and $K_{4,d_1,d_2}(x_1,x_2)$ 
defined in \eqref{kernal_mo} and \eqref{kernal_di}.
\mm{We obtain from 
\eqref{eqn_ipcir3}} and \eqref{eqn_ipcir2} which provide  explicit representations of numerators of $K_{2,d_2}$ and $K_{3,d_1}$
that 
\begin{align}
|\langle\zeta_{x_1}, \, d_2\cdot \nabla G_{x_2} \rangle_{H^{1}(\partial B_1)} | =& \left|  \frac{1}{\pi R}\sum_{n=1}^{\infty}n^{2}\,\tro^{n}\,\trt^{n-1}\,\sin(n\theta_1-(n-1)\theta_2-\alpha_2) \right | \\\notag
 =& \,\,\frac{r_1}{\pi R^2} \left| \text{Im}\left\{\frac{e^{i(\theta_1-\alpha_2)} ( 1+ \tro\,\trt\,e^{i(\theta_1-\theta_2)} ) }{(1-\tro\,\trt\,e^{i(\theta_1-\theta_2)})^3}\right\} \right|\,,\\
\label{etag_cir} 
|\langle\eta_{x_1,d_1}, \, G_{x_2} \rangle_{H^{1}(\partial B_R)} | =  &\left|\frac{1}{\pi R}\sum_{n=1}^{\infty}n^{2}\,\tro^{n-1}\trt^{n}\,\sin(n\theta_2-(n-1)\theta_1-\alpha_1) \right | \\\notag
 = &\,\,\frac{r_2}{\pi R^2} \left| \text{Im}\left\{\frac{e^{i(\theta_2-\alpha_1)} ( 1+ \tro\,\trt\,e^{-i(\theta_1-\theta_2)} ) }{(1-\tro\,\trt\,e^{-i(\theta_1-\theta_2)})^3}\right\} \right| \,.
    \end{align}

We may now see a very interesting behaviour:  a minimum (i.e., zero) of 
$|K_{2,d_2}(x_1,x_2)|$ and $|K_{3, d_1} (x_1,x_2)|$ are attained when $\alpha_1 = \theta_2$, $\alpha_1 
= \alpha_2$, and $\theta_1 = \theta_2$. 
This is an ideal behaviour as 
the maximum of the numerator of $K_1$ and $K_{4,d_1,d_2}$  occur at $\theta_1 = \theta_2$ and $\alpha_1 = \alpha_2 $
by using \eqref{eqn_pwex1} and \eqref{K1_exp_V0}, 
therefore helps contrast $K_{2,d_2}$ and $K_{3,d_1}$ with $K_1$ and $K_{4,d_1,d_2}$.  

\begin{figure}
\centering
    \includegraphics[scale = 0.29]{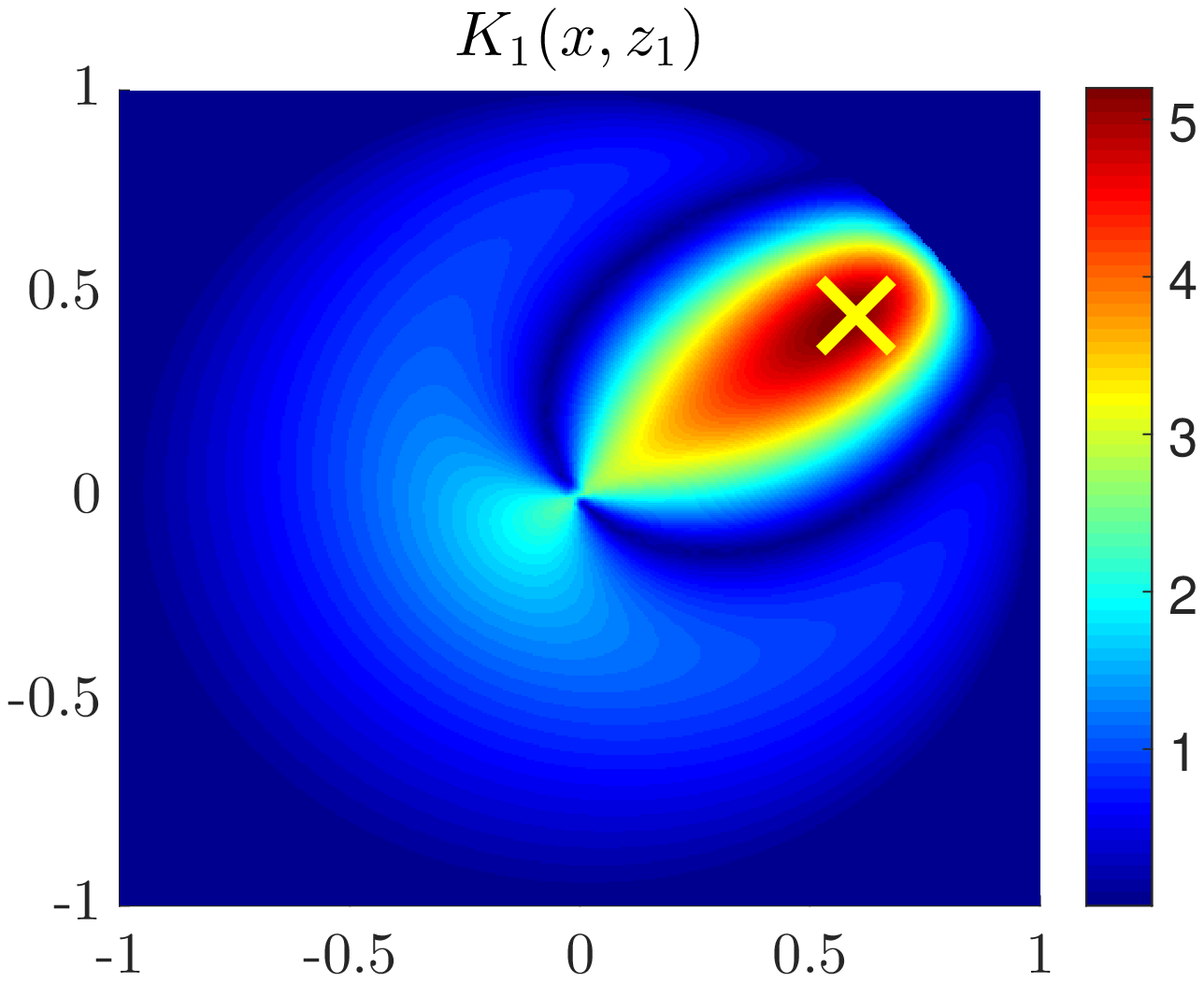}
    \includegraphics[scale = 0.29]{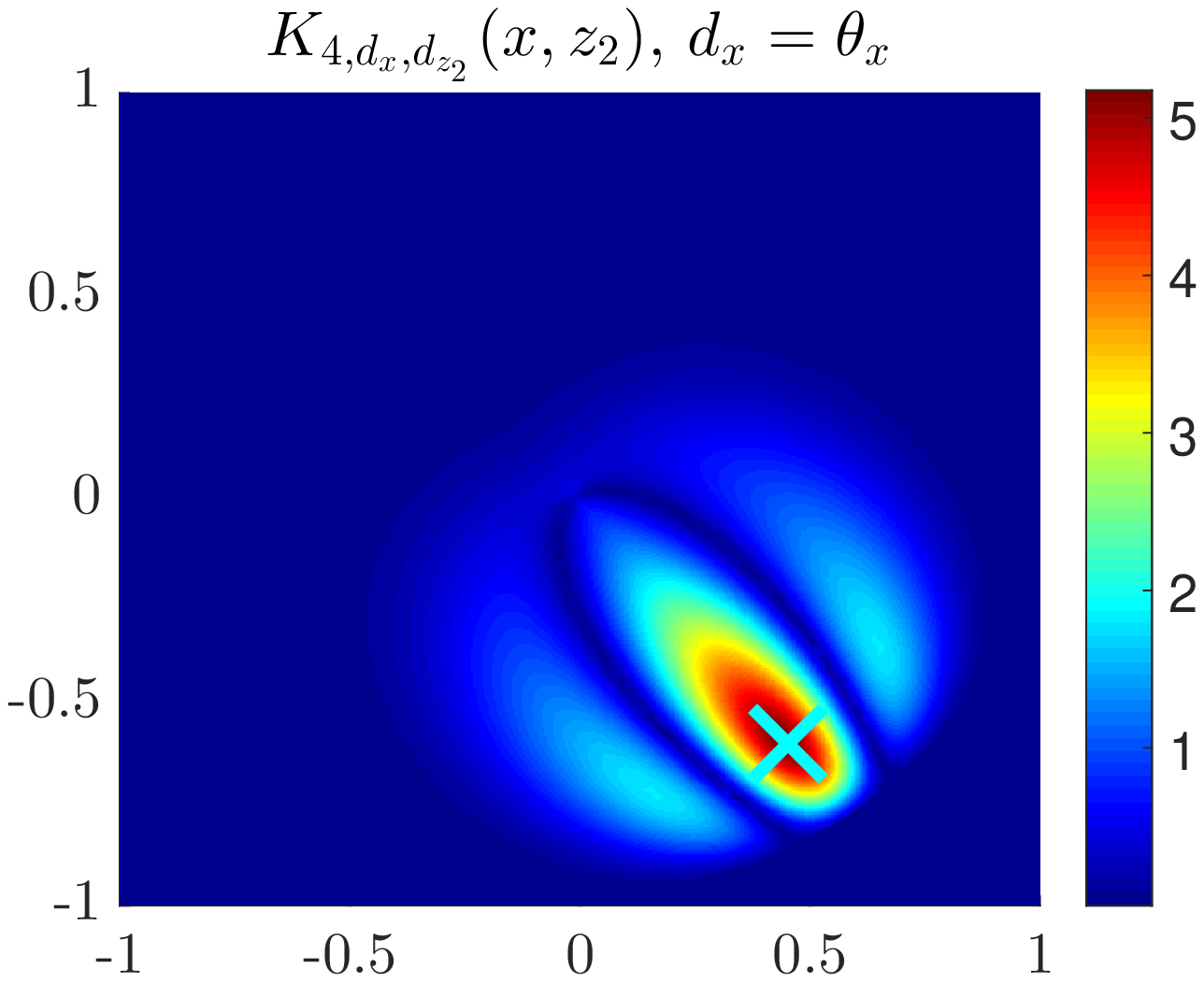}
    \includegraphics[scale = 0.29]{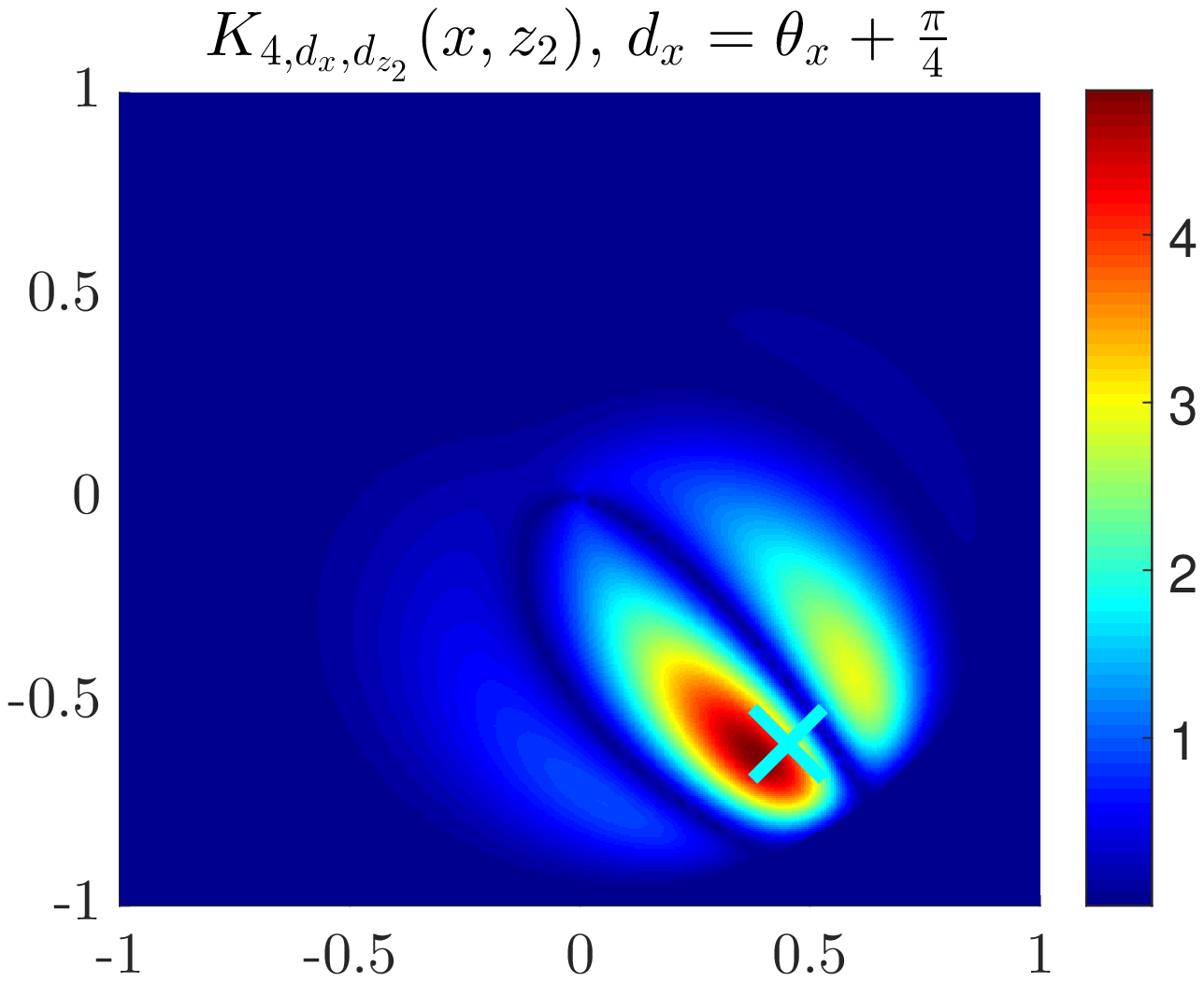}
    \includegraphics[scale = 0.29]{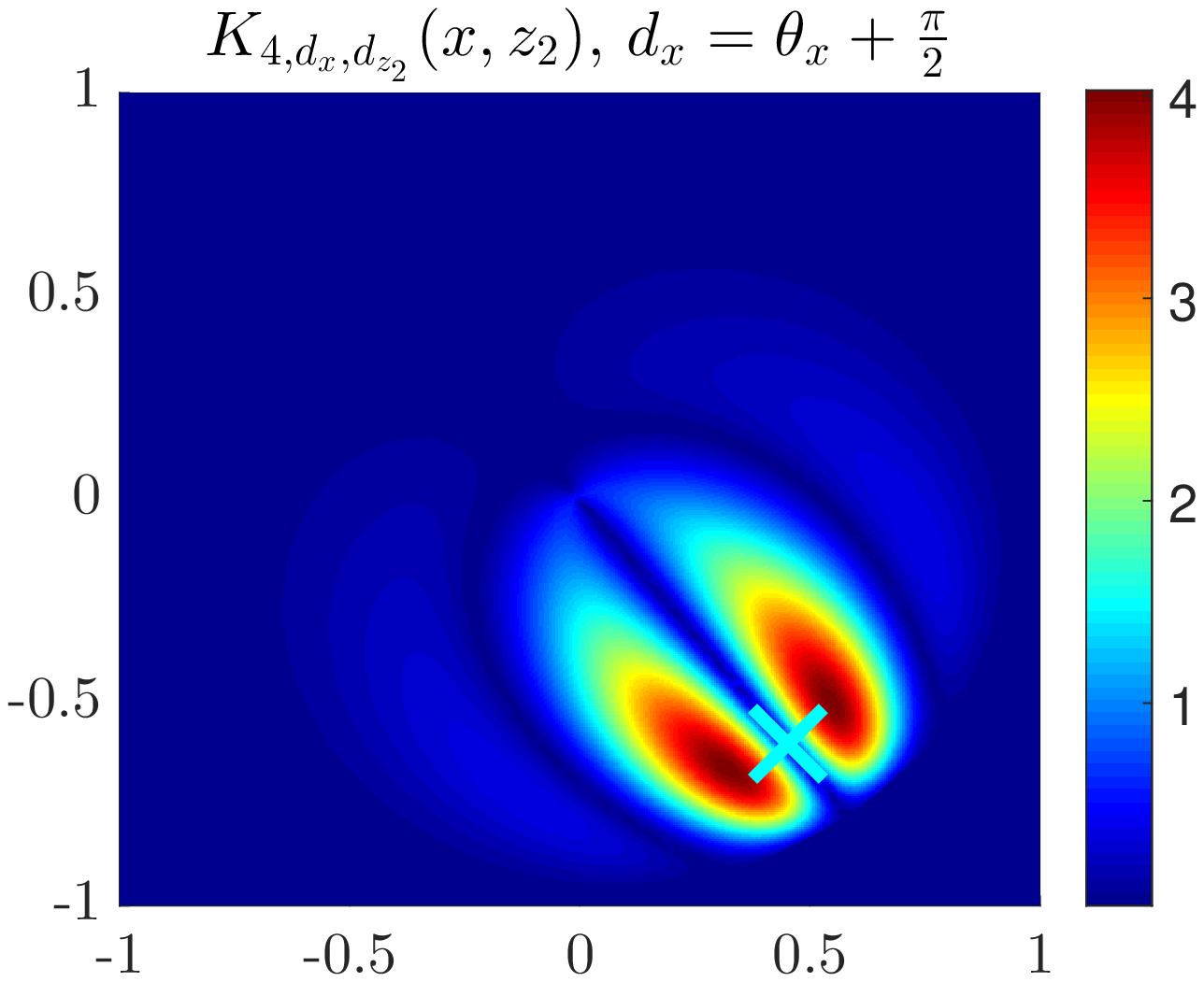}
    \caption{Almost orthogonality property of $K_1(x,z_1)$ and $K_{4,d_x,d_{z_2}}(x,z_2)$ for $V_0 = 0$, with $m_i = n_i = 1/2$ ($i=1,2$) and $z_1 = (0.6,0.45)$, $z_2 = (0.45,-0.6)$. Directions in $K_{4,d_x,d_{z_2}}(x,z_2)$ are chosen as $d_x =\theta_x$, $d_x = \theta_x+\pi/4$, $d_x = \theta_x+\pi/2$ (from left to right), and $d_{z_2} = \theta_{z_2}$.}
    \label{V0_K14}
\end{figure}

\begin{figure}
\centering
    \includegraphics[scale = 0.35]{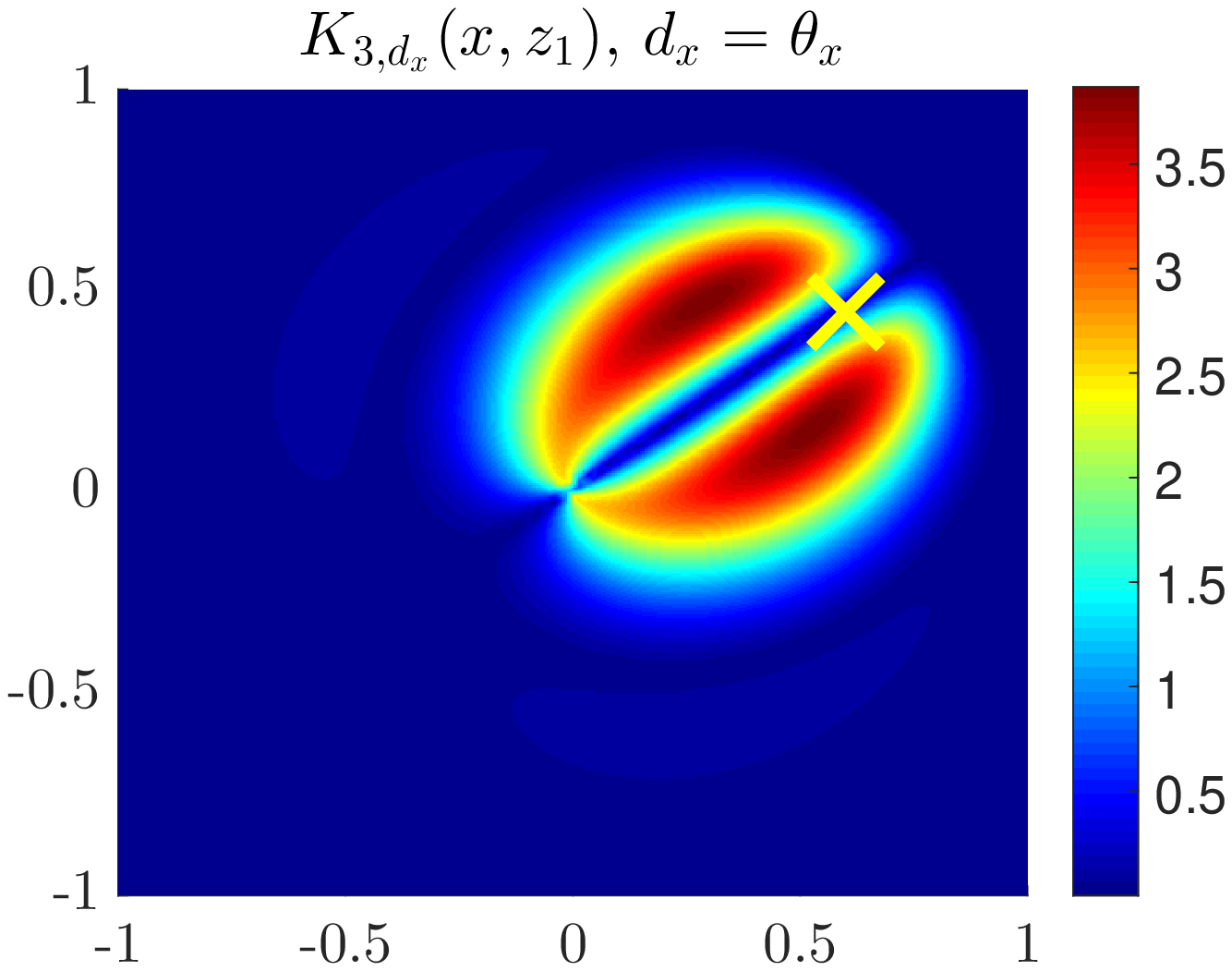}
    \includegraphics[scale = 0.35]{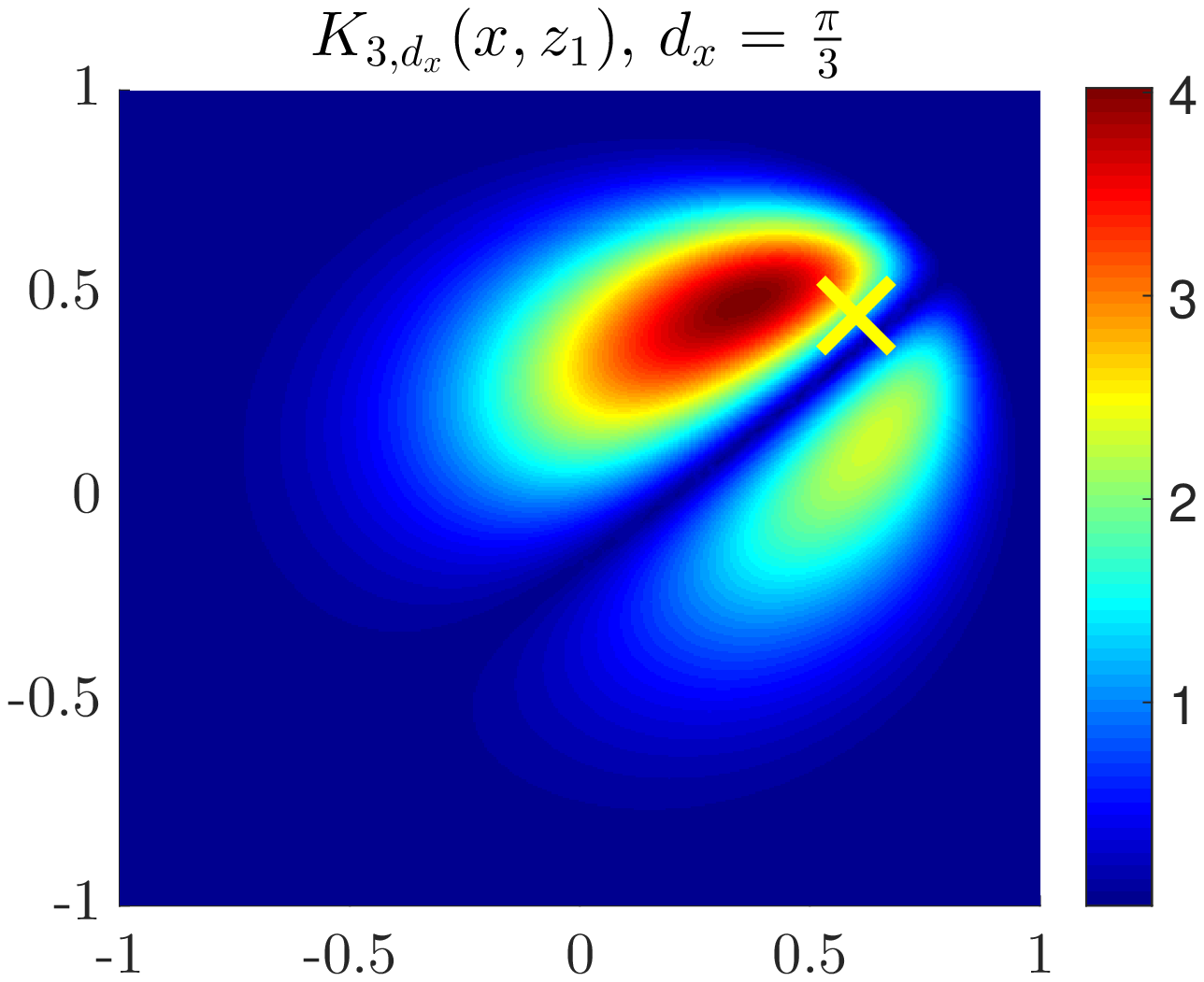}
    \includegraphics[scale = 0.35]{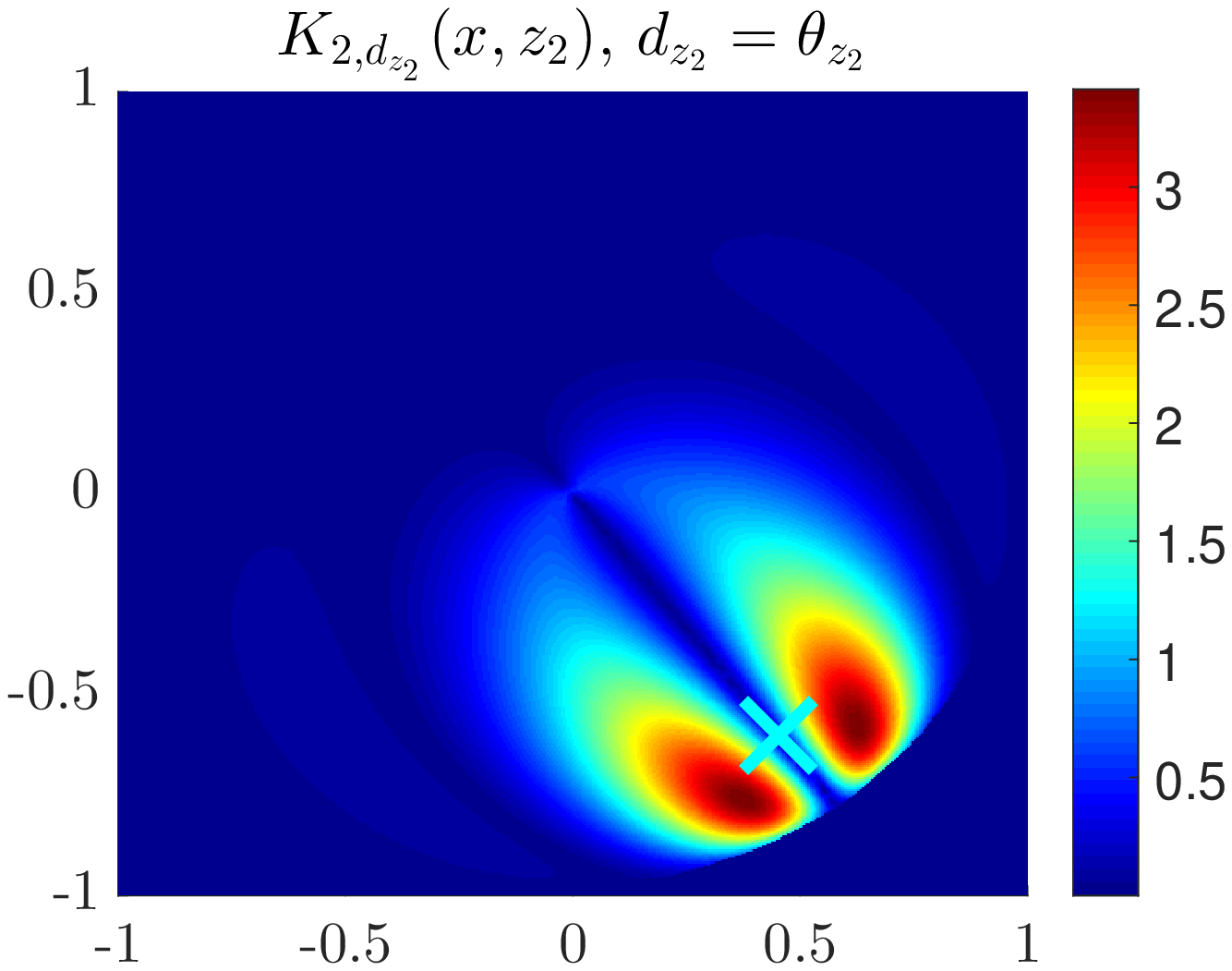}
    \caption{Mutually almost orthogonality property of $K_{3,d_x}(x,z_1)$ and $K_{2,d_{z_2}}(x,z_2)$ for $V_0= 0$, with  $m_i = n_i = 1/2$ ($i=1,2$), and $z_1 = (0.6,0.45)$, $z_2 = (0.45,-0.6)$. Directions are chosen as $d_{x}= \theta_{x}$, $d_x = \pi/3$, and $d_{z_2} = \theta_{z_2}$ (from left to right).}
    \label{V0_K23}
\end{figure}

\begin{figure}
\centering
         \includegraphics[scale = 0.37]{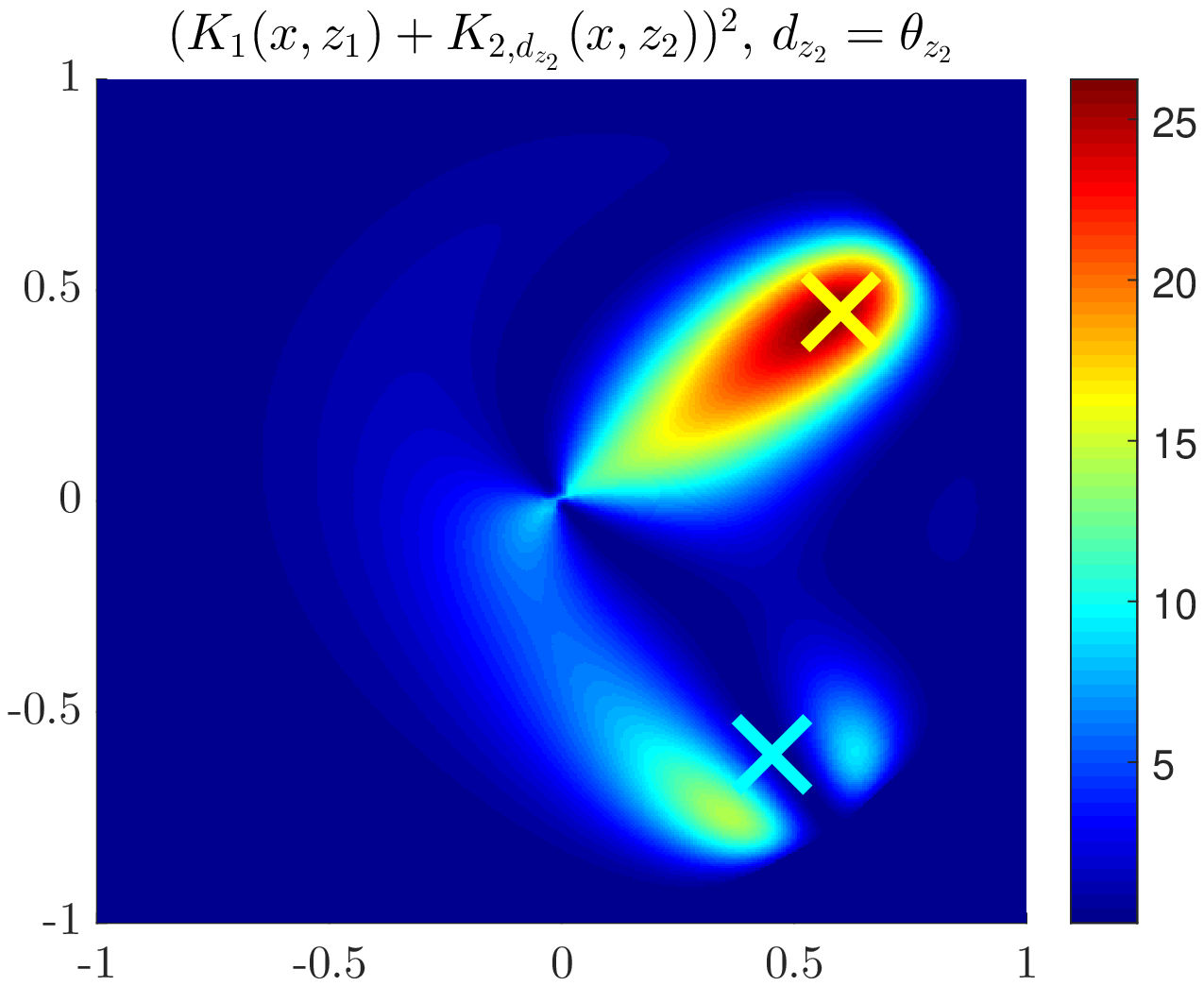}
    \includegraphics[scale = 0.37]{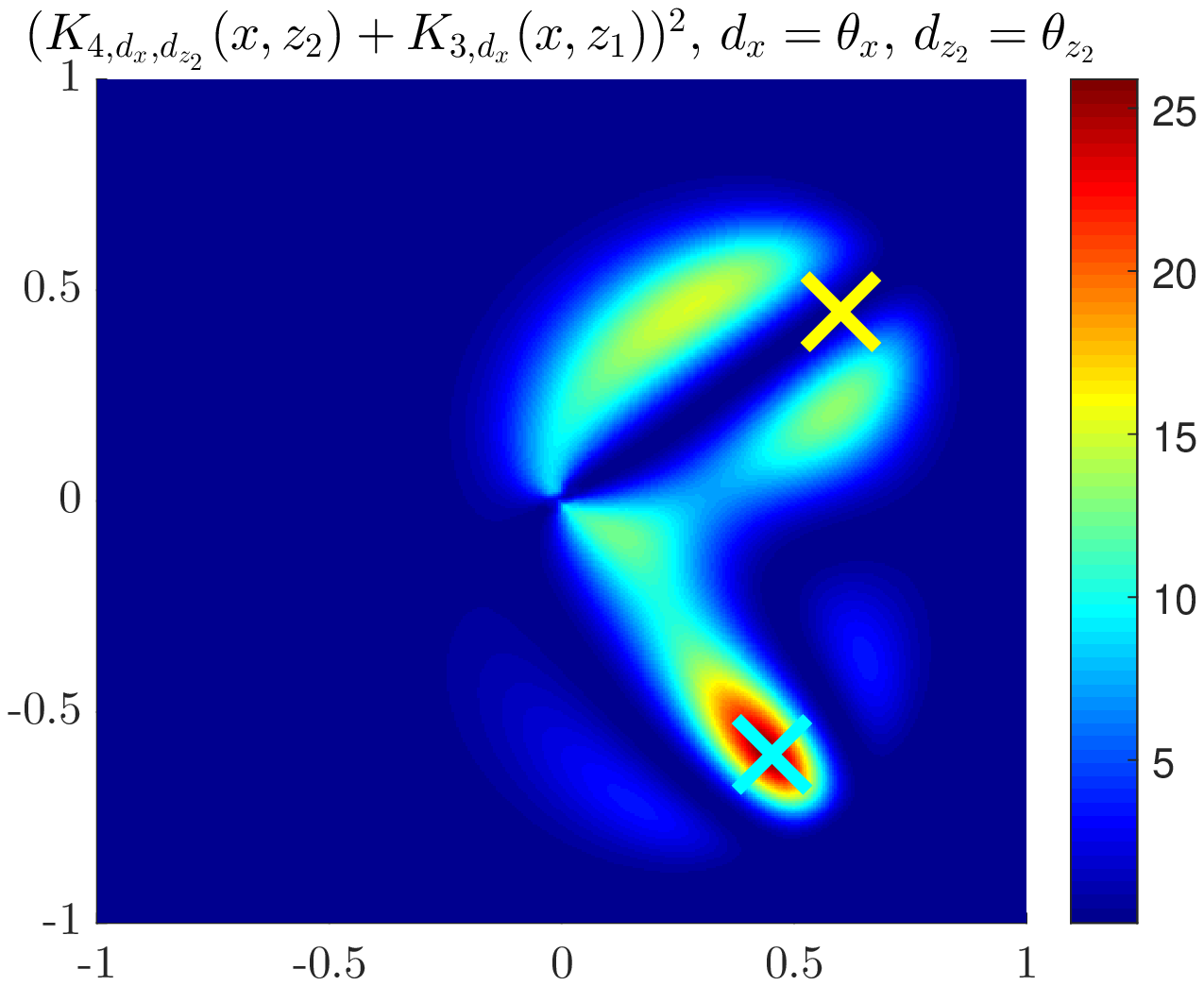}
     \includegraphics[scale = 0.37]{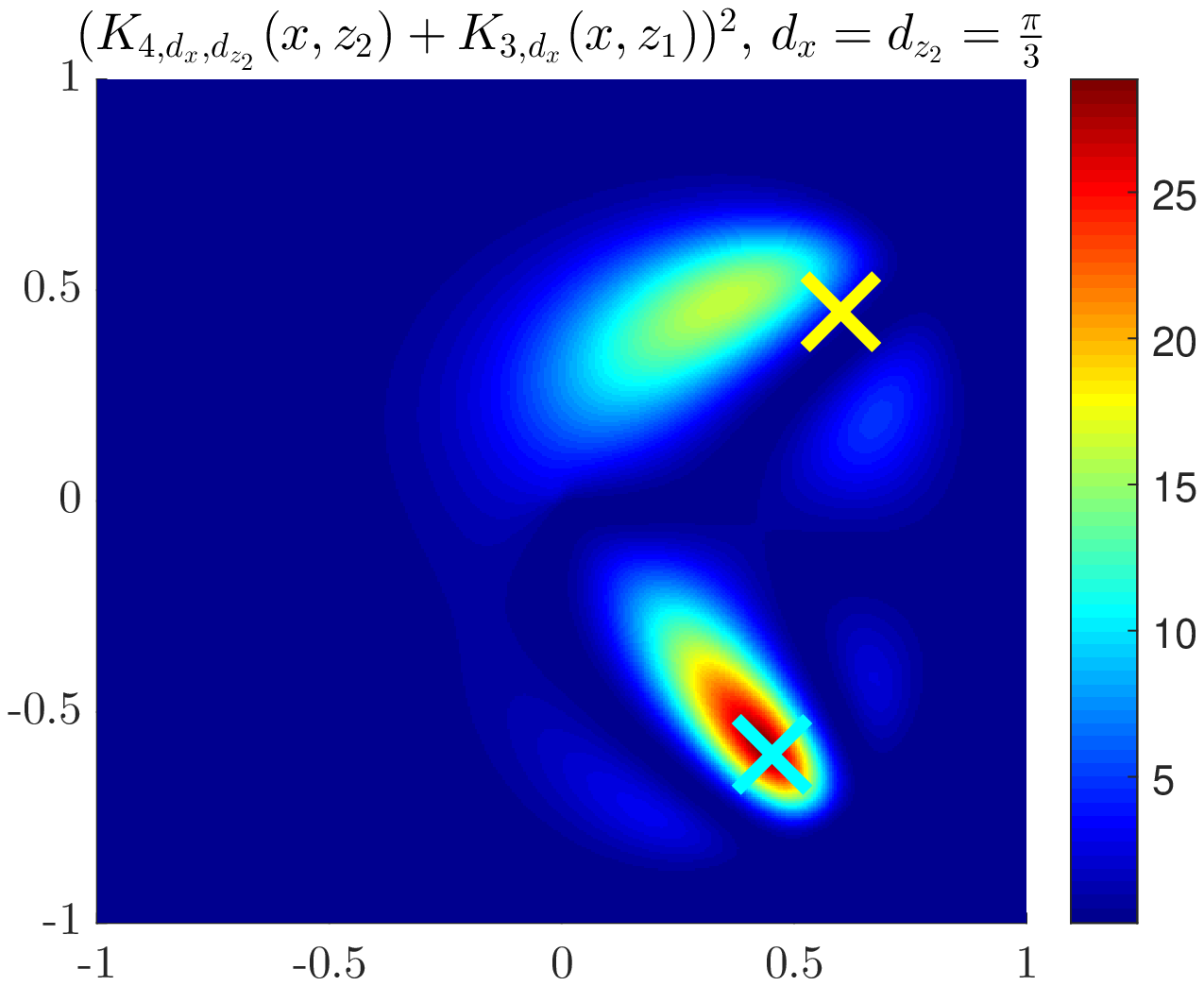}
    \caption{Mutually almost orthogonality property of $K_{1}(x,z_1)+K_{2,d_{z_2}}(x,z_2)$ (the left plot) and $K_{4,d_x,d_{z_2}}(x,z_2)+K_{3,d_x}(x,z_1)$ (the middle and right plots) for $V_0= 0$, with $m_i = n_i = 1/2$ ($i=1$, $2$), and $z_1 = (0.6,0.45)$, $z_2 = (0.45,-0.6)$. Directions are chosen as $d_{z_2} = \theta_{z_2}$, $ d_x = \theta_x$, and $d_x = d_{z_2} = \pi/3$ (from left to right).}
    \label{V0_K}
\end{figure}

In Figs.\,\ref{V0_K14}-\ref{V0_K}, mutually almost orthogonality properties are further studied through numerical experiments 
for $R=1$. From these results, we may see that there is a monopole located at $z_1 = (0.6,\,0.45)$ and a dipole located at $z_2 = (0.45,\, -0.6)$. To clearly illustrate the decoupling effect by considering the situation when the influence of the monopole and the dipole on the boundary are comparable, the monopole $G_{z_1}$ is multiplied by a constant $6$ with respect to our expressions in \eqref{K1_exp_V0} and \eqref{etag_cir}. We also take $m_i = n_i = 1/2$ ($i=1$, $2$) and denote the locations of $z_1$ and $z_2$ using a yellow cross and a blue cross respectively. In what follows, $d = \theta_x$ represents $d = (-\sin(\theta_x),\,\cos(\theta_x))^T$, where $\theta_x$ is the angular coordinate in polar coordinates for $x$.
\begin{enumerate}
\item In Fig.\,\ref{V0_K14}, the first plot is $K_1(x,z_1)$ for $x\in\Omega$. This plot demonstrates the desired property of $K_1$, 
and we notice that the maximum occurs when $x$ is very close to $z_1$. We then assume $d_{z_2} = \theta_{z_2}$; 
the second plot in Fig.\,\ref{V0_K14} is $K_{4,d_x,d_{z_2}}(x,z_2)$, with $d_x= \theta_x$. 
We can observe that the maximum occurs when $x\approx z_2$, given the appropriate probing direction. The third plot is for 
$K_{4,d_x,d_{z_2}}(x,z_2)$ with $d_{x} = \theta_x + \pi/4$. 
We notice that even if there is a moderate perturbation from the best probing direction ($\theta_x = \theta_{z_2}$), the maximum of the kernel function is not very far away from the point $z_2$. The last plot is the case when $d_x = \theta_{x} + \pi/2$. In this case, two peaks of the kernel function appear around the point with a dipole, and the maximum value in the figure is smaller than the case when $d_x = \theta_x$. This illustrates that a reasonable probing direction is essential for the accurate determination of the location of a dipole.

\item In Fig.\,\ref{V0_K23}, we demonstrate behaviours of $K_{3,d_x}(x,z_1)$ with $d_x = \theta_{x}$, $d_x = \pi/3$ and $K_{2,d_{z_2}}(x,z_2)$ with $d_{z_2} = \theta_{z_2}$ from left to right. There are two important observations: 
\mm{the maxima of $K_{2,d_z}$ and $K_{3,d_x}$ are smaller than that of $K_1$ and $K_{4,d_x,d_z}$; 
for the case $d_x= \theta_x$, the maximum appears at two sides of the point $z_i$ instead of being right at the spot.
}

\item 
In Fig.\,\ref{V0_K}, we examine the coexistence of a monopole at $z_1 = (0.6,0.45)$ and a dipole at $z_2 = (0.45,-0.6)$. 
The first plot can be considered as probing by $\zeta_x$, while the second and third plots can be considered as probing by $\eta_{x,d_x}$ under different probing directions. We may conclude that the monopole probing function $\zeta_x$ 
	interacts better with the monopole located at $z_1$, while the dipole probing function $\eta_{x,d}$ 
	interacts better with the dipole located at $z_2$, under an appropriate probing direction.
\end{enumerate}

%

%


%

\noindent \textbf{Case 2: $V_0 \neq 0$.} 
In this case, \mm{the kernel functions are expressed in terms of Bessel functions. 
A closed formula is hard to obtain, so we will verify the mutually almost orthogonality property 
mainly through numerical experiments. 
}
We first derive \mm{the explicit representations of the} numerators of $K_1$, $K_{2,d_z}$, $K_{3,d_x}$, $K_{4,d_x,d_z}$ through \eqref{eqn_ipcir1} to \eqref{eqn_ipcir4}
\begin{align}
   &     | \langle\zeta_{x_1},\, G_{x_2}\rangle_{H^1(\partial B_1)} | = \frac{1}{2\pi R k} \left | \sum_{n\in \mathbb{Z}}\bigg[e^{in(\theta_2-\theta_1)}|n|^{2}\frac{I_n(kr_1)I_n(kr_2)}{I_n'(kR)I_n(kR)}\bigg] \right |\,. \\
 &  | \langle \zeta_{x_1},\, d_2\cdot \nabla G_{x_2}\rangle_{H^1(\partial B_1)}| = \frac{1}{2 \pi Rk}\bigg|\sum_{n\in \mathbb{Z}}\frac{e^{in(\theta_2 - \theta_1)}{|n|^{2}} I_n(kr_1)}{I_n(kR)I_n'(kR)}\binom{\sin(\theta_2 - \alpha_2)}{\cos(\theta_2 - \alpha_2)}^T\binom{kI_n'(kr_2)}{inI_n(kr_2)/r_2}\bigg|\,.\\
 &  | \langle \eta_{x_1,d_1}, \,G_{x_2}\rangle_{H^1(\partial B_1)}  | = \frac{1}{2 \pi Rk} \left |\sum_{n\in \mathbb{Z}}\bigg[\frac{e^{in(\theta_2 - \theta_1)}{|n|^{2}} I_n(kr_2)}{I_n(kR)I_n'(kR)}\binom{\sin(\theta_1 - \alpha_1)}{\cos(\theta_1 - \alpha_1)}^T\binom{kI_n'(kr_1)}{-inI_n(kr_1)/r_1}\bigg] \right| \,.
\end{align}
\begin{align}
      & | \langle \eta_{x_1,d_1},\, d_2\cdot \nabla G_{x_2}\rangle_{H^1(\partial B_1)} |\\ 
        =& \,\,\frac{1}{2\pi Rk}\bigg|\sum_{n\in \mathbb{Z}}\bigg[\frac{e^{in(\theta_2-\theta_1)}|n|^2 }{I_n(kR) I_n'(kR)}\binom{\sin(\theta_1 - \alpha_1)}{\cos(\theta_1 - \alpha_1)}^T\binom{kI_n'(kr_1)}{-inI_n(kr_1)/r_1}\binom{\sin(\theta_2 - \alpha_2)}{\cos(\theta_2 - \alpha_2)}^T{\binom{kI_n'(kr_2)}{inI_n(kr_2)/r_2}}\bigg]\bigg|\,.\notag
\end{align}

Similarly, the explicit expressions for $H^{\gamma}$ semi-norms can be derived from \eqref{eqn_seminor1} and \eqref{eqn_seminor3} as
\begin{align}
    |\eta_{x_1,d_1}|^2_{H^1} &= \sum_{n\in \mathbb{Z}} \frac{|n|^2\bigg[(\cos(\theta_1-\alpha_1)\frac{n}{r_1}I_n(kr_1))^2+(\sin(\theta_1-\alpha_1)kI_n'(kr_1))^2\bigg]}{2\pi R I_n(kR)^2}\,,\\
    |d_1 \cdot \nabla G_{x_1} |^2_{H^1} &= \sum_{n\in \mathbb{Z}} \frac{|n|^2\bigg[(\cos(\theta_1-\alpha_1)\frac{n}{r_1}I_n(kr_1))^2+(\sin(\theta_1-\alpha_1)kI_n'(kr_1))^2\bigg]}{ 2\pi R k^2 I_n'(kR)^2}\,.
\end{align}
\begin{equation}
\label{zetanorm_V}
    |\zeta_{x_1}|_{H^1}^2 = \sum_{n=1}^{\infty}\frac{n^2}{\pi R}\frac{I_n(kr_1)^2}{I_n(kR)^2}\,,\qquad |G_{x_2}|_{H^1}^2 = \sum_{n=1}^{\infty}\frac{n^2}{\pi Rk^2}\frac{I_n(kr_2)^2}{I_n'(kR)^2}\,.
\end{equation}

\begin{figure}
\centering
    \includegraphics[scale = 0.35]{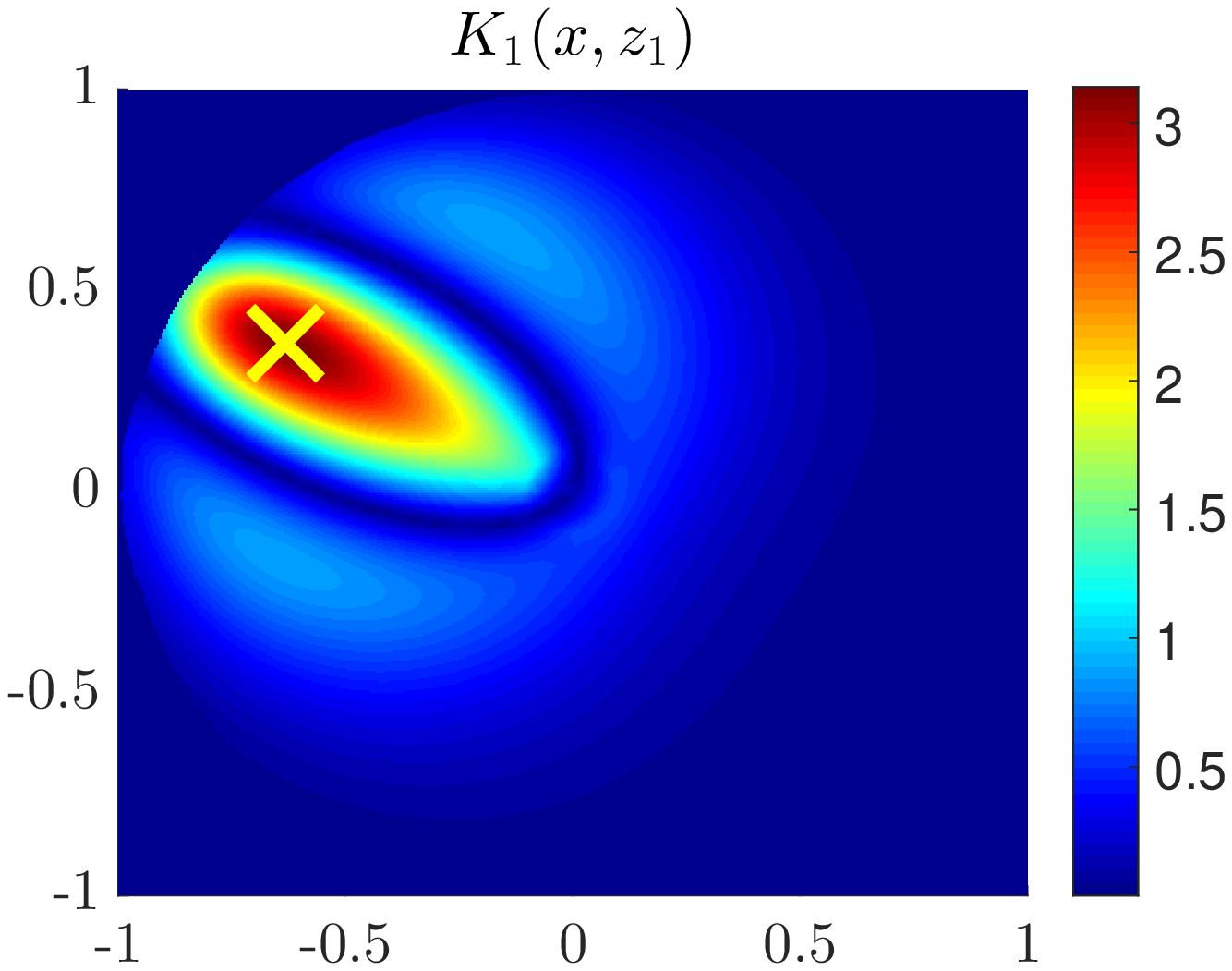}
    \includegraphics[scale = 0.35]{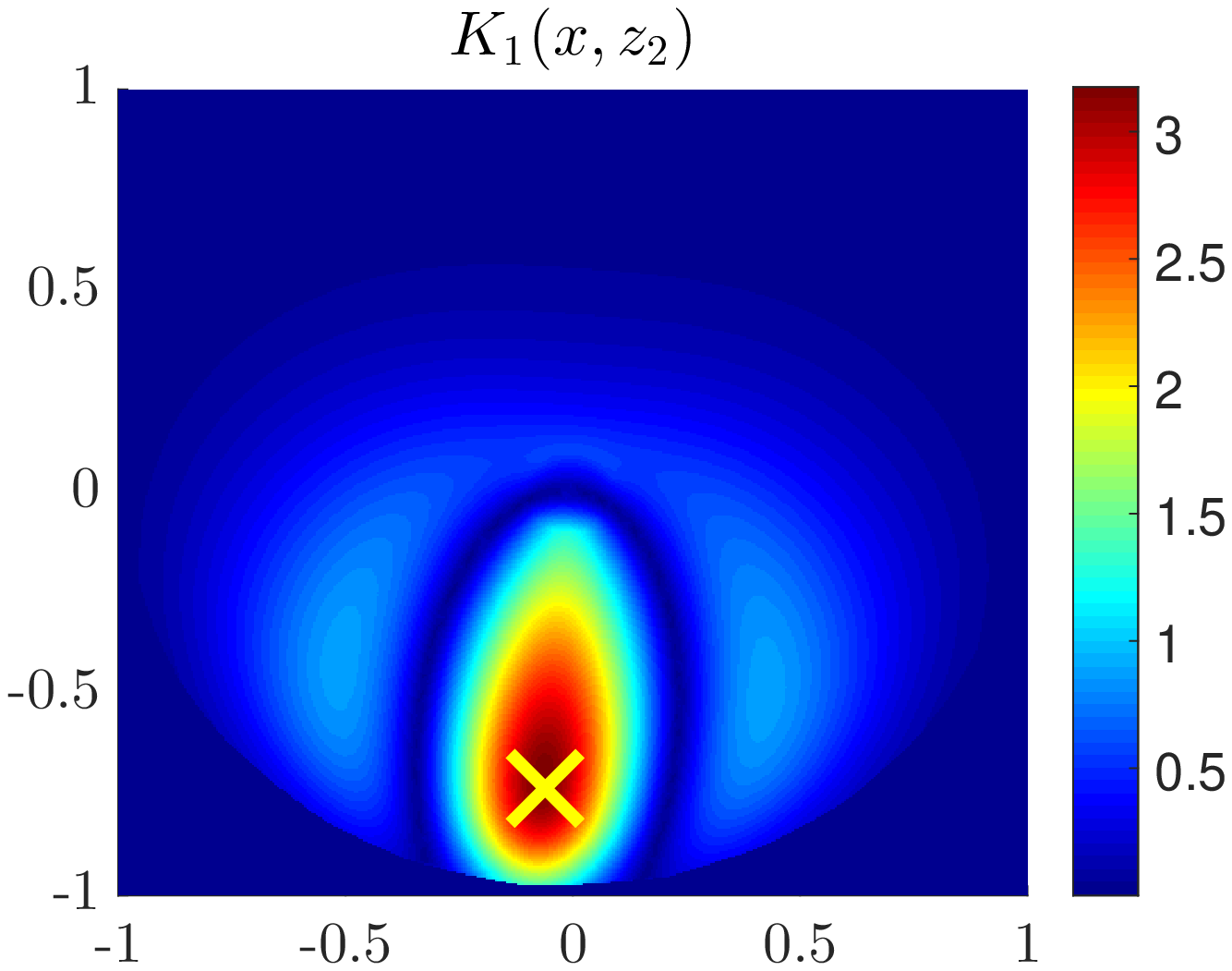}
    \includegraphics[scale = 0.35]{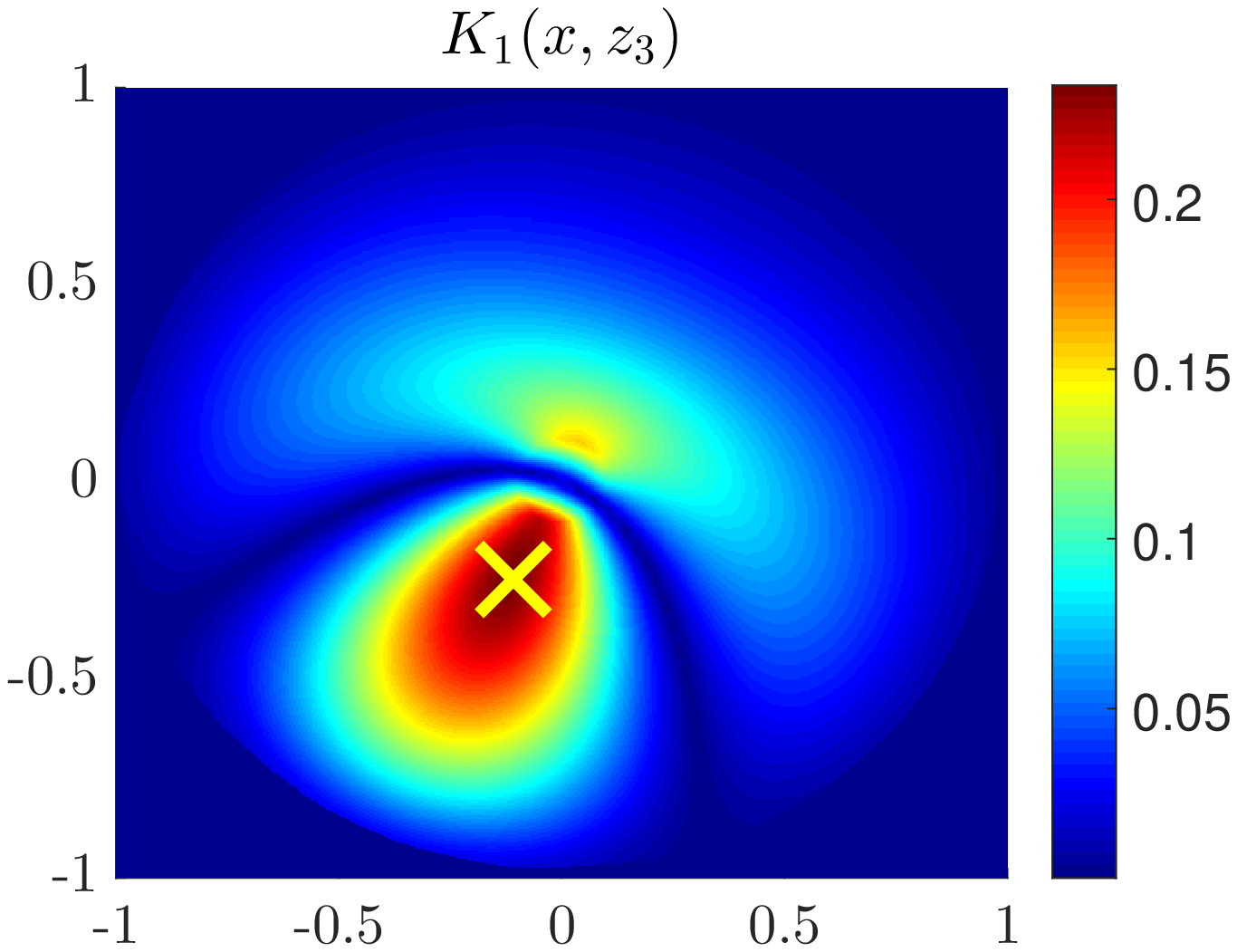}
    \caption{\mm{Almost orthogonality property of $K_1(x,z_i)$ for $V_0\neq 0$, with $n_1 = n_2 = 1/2$, and $z_1 = (-0.63,0.37)$, $z_2 = (-0.06,-0.73)$, $z_3 = (-0.11, -0.24)$ (from left to right).}}
    \label{V_K1}
\end{figure}

\begin{figure}
\centering
    \includegraphics[scale = 0.35]{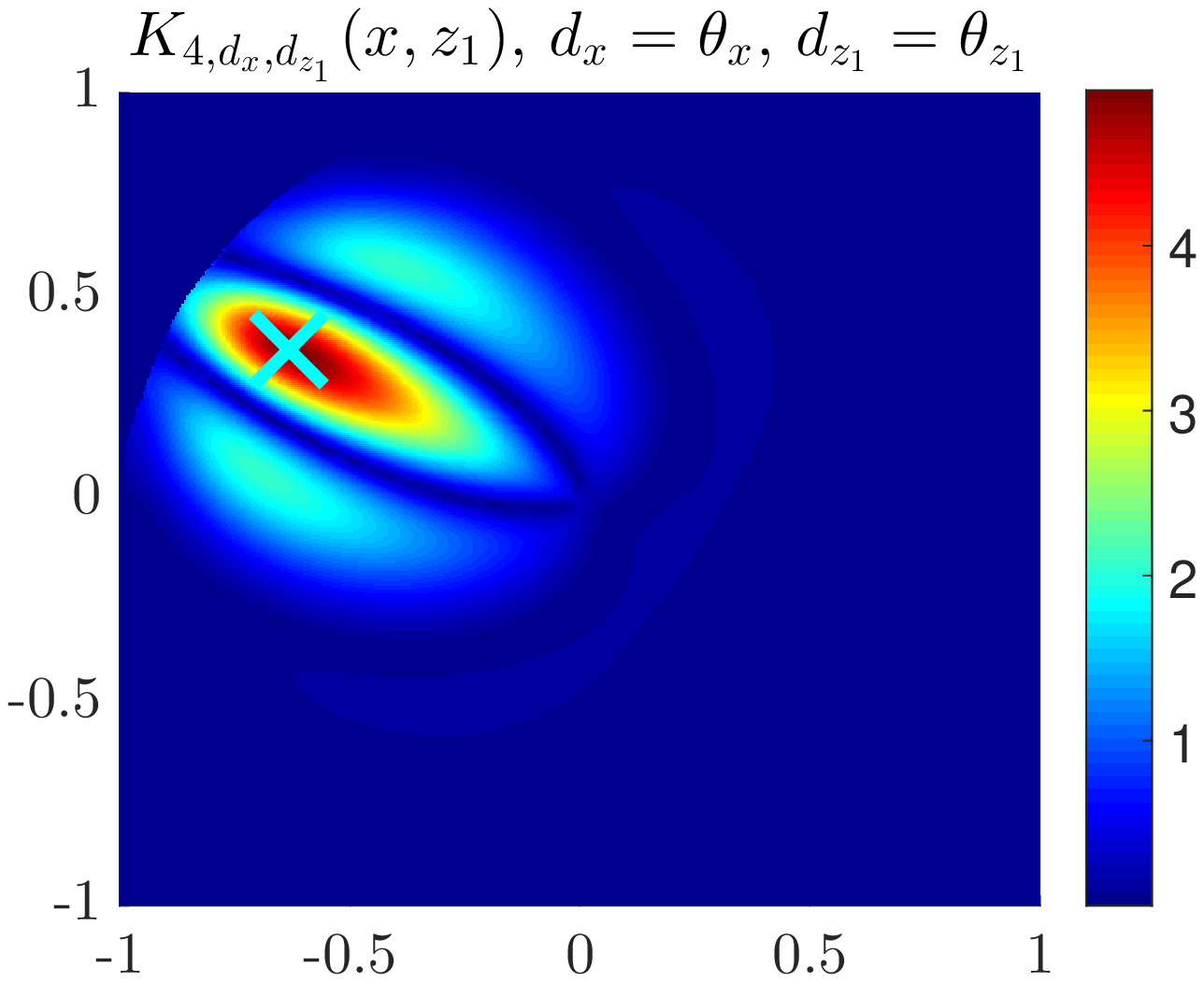}
    \includegraphics[scale = 0.35]{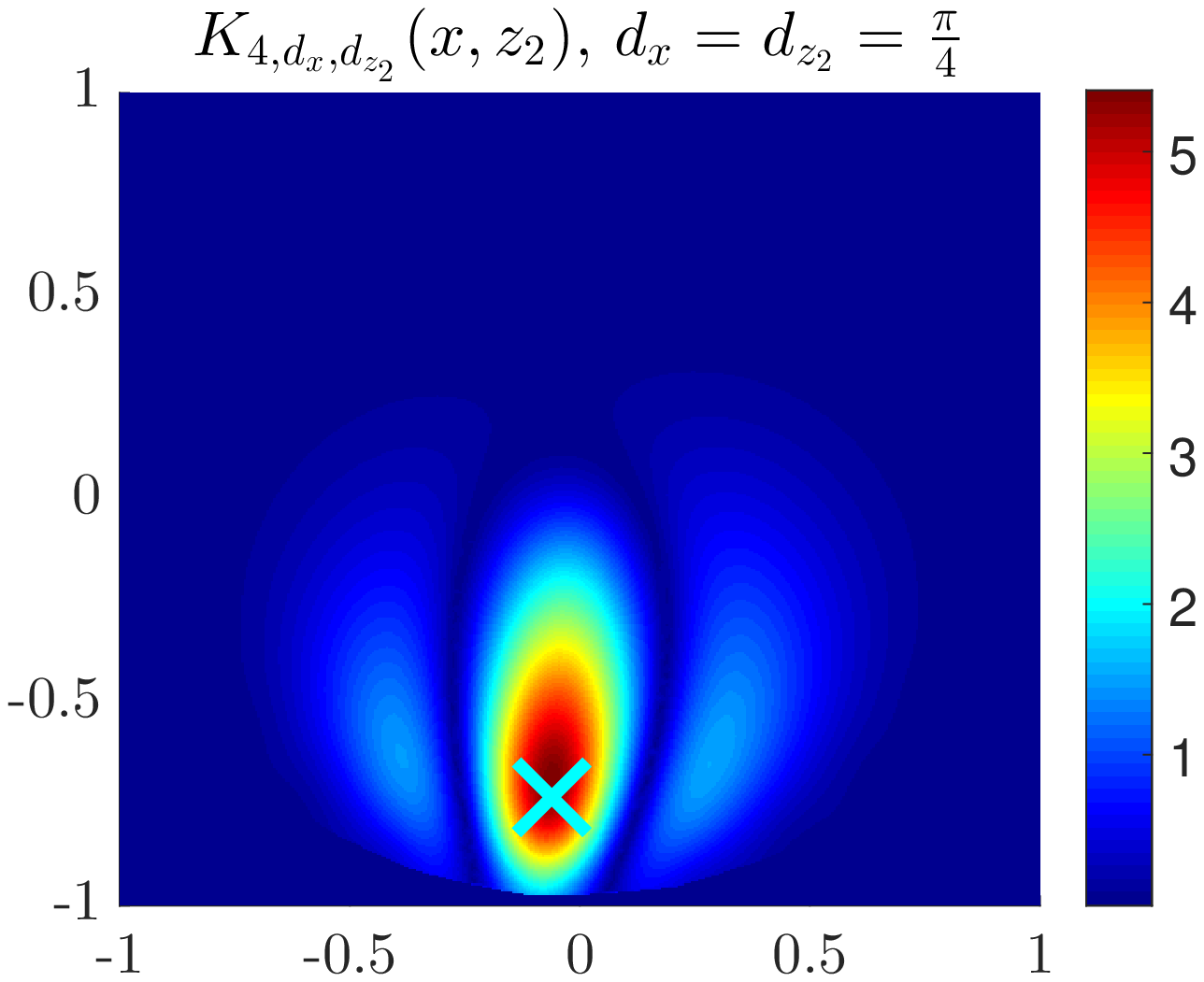}
    \includegraphics[scale = 0.35]{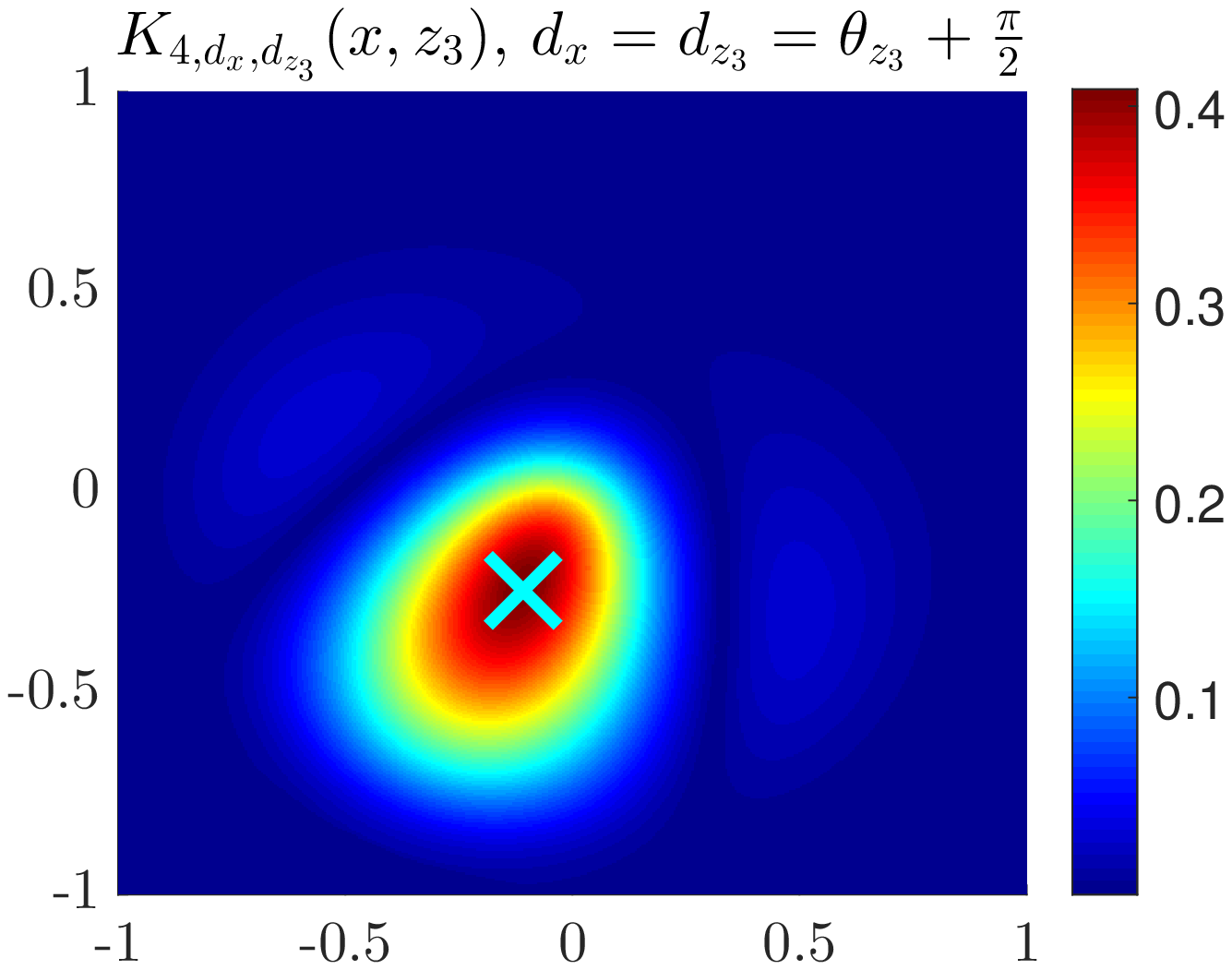}
    \caption{Almost orthogonality property of $K_{4,d_x,d_{z_i}}(x,z_i)$ for $V_0\neq 0$, with $m_1 = m_2 = 1/2$, $d_x = d_{z_i}$, and $z_1 = (-0.63,0.37)$, $z_2 = (-0.06,-0.73)$, $z_3 = (-0.11, -0.24)$ (from left to right).}
        \label{V_K4}
\end{figure}

\begin{figure}
\centering
    \includegraphics[scale = 0.35]{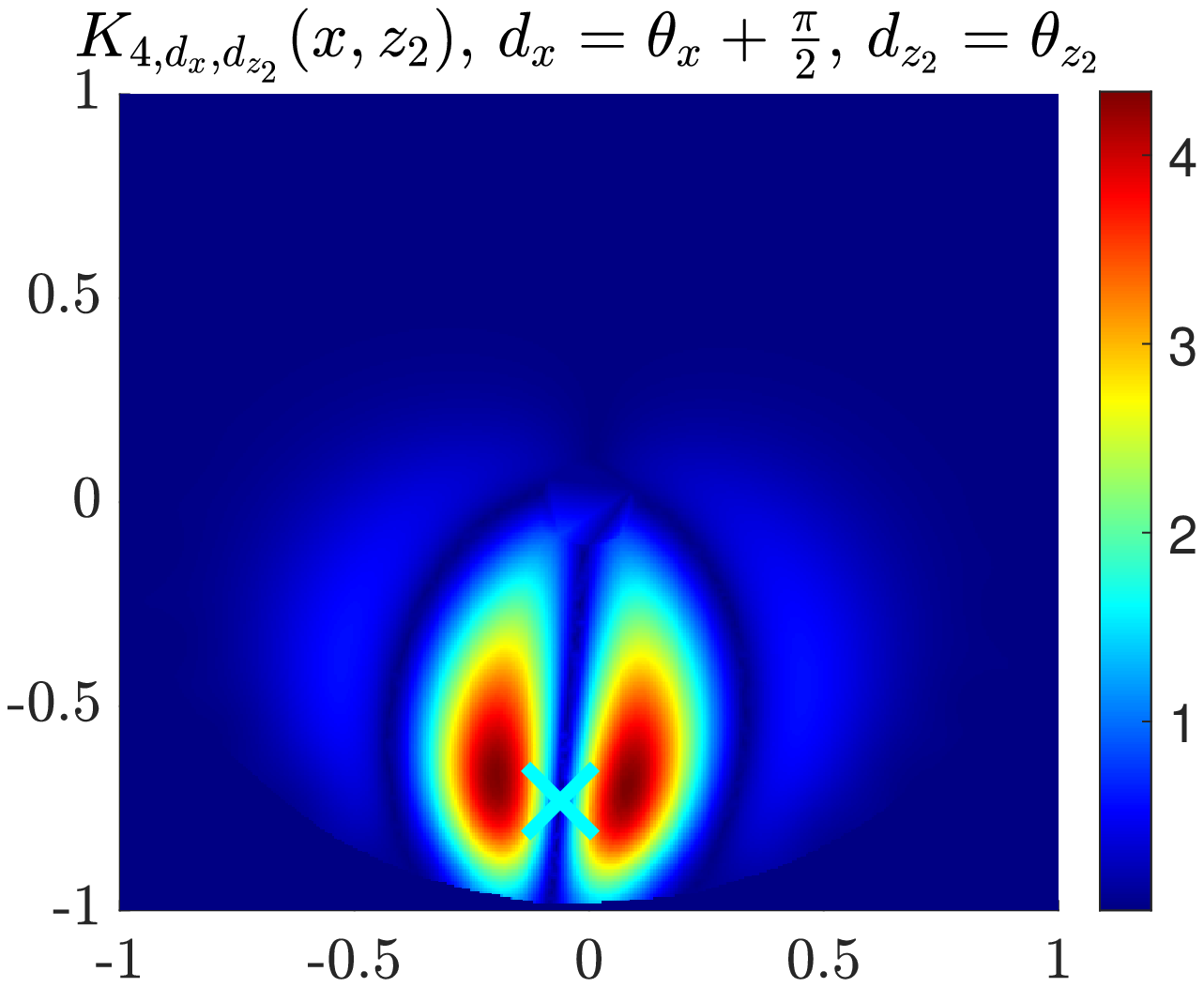}
    \includegraphics[scale = 0.35]{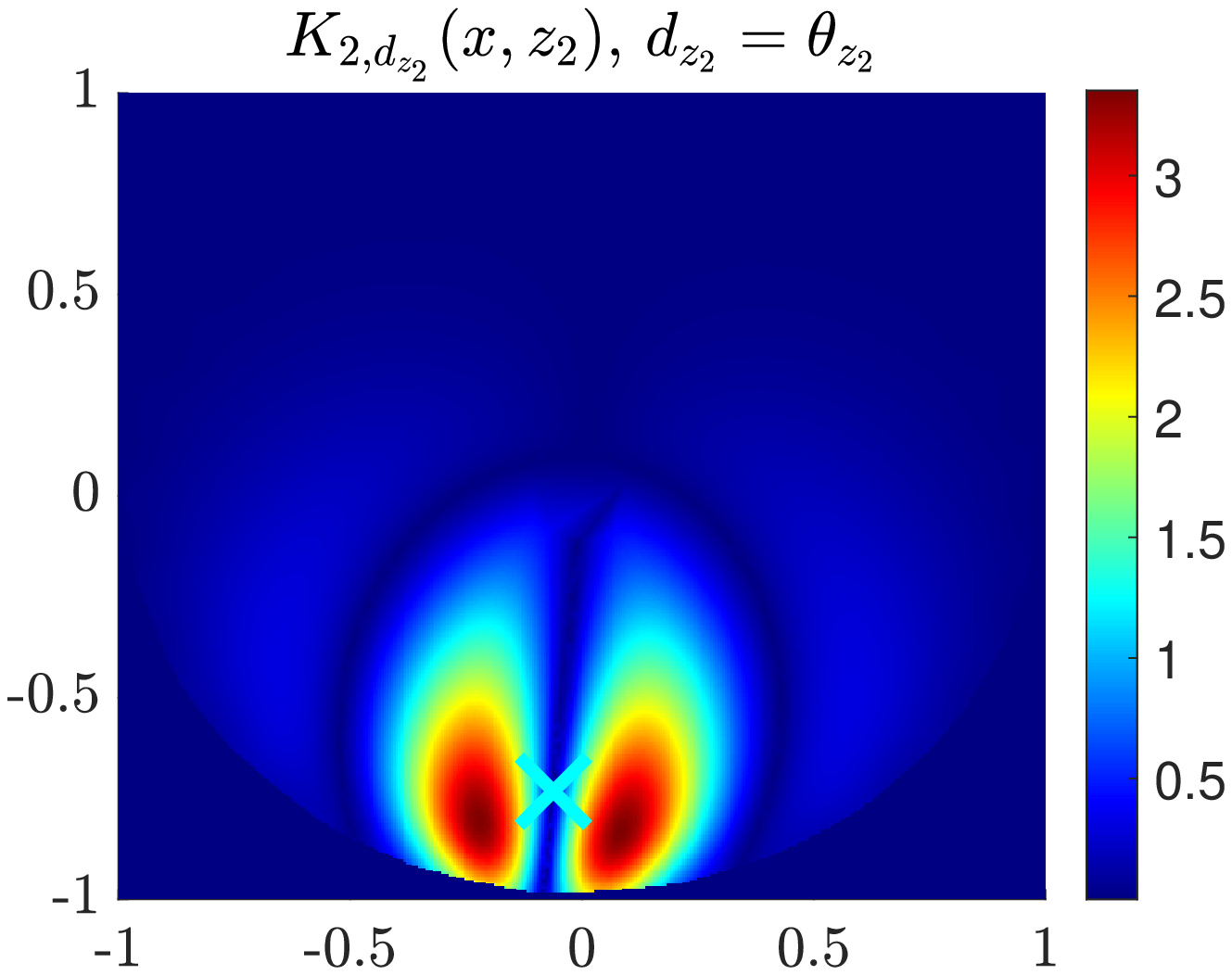}
     \includegraphics[scale = 0.35]{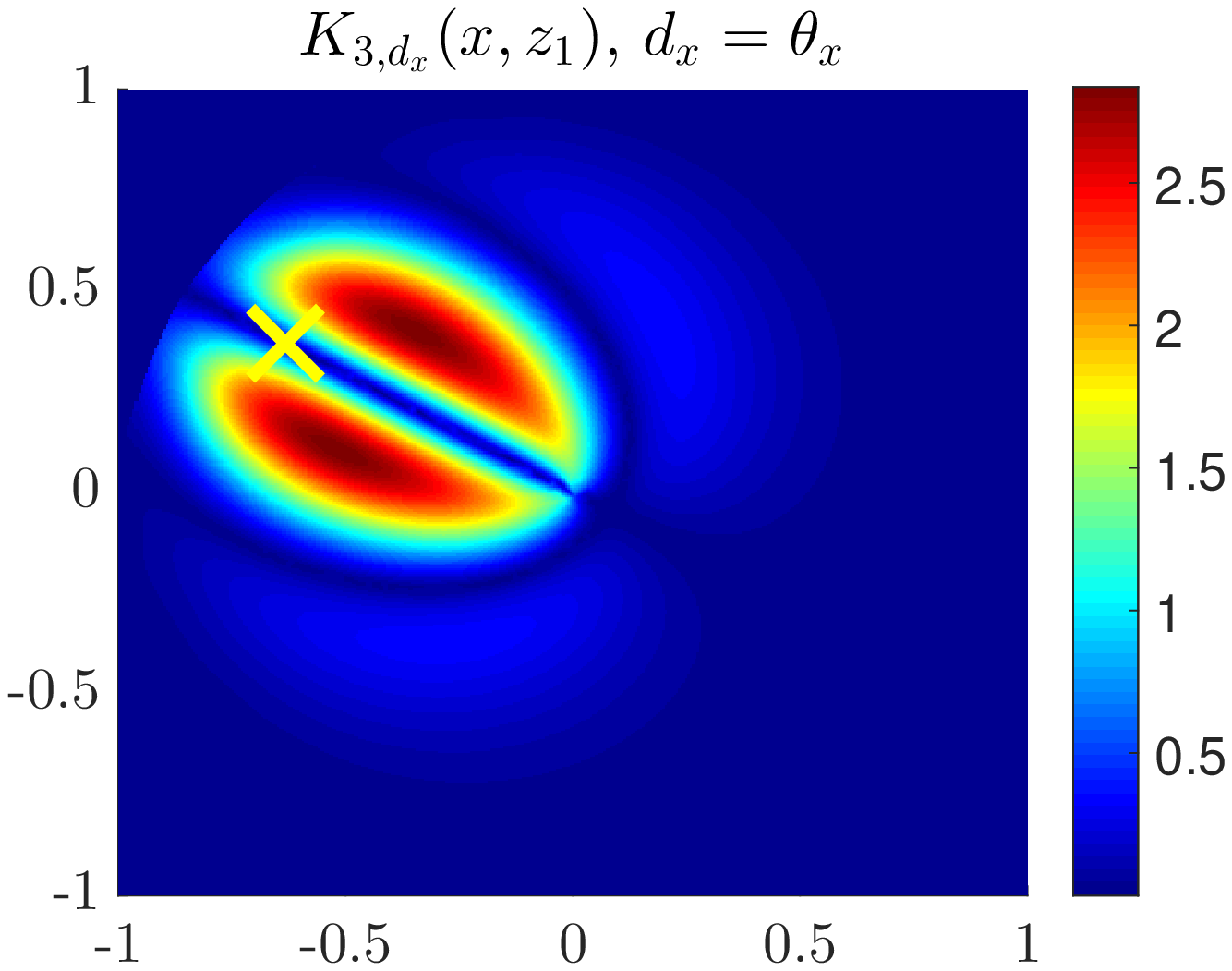}
    \caption{Mutually almost orthogonality property of $K_{4,d_x,d_{z_2}}(x,z_2)$, $K_{3,d_x}(x,z_1)$, and $K_{2,d_{z_2}}(x,z_2)$ for $V_0\neq 0$, with $m_i = n_i = 1/2$ ($i=1$, $2$), and $z_1 =(-0.63,0.37) $, $z_2 =  (-0.06,-0.73)$. Directions are chosen as $d_x\cdot d_{z_2}=0$, $d_{z_2} = \theta_{z_2}$, and $d_{x} = \theta_{x}$ (from left to right)}.
    \label{V_K23}
\end{figure}

\begin{figure}
\centering
    \includegraphics[scale = 0.35]{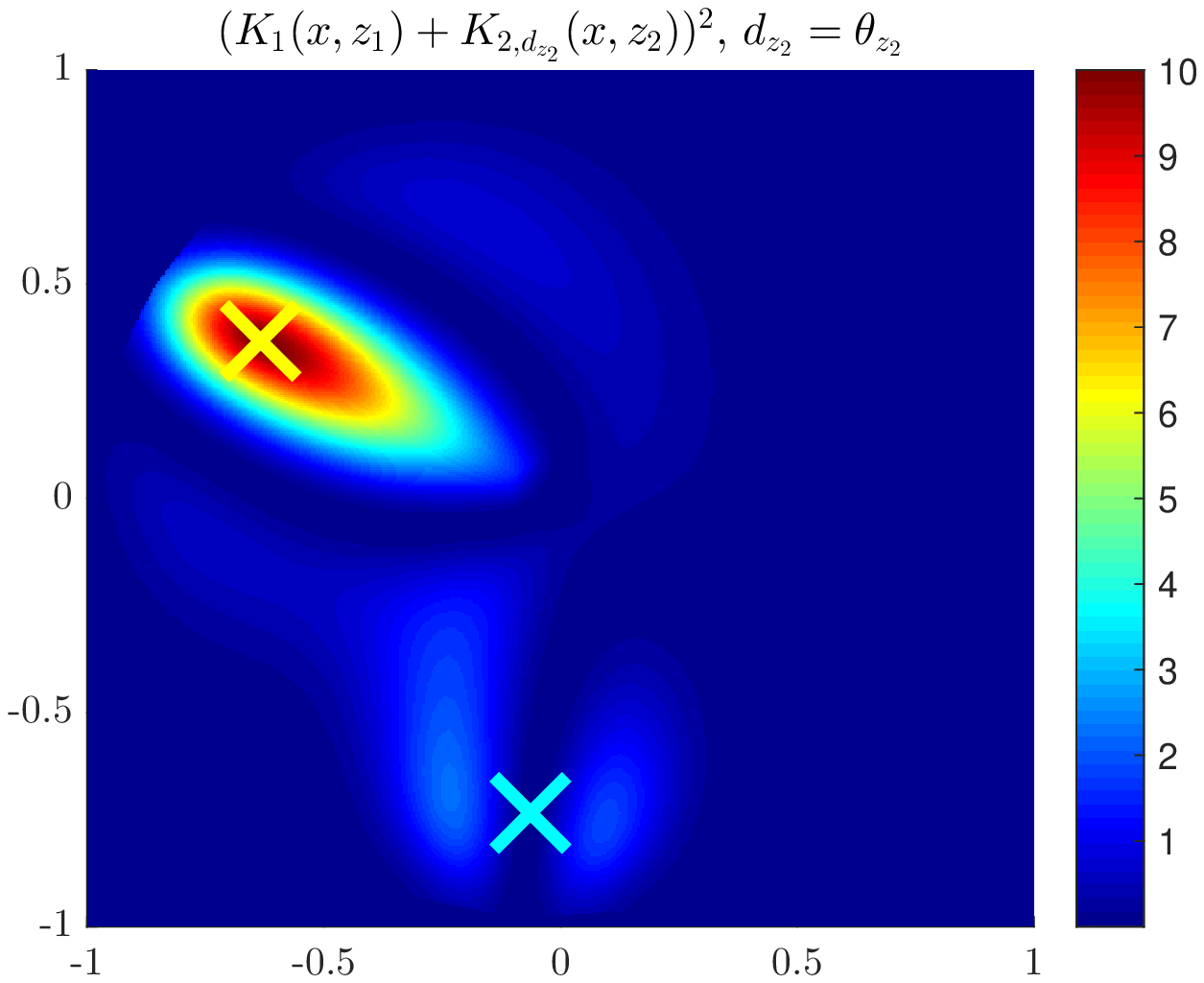}
    \includegraphics[scale = 0.35]{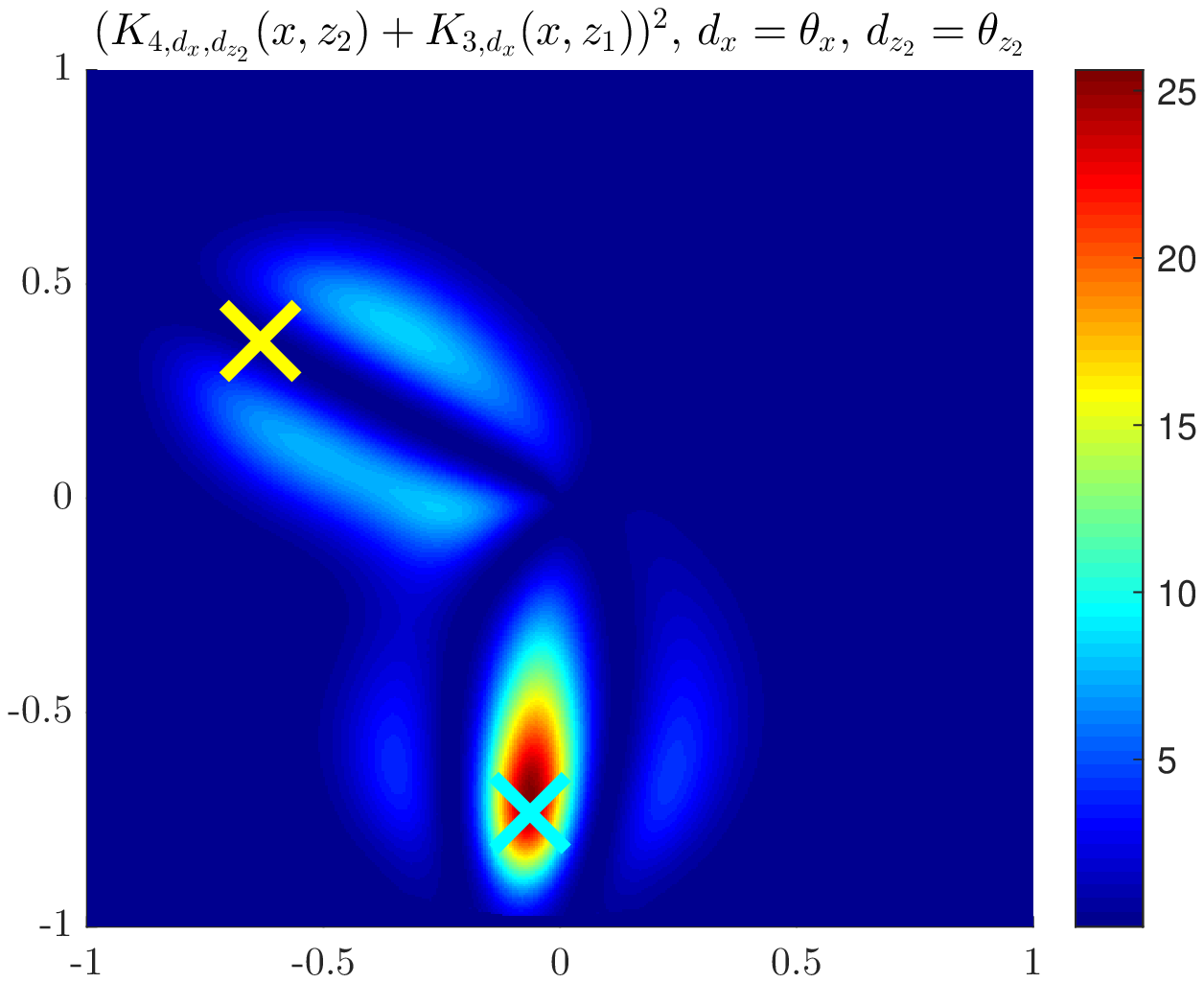}
 	\includegraphics[scale = 0.35]{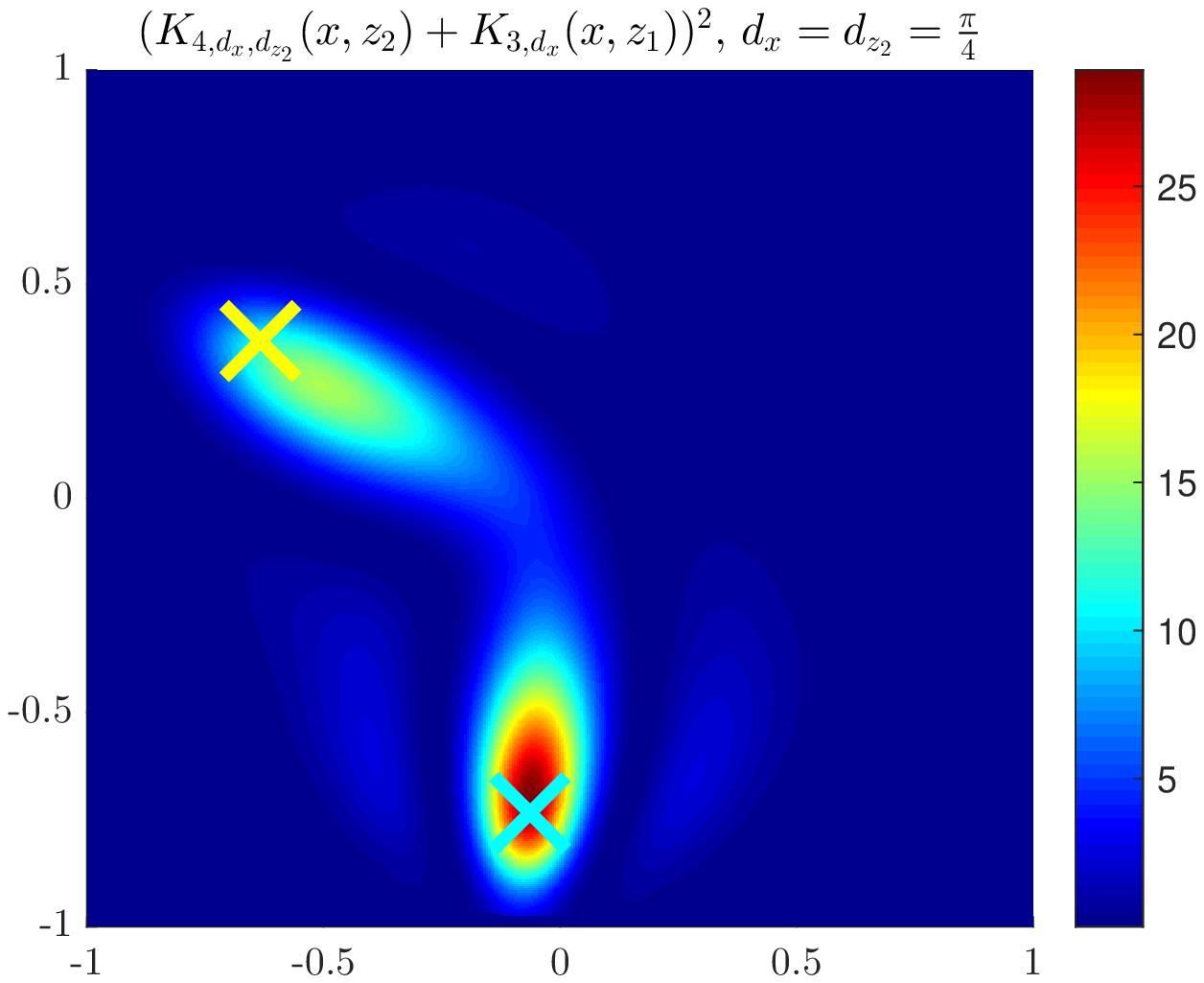}
    \caption{Mutually almost orthogonality property of $K_{1}(x,z_1)+K_{2,d_{z_2}}(x,z_2)$ (the left plot) and $K_{4,d_x,d_{z_2}}(x,z_2)+K_{3,d_x}(x,z_1)$ (the middle and the right plots) for $V_0\neq 0$, with $m_i = n_i = 1/2$ ($i=1$, $2$), and $z_1 = (-0.63,0.37)$, $z_2 = (-0.06,-0.73)$. Directions are chosen as $d_{z_2} = \theta_{z_2}$, $d_x = \theta_x$,  and $d_x = d_{z_2} = \pi/4$ (from left to right).}
    \label{V_K}
\end{figure}

\mm{Numerical experiments are conducted again to verify the mutually almost orthogonality property of 
the kernel functions in Figs.\,\ref{V_K1}-\ref{V_K},  with $k^2 =10$ and $R=1$. Three points are chosen 
in $\Omega$, i.e.,  $z_1 = (-0.63,0.37)$, $z_2 = (-0.06,-0.73)$, $z_3 = (-0.11,-0.24)$, 
and the constants $m_i = n_i = 1/2$ ($i=1$, $2$) are selected 
as the normalizations which are used in \eqref{eqn_defindexmo} and \eqref{eqn_defindexdi}. 
}
In the following figures, the yellow cross and the blue cross represent the location of a monopole and a dipole respectively.

\begin{enumerate}
	\item \mm{Fig.\,\ref{V_K1} plots the kernel $K_1(x, z_i)$ for $i = 1,2,3$. We can clearly see  its maximum 
	is attained when $x \approx z_i$, hence verifies the almost orthogonality property of $K_1(x, z_i)$.}
	\item \mm{Fig.\,\ref{V_K4} plots the kernel $K_{4,d_x,d_{z_i}}(x, z_i)$ for $i = 1,2,3$. With an appropriate probing direction,
	we can clearly see its maximum is attained when $x\approx z_i$ and $d_x= d_{z_i}$, hence verifies the almost orthogonality property of $K_{4,d_x,d_{z_i}}(x,z_i)$.}
	\item \mm{We show in Fig.\,\ref{V_K23} the effect of the probing direction.  
	In the first plot, we examine the special choice of the probing direction such that 
	$d_x\cdot d_{z_2} = 0$ at $z_2$, and see the kernel function 
	$K_{4,d_x,d_{z_2}}(x,z_2)$
	can not properly indicate the location of the dipole. The second and third plots demonstrate the behaviours of $K_{2,d_z}$ and $K_{3,d_x}$ when $d_z = \theta_z$, $d_x = \theta_x$. We notice that as in the case $V_0 = 0$, the peaks of the kernel functions appear to be very close to the location of the dipole or the monopole.  
	Meanwhile we see clearly that the value of $K_{4,d_x,d_z}$ is larger than the peak values of $K_{2,d_z}$ and $K_{3,d_x}$.
	}
	\item \mm{In Fig.\,\ref{V_K}, we examine the coexistence of a monopole at $z_1 = (-0.63,0.37)$ and a dipole at $z_2 = (-0.06,-0.73)$. To consider the case when the influence of the monopole and the dipole are comparable on the boundary, we enhance the strength of the monopole by multiplying a constant $1.5$.  
The first plot can be considered as probing by $\zeta_x$, while the second and third plots can be considered as probing by $\eta_{x,d_x}$ under different probing directions. We may conclude that the monopole probing function $\zeta_x$ 
	interacts better with the monopole located at $z_1$, while the dipole probing function $\eta_{x,d}$ 
	interacts better with the dipole located at $z_2$, under an appropriate probing direction.
	}
\end{enumerate}


%
%
\subsubsection{Explicit representations of probing functions in terms of Bessel function}\label{subsec_explicit}
Before we continue \mm{to explore} the mutually almost orthogonality property in other special domains, 
\mm{we present some explicit representations of the probing functions on the boundary of the unit disk.   
This will help us efficiently evaluate the inner products involved 
in the index functions \eqref{eqn_defindexmo} and \eqref{eqn_defindexdi}. 
Note that the corresponding norms of the probing functions used as the weights in the index functions 
were already given in the previous subsection.}

We first compute an explicit expression for $\zeta_x$.
Via a separation of variables, the solution to \eqref{eqn_defmonoknown2} can be represented by
\begin{equation}
    v_x^{(2)}(y) = \sum_{n=-\infty}^{\infty}C_n(k,r_x)I_n(kr_y)e^{in(\theta_y - \theta_x)}\,,
\end{equation}
where $x=(r_x,\theta_x)$, $y = (r_y,\theta_y)$ in polar coordinates, and $C_n(k,r_x)$ are coefficients determined by the boundary condition.
Now let us consider one special solution to \eqref{eqn_defmonoknown1}, which we may choose as $K_0(k |y-x|)$, where $K_0$ is the modified Bessel function of the second kind of order $0$. 
 Note that $x$ represents a point inside $\Omega$ and $y$ represents a point on $\partial \Omega$, hence we always have $r_y>r_x$. Applying the Graf's formula \cite{handbook}, we obtain 
\begin{equation}
    K_0(k|y-x|) =\sum_{n=-\infty}^{\infty}I_n(kr_x)K_n(kr_y)e^{in(\theta_y - \theta_x)}\,.
\end{equation}

Furthermore, we may determine $C_n(k,r_x)$ by a comparison of coefficients, and derive
\begin{equation}
\label{v_x(y)}
    v_x(y) = \sum_{n\in \mathbb{Z}}\bigg(I_n(kr_x)K_n(kr_y)-\frac{I_n(kr_x)K_n(k)}{I_n(k)}I_n(kr_y)\bigg)e^{in(\theta_y-\theta_x)}\,.
\end{equation}
Employing the relationship on the Wronskian between $K_n$ and $I_n$ \cite{handbook}, we then get 
the expression of $\zeta_z$ when $r_y = 1$:
\begin{equation}
    \zeta_x(y) = \frac{\partial v_x(y)}{\partial r_y} 
     = k\sum_{n\in\mathbb{Z}}\frac{I_n(kr_x)}{I_n(k)}e^{in(\theta_y-\theta_x)}\,.
\end{equation}

To compute $\eta_{x,d}$, we first note that $\eta_{x,d}$ is linear with respect to different choices of $d$, so it suffices to compute $\eta_{x,e_i}$ ($i=1,2$)  for two canonical basis vectors $e_1$ and $e_2$ in $\mathbb{R}^2$. For simplicity, we set
\begin{align}
    a_n(r_x,\,r_y) &= \frac{I_n(kr_x)}{I_n(k)}\big[I_n(k)K_n(kr_y)-K_n(k)I_n(kr_y)\big]\,,\\
    b_n(r_x,\,r_y) &= \,k\frac{I_n'(kr_x)}{I_n(k)}\big[I_n(k)K_n(kr_y)-K_n(k)I_n(kr_y)\big]\,.
\end{align}
A particular solution to $w_{x,e_1}$ defined in \eqref{eqn_defwxd} can be obtained by taking the partial derivative of $v_x(y)$ in \eqref{v_x(y)} with respect to $y\cdot e_1$: 
\begin{equation}\label{eq:wxe}
    w_{x,e_1}(y)  =\sum_{n \in \mathbb{Z}}\bigg[\cos(\theta_x)b_n(r_x,r_y)- in\frac{\sin(\theta_x)}{r_x}a_n(r_x,r_y) \bigg]e^{in(\theta_y-\theta_x)}\,.
\end{equation}
Then the probing function $\eta_{x,e_1}(y)$ in \eqref{eqn_defeta} with $r_y=1$
is obtained by applying the partial derivative with respect to $r_y$ 
\begin{equation}\label{eq:eta_xe}
    \eta_{x,e_1}(y) = \sum_{n \in \mathbb{Z}}\bigg[k\cos(\theta_x) \frac{I_n'(kr_x)}{I_n(k)}-in\frac{\sin(\theta_x)}{r_x}\frac{I_n(kr_x)}{I_n(k)} \bigg]e^{in(\theta_y-\theta_x)}\,.
\end{equation}
Similarly, $\eta_{x,e_2}$ can be given by
\begin{equation}
	 \eta_{x,e_1}(y) = \sum_{n \in \mathbb{Z}}\bigg[k\sin(\theta_x) \frac{I_n'(kr_x)}{I_n(k)}+in\frac{\cos(\theta_x)}{r_x}\frac{I_n(kr_x)}{I_n(k)} \bigg]e^{in(\theta_y-\theta_x)}\,.
\end{equation}

\subsection{\mm{Spherical domains in $\mathbb{R}^d$ for $d > 2$}}
\begin{figure}
\centering
    \centering
    \begin{subfigure}{2.5in}
    \includegraphics[scale = 0.4]{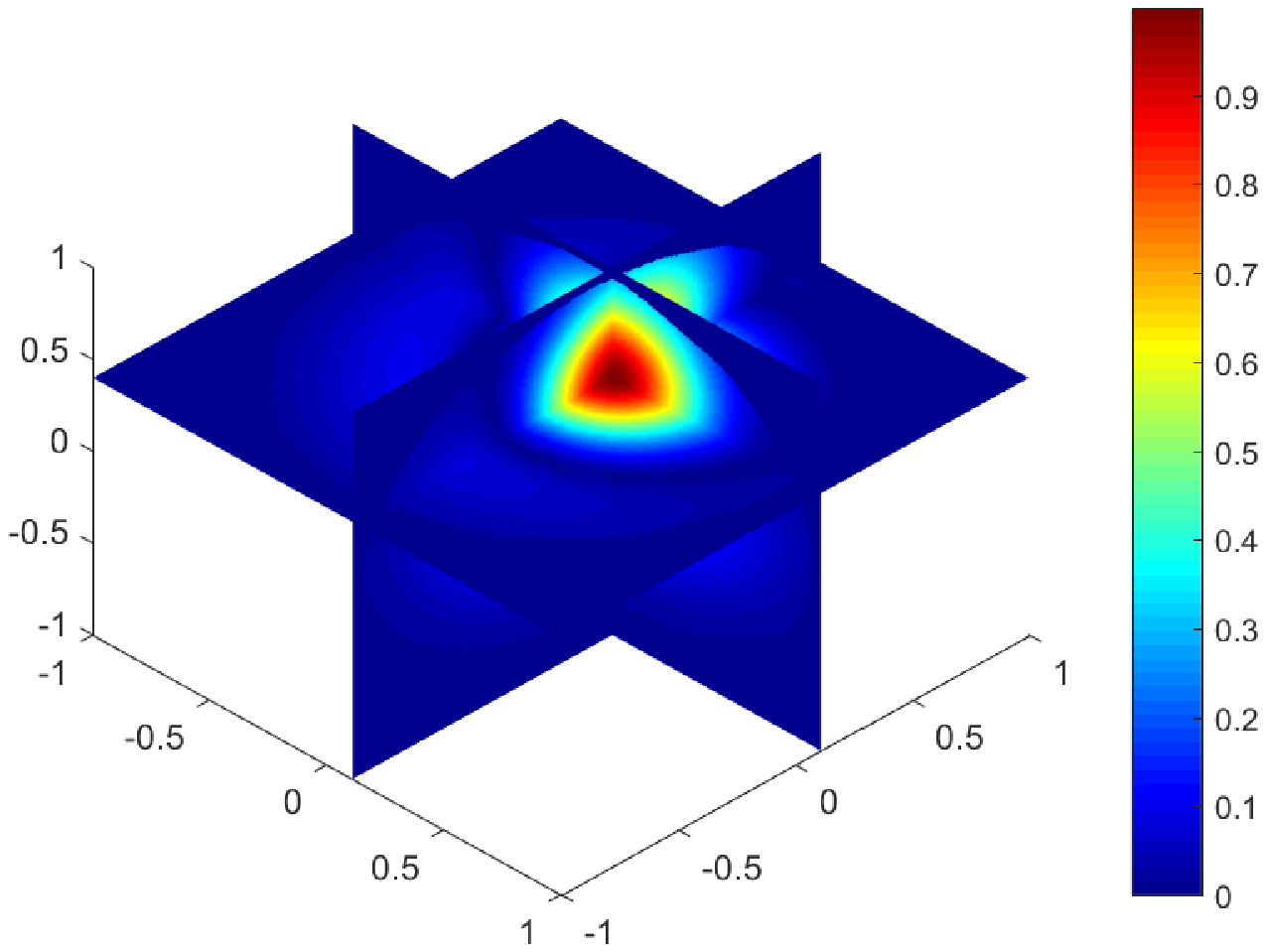}
    \centering
    \caption{$K_1$ defined in \eqref{kernal_mo}.}
    \end{subfigure}
    \begin{subfigure}{2.5in}
    \includegraphics[scale = 0.4]
        {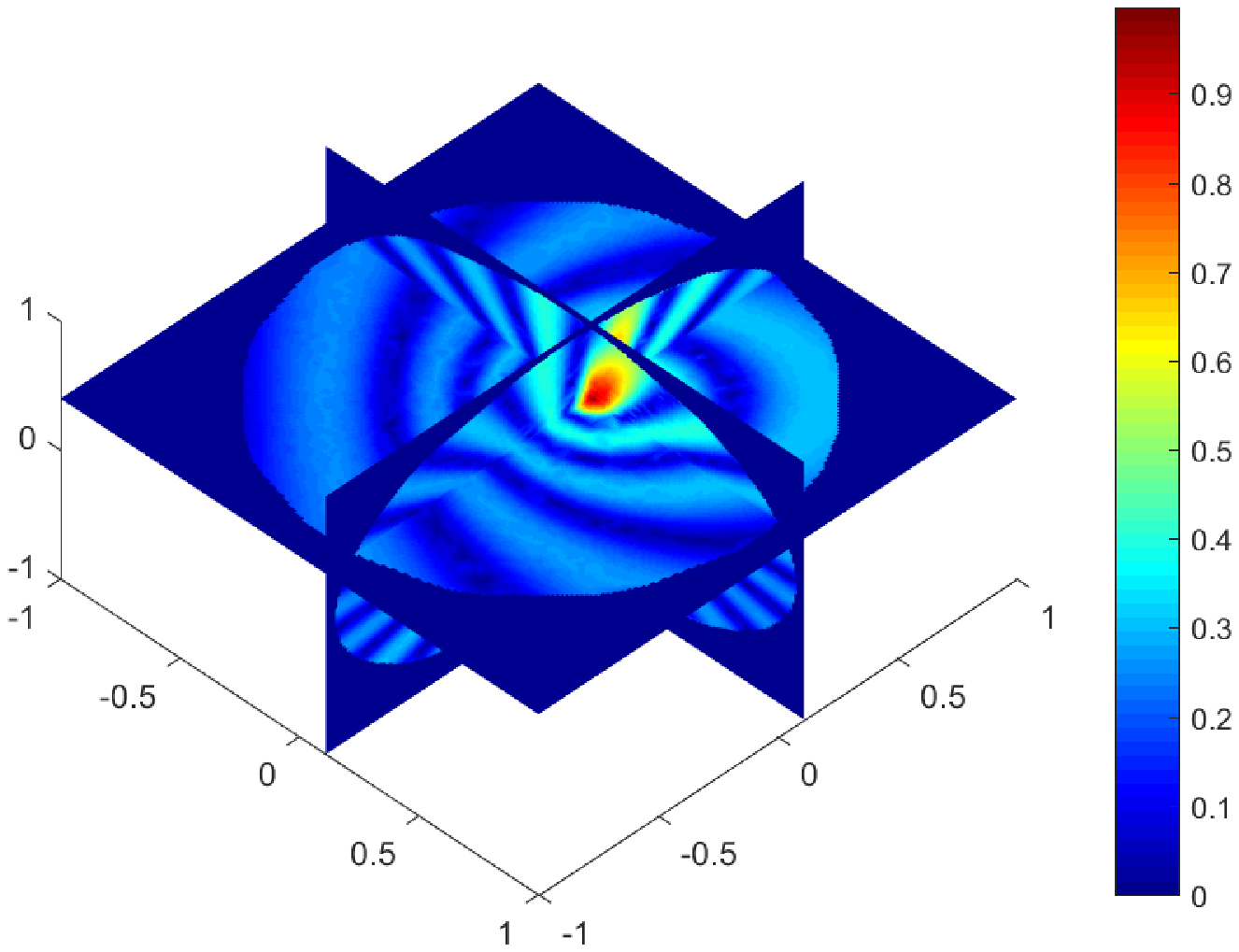}
            \centering
    \caption{$K_{4,d_{x},d_z}$ defined in \eqref{kernal_di}.}
    \end{subfigure}
    \caption{Almost orthogonality property of $K_1(x,z)$ and $K_{4,d_{x},d_z}(x,z)$ with $\gamma = 1$, $m_i = n_i = 1/2$ ($i=1$, $2$), $d_{x} = d_z = (0, 0, 1)$, $x = (0.114,0.114,0.396) $, and $z \in B(0,1)$.}
    \label{ver_sph}
\end{figure}

\mm{We now derive the explicit expressions of kernels $K_i$ defined in \eqref{kernal_mo} and \eqref{kernal_di} and 
the probing functions for the case of open balls in $\mathbb{R}^d$ for $d >2$. 
The analyses are quite similar to the circular case in the previous two subsections, so 
we will give a sketch only for
$d=3$ and emphasize some main differences.
}
Let $\Omega$ be a unit ball centered at $0$ in $\mathbb{R}^3$, 
and $\Gamma_n$ and $Y_n^{m} $ satisfy equations
\begin{equation}
        \label{sep_var_ball}\frac{r^2}{\Gamma_n}\frac{\partial^2 \Gamma_n}{\partial r^2}+\frac{2r}{\Gamma_n}\frac{\partial \Gamma_n}{\partial r}-(k^2r^2+ n(n+1))= 0\,; \quad -\Delta_{S^{2}}Y_n^{m} = n(n+1) Y_n^{m}\,.
    \end{equation}
Then by a separation of variables, the kernel of $- \Delta + k^2$ can be spanned by the Schauder basis 
$\{ \Gamma_n (r)  Y_n^{m} (\theta, \phi)\,,$ $~{n \in \mathbb{N}\,,~ |m|\leq n}\}$. And we can readily check that $\Gamma_n$ can be solved by the spherical Bessel function of the first kind $j_n$ while $Y_n^m$ can be solved by the spherical harmonic function.
\mm{The eigenpairs defined in \eqref{eqn_defNtD} for $d=3$ can be given by
\begin{equation}
        \varphi_n^m = \frac{j_n(ikr)}{j_n(ik)}Y_n^m(\theta, \omega)\,, ~~
         \lambda_n = \frac{j_n(ik)}{ik j_n'(ik)}\,, ~~
 n \in \mathbb{N}\,,  ~m = -n,\dots,n\,.
\end{equation}

Since the spherical harmonics form a complete orthogonal basis in $L^2(\mathbb{S}^2)$, we may rewrite 
the duality product, the $H^{\gamma} $ semi-norm, and probing functions in terms of this basis.  For instance, 
we can write the $H^\gamma$ duality product as
}
\begin{equation}
\label{3d_Hgamma}
    \langle f, g\rangle_{H^\gamma} = \sum_{n\in \mathbb{N}} \sum_{m = -n}^n n^\gamma(n+1)^{\gamma}\overline{\hat{f}(n,m)}{\hat{g}(n,m)}\,,
\end{equation}
\mm{where $\hat{f}(n,m) = \int_{\mathbb{S}^2}f(\theta, \omega)Y_n^m(\theta, \omega)ds$ is the corresponding coefficient. 
Then using the addition formula for Legendre polynomials, 
we can obtain all we need for an explicit expression of $K_1$ (with $\gamma = 1$):
}
\begin{equation*}
\begin{split}
    &\langle \zeta_{x_1}, G_{x_2}\rangle_{H^1}
    = \sum_{n\in \mathbb{N}}\frac{n(n+1)(2n+1)^2I_{n+\frac{1}{2}}(kr_1)I_{n+\frac{1}{2}}(kr_2)P_n(\frac{x_1\cdot x_2}{r_1r_2})}{4\pi k I_{n+\frac{1}{2}}(k)(r_1r_2)^{1/2}[nI_{n-\frac{1}{2}}(k)+(n+1)I_{n+\frac{3}{2}}(k)]}\,;
    \end{split}
\end{equation*}
\begin{equation*}
    |\zeta_{x_1}|_{H^1}^2 = \sum_{n\in \mathbb{N}} \frac{(n)(n+1)(2n+1)(I_{n+\frac{1}{2}}(kr_1))^2}{ 4\pi r_1 (I_{n+\frac{1}{2}}(k))^2}\,, \qquad 
    |G_{x_1}|_{H^1}^2 = \sum_{n\in \mathbb{N}} \frac{(n)(n+1)(2n+1)^3(I_{n+\frac{1}{2}}(kr_1))^2}{4\pi k^2r_1 [nI_{n-\frac{1}{2}}(k)+(n+1)I_{n+\frac{3}{2}}(k)]^2}\,.
    \end{equation*} 
The explicit expressions for $K_{2,d_z}, K_{3,d_x}, K_{4,d_x,d_z}$, as well as that of the probing functions, are similar.   
As an example, Fig.\,\ref{ver_sph} shows the almost orthogonality property for the kernel $K_1(x,z)$ and $K_{4,d_{x},d_z} (x,z)$ defined in \eqref{kernal_mo} and \eqref{kernal_di}, 
with $\gamma = 1$, $m_i = n_i = 1/2$, ($i = 1$, $2$), $d_x = d_z = (0, 0, 1)$, $x = (0.114,0.114,0.396)$ and $z \in \Omega$.



\subsection{A decoupling strategy based on the frequency of the boundary influx}
\label{sec_decoupling}
\mm{In this subsection, we investigate a decoupling strategy that makes use of the effect from changing the frequency of the boundary influx. This strategy is a very reliable and effective decoupling technique when we implement our DSM.
For illustrations, we consider two different cases: the first one for two small inhomogeneous inclusions, 
each inhomogeneity from one of two parameters $\sigma$ and $V$ in \eqref{eqn_prostso}; 
the second one for one inhomogeneous inclusion.
}

\subsubsection{\mm{Two small inhomogeneous inclusions}}
\mm{Let us consider a simplified situation when there are two small inhomogeneous inclusions $D_1$, $D_2$ in $\Omega = B_1$.  
We write $D_1 = z_1 + \delta B_1$, $D_2 = z_2 + \delta B_1$, with $z_1$, $z_2 \in \Omega$ and $|\delta|<<1$.
We further assume in \eqref{eqn_prostso} that 
$\sigma = \sigma_1$ in $D_1$ and $\sigma = \sigma_0$ otherwise; and $V = V_1$ in $D_2$ and $V = V_0$ otherwise.  
Under this setting, we can readily obtain the asymptotic expansion of $u-u_0$ for $x \in \partial \Omega$, 
uniformly as $k\delta \rightarrow 0$ \cite{hansen2008high}:
}
\begin{equation}
\label{asy_exp}
\begin{split}
        (u-u_0)(x)  \approx \delta^2 \Big\{C_1(\sigma,\sigma_0,\Omega)\nabla G_{z_1}(x)\cdot \nabla u_0(z_1) +  C_2(V,V_0,\Omega)G_{z_2}(x)u_0(z_2)\Big\}\,,
\end{split}
\end{equation}
where constants $C_1$ and $C_2$ depend only on \mm{the domain}. Supposing the boundary influx is of the form $f = e^{im\theta}$ on $\partial \Omega$, \mm{we can get} the following expressions of $u_0$ satisfying \eqref{eqn_prostsh} and its gradient:
\begin{equation*}
    u_0(x) = \frac{I_m(kr_x)}{I_m'(k) k}e^{im\theta_x}\,, \qquad \frac{\partial  u_0(x)}{\partial r} = \frac{I_m'(kr_x)}{I_m'(k)}e^{im\theta_x}\,, \qquad \frac{\partial u_0(x)}{\partial \theta} = imu_0(x)\,;
\end{equation*}
\begin{equation*}
\begin{split}
        &\nabla u_0(z_1)
     = \begin{pmatrix}
        \cos(\theta_{z_1}) & -\sin(\theta_{z_1})\\
        \sin(\theta_{z_1}) & \cos(\theta_{z_1})
    \end{pmatrix} \binom{I_m'(k r_{z_1})}{\frac{imI_m(kr_{z_1})}{k r_{z_1}}}\frac{e^{im\theta_{z_1}}}{I_m'(k)} = \sqrt{\bigg(\frac{I_m'(kr_{z_1})}{I_m'(k)}\bigg)^2+\bigg(\frac{mI_m(kr_{z_1})}{kr_{z_1}I_m'(k)}\bigg)^2} \vec{d}_{z_1}\,,
    \end{split}
\end{equation*}
where $|\vec{d}_{z_1}| = 1$. Denoting $\tilde{\beta}_{m}(z_1)  = \Big\{\Big(\frac{I_m'(kr_{z_1})}{I_m'(k)}\Big)^2+\Big(\frac{mI_m(kr_{z_1})}{kr_{z_1}I_m'(k)}\Big)^2\Big\}^{1/2}$, $\beta_m(z_2) = \frac{I_m(kr_{z_2})}{I'_m(k)k}$, 
we can readily derive
\begin{equation}
\label{ratio_ratio}
    \frac{|\nabla u_0(z_1)|}{| u_0(z_2)|} = \frac{\tilde{\beta}_m(z_1)}{\beta_m(z_2)} = \sqrt{\bigg(\frac{kI_m'(kr_{z_1})}{I_m(kr_{z_2})}\bigg)^2 + \bigg(\frac{m}{r_{z_1}}\frac{I_m(kr_{z_1})}{I_m(kr_{z_2})} \bigg)^2 }\geq \frac{m}{r_{z_1}}\frac{I_m(kr_{z_1})}{I_m(kr_{z_2)}}\,.
\end{equation}
\mm{The above comparison hints that the inhomogeneity associated with $\sigma$ is more sensitive to the change of frequency around the local maxima of $K_1$, $K_{2,d_x}$, $K_{3,d_z}$, $K_{4,d_x,d_z}$ when $r_{z_1}\approx r_{z_2}$. To see this, 
let us consider the index function in \eqref{eqn_defindexdi} when Sobolev scale $\gamma = 0$, then we can approximate $I_{di}$ in \eqref{kernal_di} by 
}
\begin{equation}
\begin{split}
        I_{di}(x, d_x)  \approx & \,\,C_1\tilde{\beta}_m(z_1)K_{4,d_x, d_{z_1}} (x,z_1) + C_2 \beta_m(z_2)e^{im\theta_{z_2}}K_{3, d_x }(x,z_2)\,.
        \end{split}
\end{equation}
Now from \eqref{ratio_ratio}, \mm{it is ready to see} that the coefficient associated with $K_{4,d_x,d_z}$ will be more significant as $m$ becomes larger compared with the coefficient associated with $K_{3,d_x}$. Therefore, \mm{we should expect 
a much larger value of the index function around $D_1$} when the boundary influx has a higher frequency.


\subsubsection{\mm{A single inhomogeneous extended inclusion}}
\mm{We now consider the case when there is a single inhomogenous inclusion that is not necessarily small.
We compare the effects of varying two inhomogeneous coefficients $\sigma_1$ and $V_1$ in the same inclusion.
For the sake of exposition, we assume that the inhomogeneity is located in a disk $B_R$ 
with radius $R$, 
and take $u_0 = {I_m(kr)e^{im\theta}}/{I_m(k)}$ 
in polar coordinates. 
}

\mm{\textbf{Case 1:} $V$ is constant, but $\sigma$ is piecewise constant, i.e., 
$\sigma = \sigma_1$ in $B_R$, and $\sigma = \sigma_0$ otherwise.
Letting $k_s^2 := {V_0}/{\sigma_1}$, then the scattered wave $u^s := u - u_0$ and the total wave $u$ satisfy the equations
}
\begin{equation}
\begin{cases}
        - \Delta u+ k_s^2 u = 0 \qquad &|x|<R\,, \\
        - \Delta u^{s} + k^2 u^{s} = 0 \qquad &|x|>R\,, \\
        u^s + u_0 = u \qquad &\text{on} \quad \partial B_R\,, \\
        \sigma_0 \frac{\partial (u^s + u_0) }{\partial \nu} = \sigma \frac{\partial u}{\partial \nu} \qquad &\text{on} \quad \partial B_R\,.
\end{cases}
\end{equation}
\mm{As we expect no singularity for $u$ around the origin, we may assume $u (r,\theta) = \sum_{n=1}^{\infty} \alpha_n I_n(k_s r)e^{in\theta}$ for some $\alpha_n$. Similarly, we write 
$u^s(r,\theta)= \sum_{n=1}^{\infty} \beta_n K_n(k r)e^{in\theta}$ for some $\beta_n$. 
By comparing Fourier coefficients, we easily see 
$\alpha_n = \beta_n = 0$ if $n \neq m$. Therefore it suffices to consider the Fourier coefficient 
associated with $e^{im\theta}$. Using the transmission condition on $\partial B_R$, we derive
\begin{equation}
\label{bm}
\begin{split}
    |\beta_m| =& \,\,\bigg|\frac{{k_s}I_m(k_sR)I_m'(kR) - k I_m'(k_sR)I_m(kR) }{k I_m'(k_sR)K_m(kR) -  k_s I_m(k_sR)K_m'(kR)}\bigg|\frac{1}{I_m(k)}
    \geq  C  \bigg(\frac{I_m(k_sR)I_m(kR)}{I_{m}(k_sR)K_{m+1}(kR)I_m(k)}\bigg)\,,
    \end{split}
\end{equation}
}
\mm{
for some constant $C>0$, where we have used the following estimate for Bessel functions \cite{handbook}:
\begin{eqnarray}
         \bigg|k_sI_m(k_sR)I_{m+1}(kR) - kI_{m+1}(k_sR)I_m(kR)\bigg| 
 &= &\bigg|\bigg[I_m(kR)I_{m}(k_sR)kk_s\bigg] \bigg[\frac{I_{m+1}(kR)}{kI_m(kR)} - \frac{I_{m+1}(k_sR)}{k_sI_m(k_sR)}\bigg]\bigg|
 \nonumber \\
 &\leq & \,\bigg( I_m(kR)I_{m}(k_sR)kk_s\bigg)\bigg(\frac{R}{m}\bigg)\,.
\end{eqnarray}
}
%

\textbf{Case 2:} 
\mm{$\sigma$ is constant, and $V$ is piecewise constant, i.e., 
$V = V_1$ in $B_R$, and $V = V_0$ otherwise. 
Letting $k_v^2 := {V_1}/{\sigma_0}$, we write the scattered wave $\tilde{u}^s(r,\theta)= \sum_{n=1}^{\infty} \tilde{\beta}_n K_n(k r)e^{in\theta}$ for some $\tilde{\beta}_n$.  Again, we can see that $\tilde{\beta}_n = 0$ for $n \neq m$, 
hence we need to focus only on $\tilde{\beta}_m$, which can be estimated as follows: 
\begin{equation}
\label{tbm}
\begin{split}
        \bigg|\tilde{\beta}_m\bigg| =&\,\,\bigg|\frac{kI_m(k_vr)I_m'(kr) -k_v I_m'(k_vr)I_m(kr)}{ k_v I_m'(k_vr)K_m(kr) - kI_m(k_vr)K_m'(kr)}\bigg|\frac{1}{I_m(k)}
        \leq \tilde{C}\bigg(\frac{I_m(k_vR)I_m(kR)}{mI_{m}(k_vR)K_{m+1}(kR)I_m(k)}\bigg)\,.
        \end{split}
\end{equation}
}
\mm{
\textbf{Comparison between Cases 1 and 2:}
Considering the ratio $\tau_m :={|\beta_m|}/{|\tilde{\beta}_m|}$ between the Fourier coefficients from the above two cases, 
we can readily see from \eqref{bm} and \eqref{tbm} that $\tau_m\geq c \, m$ for some constant $c$. 
Noting that $\beta_m$ and $\tilde{\beta}_m$ represent the magnitude of the scattered waves 
for two different inhomogeneous inclusions respectively, 
this infers that the measurement coming from the inhomogeneous inclusion with a different $\sigma$ 
is more sensitive than that coming from an inhomogeneous inclusion with a different $V$ at the high frequency regime of the boundary influx.  
}

\section{Numerical experiments}
\label{sec_numerical}

\mm{In this section, we present a series of typical examples to illustrate the efficiency and robustness of 
our proposed direct sampling method for solving the inverse coefficient problem \eqref{eqn_prostso}.
We take the probing domain $\Omega$ to be the unit disk in $\mathbb{R}^2$, and the coefficients 
$\sigma_0$ and $V_0$ in the homogeneous background to be $\sigma_0 = 1$, $V_0 = 10$. 
For each numerical experiment, there are several inhomogeneities of different types that are located 
separately inside the domain.}

{\bf Forward data}. 
In all the experiments, we choose a boundary influx $f= \cos(k\theta)$, with different $k \in \mathbb{N}$. 
We solve the forward problem for $u$ and $u_0$ using a finite element method of mesh size $1/100$, 
and take as the forward data the values of the potential 
$u_s= u - u_0$ at a set of discrete probing points, denoted by $\Gamma_p$, 
distributed uniformly on the boundary of $\Omega$. 
Then the noisy data is generated by adding a random noise of multiplicative form: 
\begin{equation}
    u_s^{\delta}(x) = u_s(x)(1+\varepsilon \delta)\,, \quad x\in \Gamma_p\,,
\end{equation}
where $\varepsilon$ is randomly uniformly distributed in $[-1,1]$. 
Unless it is specified otherwise, $\Gamma_p$ shall often consist of $48$ points, 
and the noise level $\delta$ is chosen to be $3\%$. 

Then we move on to address the implementation of the new DSM. 
We first compute the pointwise evaluations of the monopole and dipole probing functions using the explicit expressions 
in section \ref{subsec_explicit},  and all these are carried out off-line. 
We then compute the monopole and dipole index functions 
$I_{\text{mo}} (x)$ and $I_{\text{di}} (x,d_x)$ in \eqref{eqn_defindexmo} and \eqref{eqn_defindexdi} at each sampling point 
through appropriate numerical integrations. 
In all our numerical examples, we choose the parameters involved 
in \eqref{eqn_defindexmo} and \eqref{eqn_defindexdi} as follows:
$n_1 = n_2 = 1/2$, $m_1 = m_2 = 1/2$, $\gamma_{\text{mo}} = \gamma_{\text{di}} = 1$ (except \textbf{Example 1}).  
At each probing point $x$, the probing direction $d_x$ is chosen to be $d_x = \nabla \phi(x) / | \nabla \phi(x) |$,  
as it is described in section \ref{subsec_alternative}.    

We make a remark on the denominator of $I_{mo}$, by noting the fact 
that $|\zeta_{\vec{0}}|_{H^1} = 0$ from \eqref{zetanorm_V} and hence 
the index function $I_{mo}$ is singular around the origin when $\gamma = 1$. 
To get rid of this singularity, 
we take $|\zeta_{x_1}|_{H^1} = |\zeta_{(\eta,0)}|_{H^1}$ for all $|x_1| < \eta$, 
with $\eta$ fixed at $0.1$. 
The same modification is also applied to $|G_{x_2}|_{H^1}$.

For each example, 
we plot the exact inhomogeneous inclusions, along with 
the monopole and dipole index functions $\tilde{I}_{\text{mo}}$ and 
$\tilde{I}_{\text{di}}$, which are the squares of the respective normalized 
monopole and dipole index functions 
$I_{\text{mo}} (x) / \max_{y} I_{\text{mo}}(y)$ and $ I_{\text{di}} (x,d_x) / \max_{y} I_{\text{di}}(y, d_y)$.
The choice of squaring the index functions and normalizing by their maximum are only for the sake of 
better illustrations, and other choices can be used as well. 
In all the figures showing the exact inclusions, the orange color represents an inhomogeneity associated with 
$\sigma$, whereas the blue color represents an inhomogeneity associated with $V$. 

\subsection{Numerical tests on appropriate choices of boundary influxes and Sobolev index}
\label{subsec_numerical_ver}
%

We start first with an illustrative example  
to demonstrate the effectiveness of the decoupling strategy we proposed in 
section\,\ref{sec_decoupling} 
for choosing boundary influxes $f$ with different frequencies and the necessity of 
choosing a non-zero Sobolev scale $\gamma$ 
that appears in the index functions \eqref{eqn_defindexmo} and \eqref{eqn_defindexdi}.
We pick us a toy example, Example 1, 
that contains two inhomogeneous inclusions, arising from $\sigma$ and $V$, respectively.
With boundary influxes of different frequencies, we compare the indices $\tilde{I}_{\text{mo}}$ and 
$\tilde{I}_{\text{di}}$.
%
This helps us develop an appropriate choice of two frequencies for boundary influxes 
for the use in all the subsequent evaluations of the monopole and the dipole index functions. 

\begin{figure}
    \centering
    \begin{subfigure}{3in}
\begin{center}
    \includegraphics[scale = 0.25]{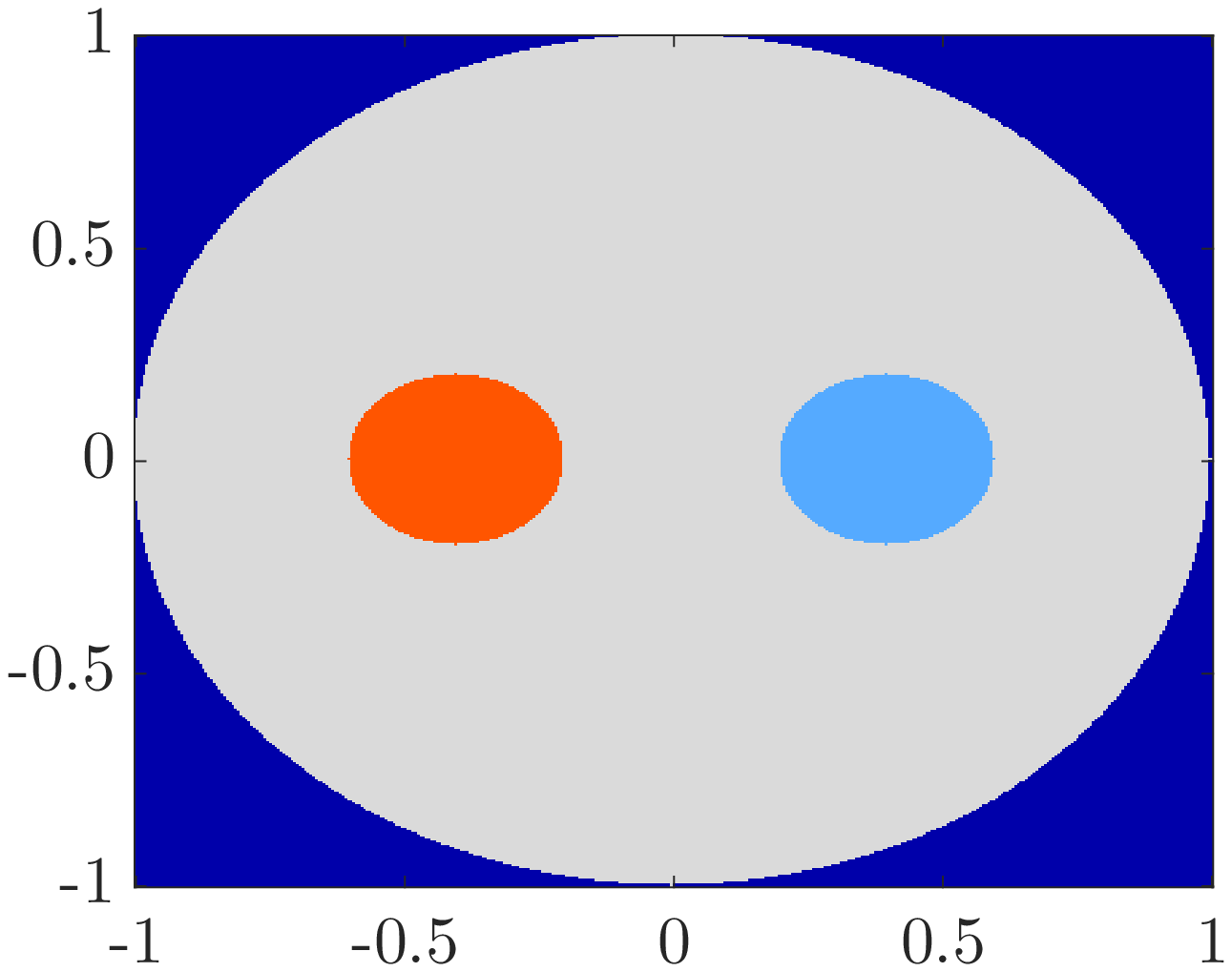}
\end{center}
    \end{subfigure}
    \begin{subfigure}{3in}
    \includegraphics[scale = 0.25]{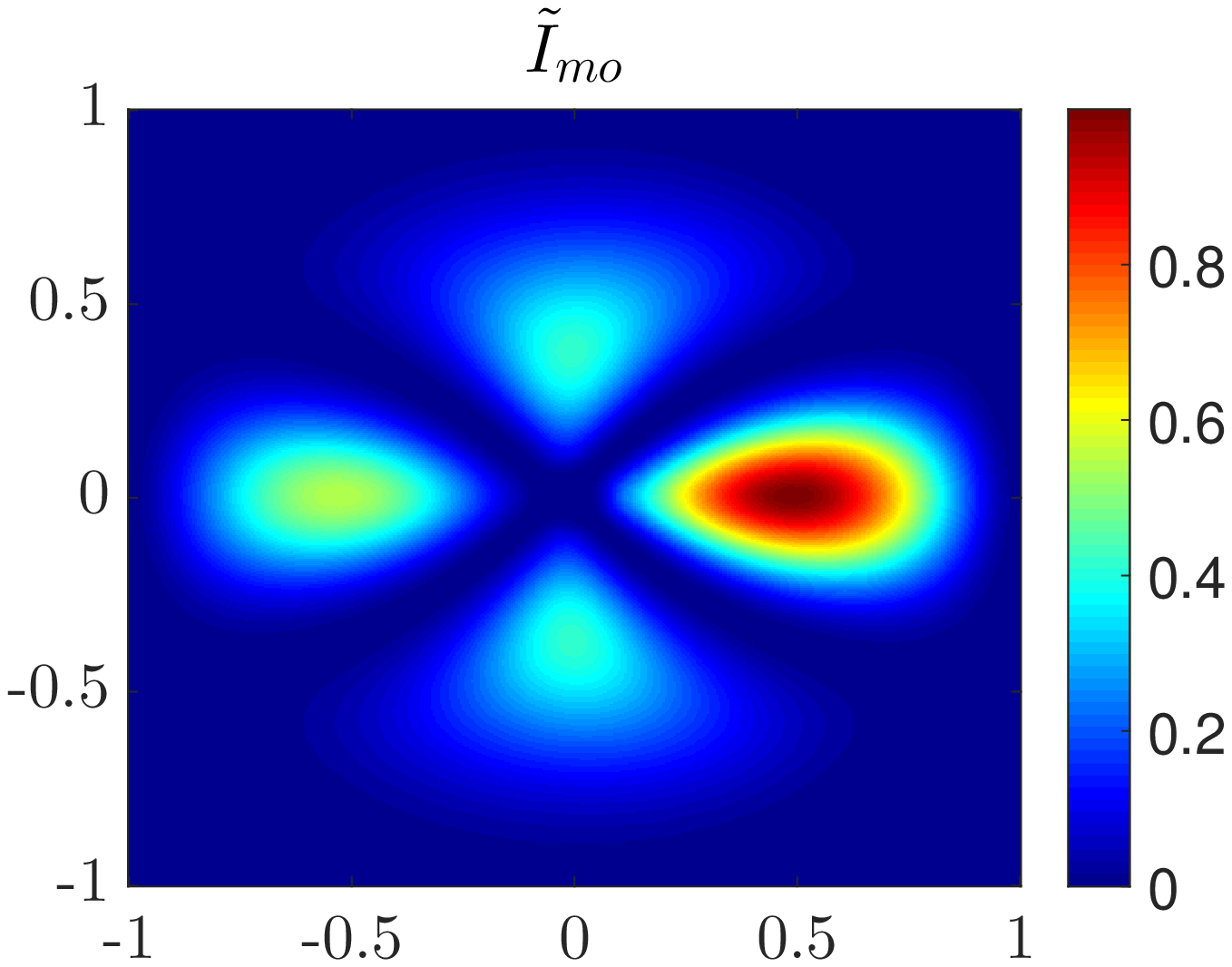}
    \includegraphics[scale = 0.25]{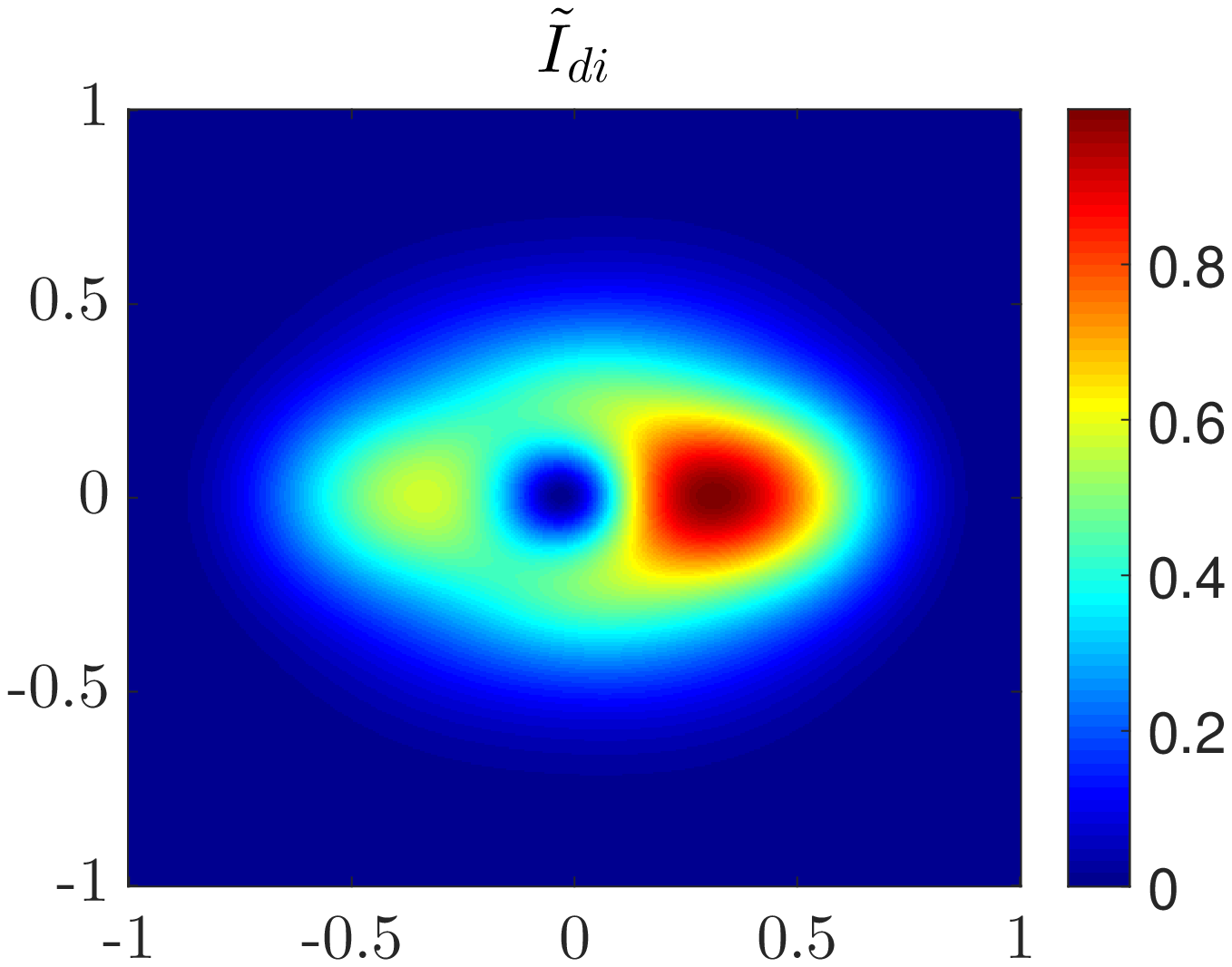}    
    \end{subfigure} \\
    \begin{subfigure}{3in}
    \includegraphics[scale = 0.25]
        {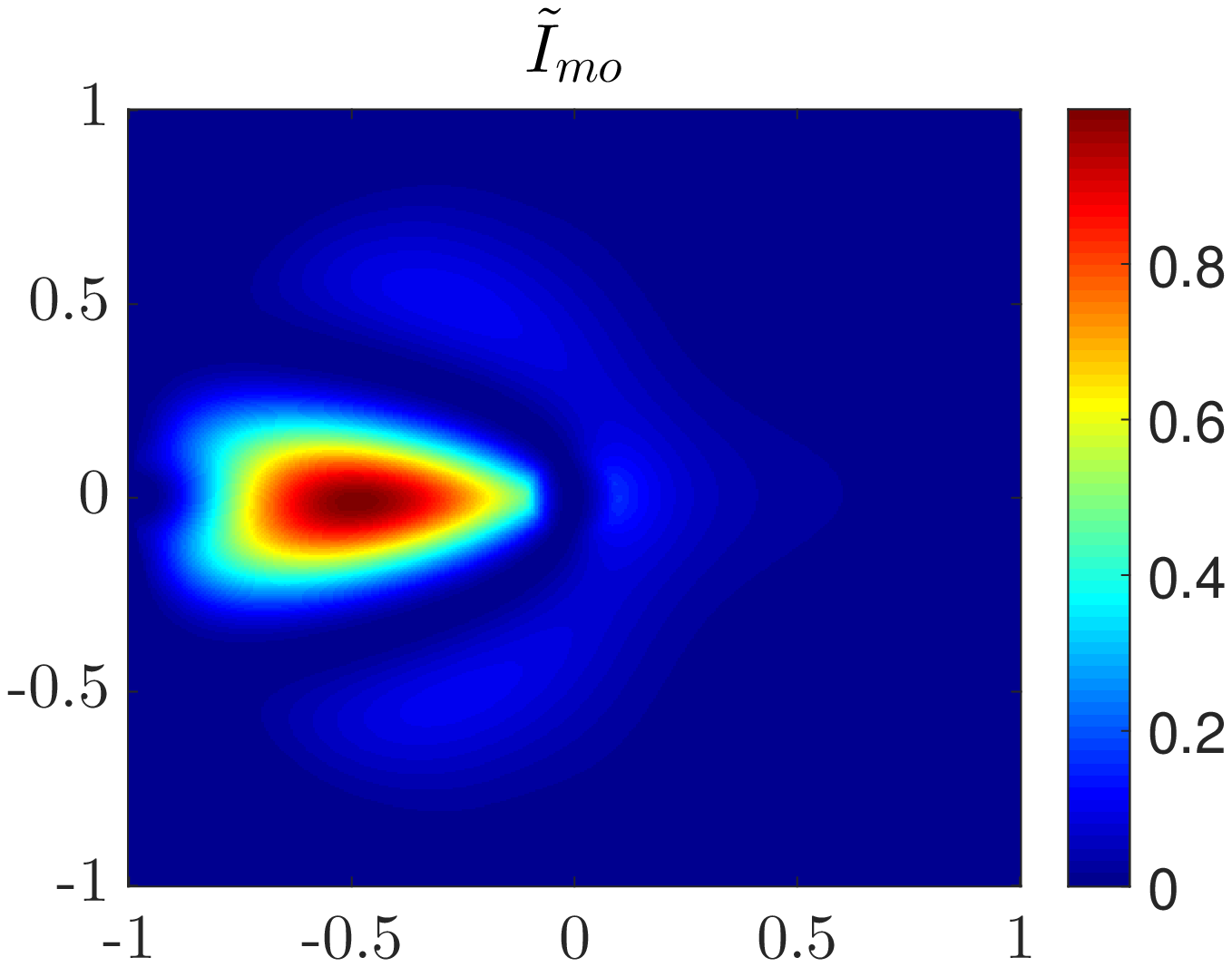}
    \includegraphics[scale = 0.25]
        {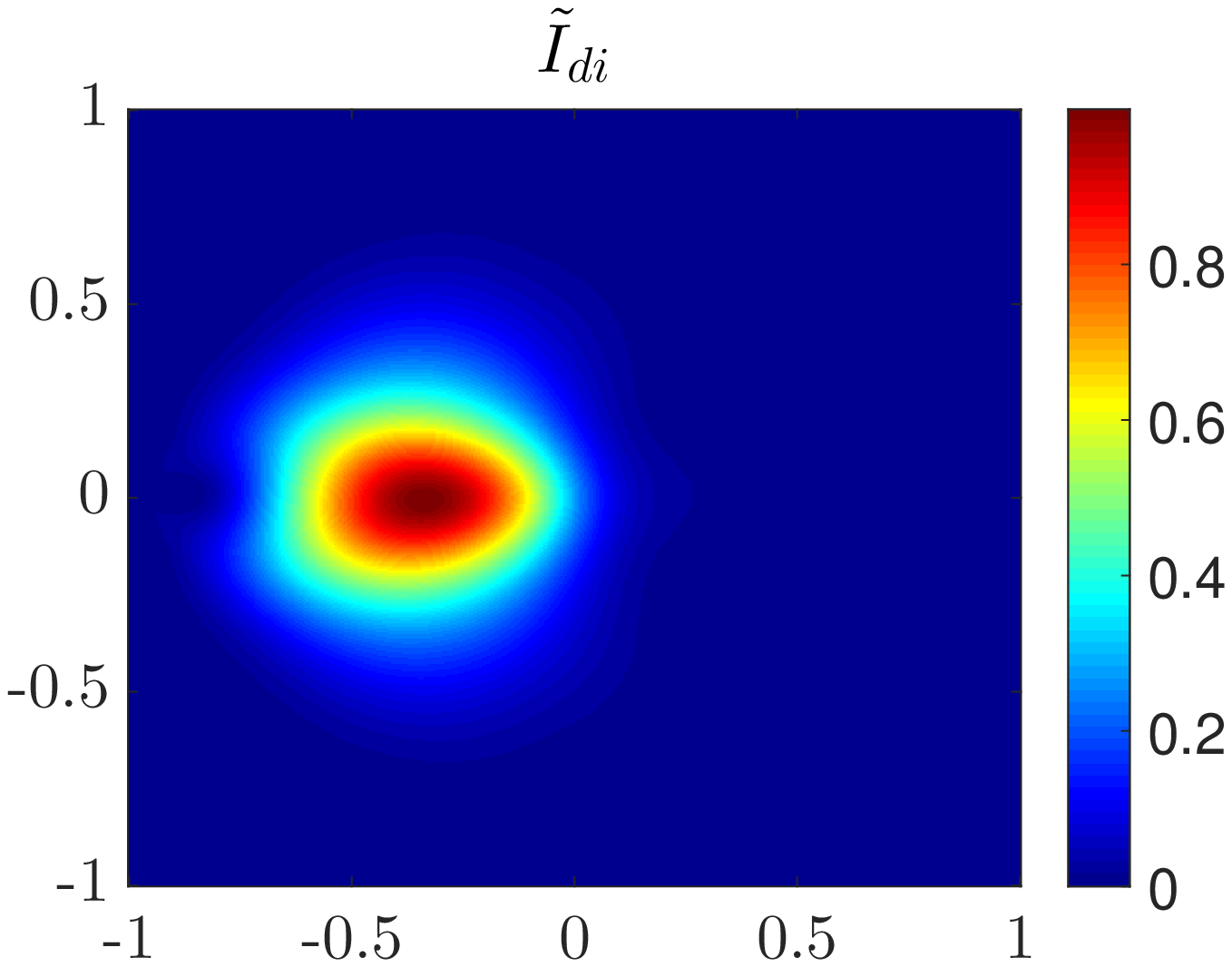}
    \end{subfigure}
        \begin{subfigure}{3in}
    \includegraphics[scale = 0.25]
        {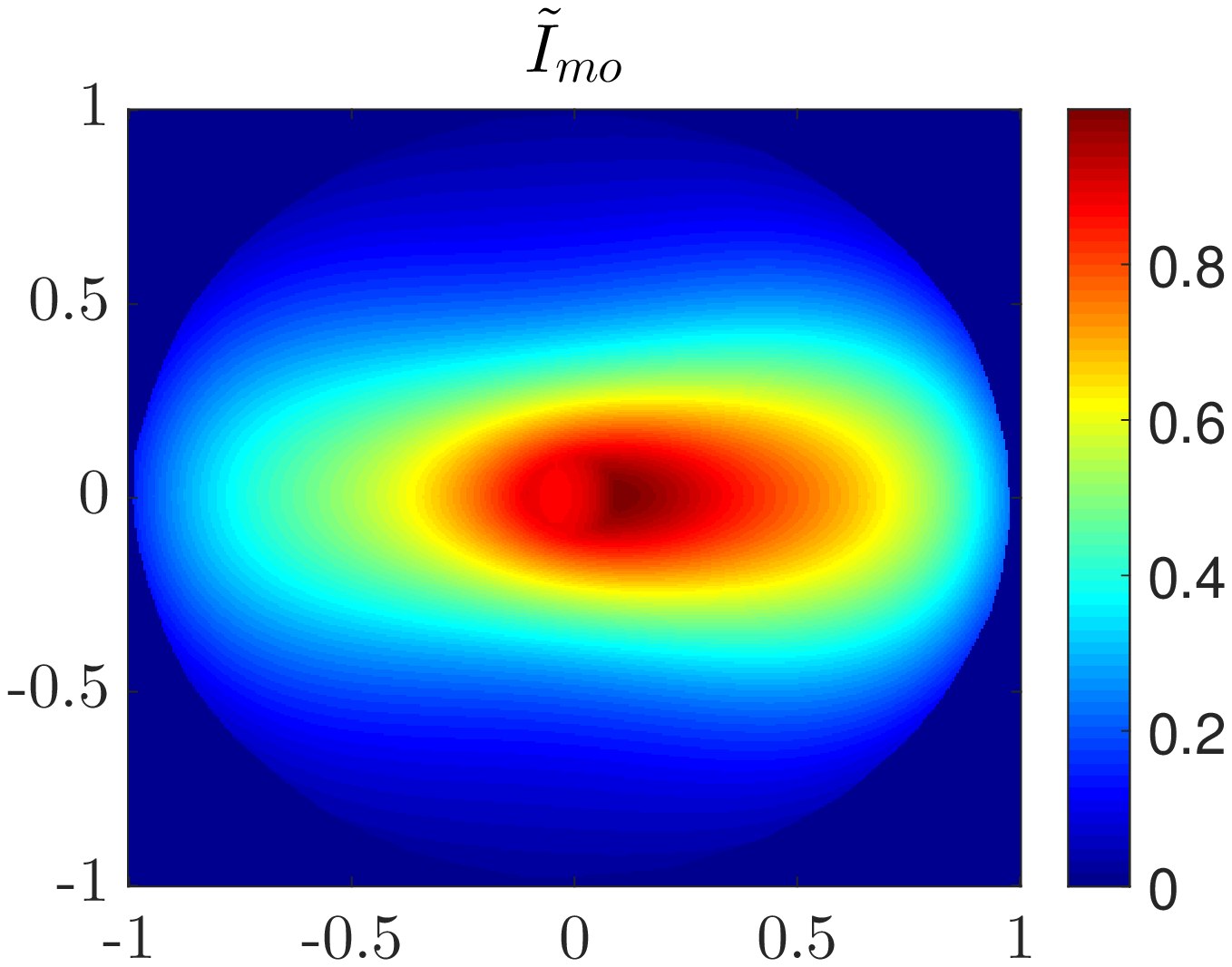}
    \includegraphics[scale = 0.25]
        {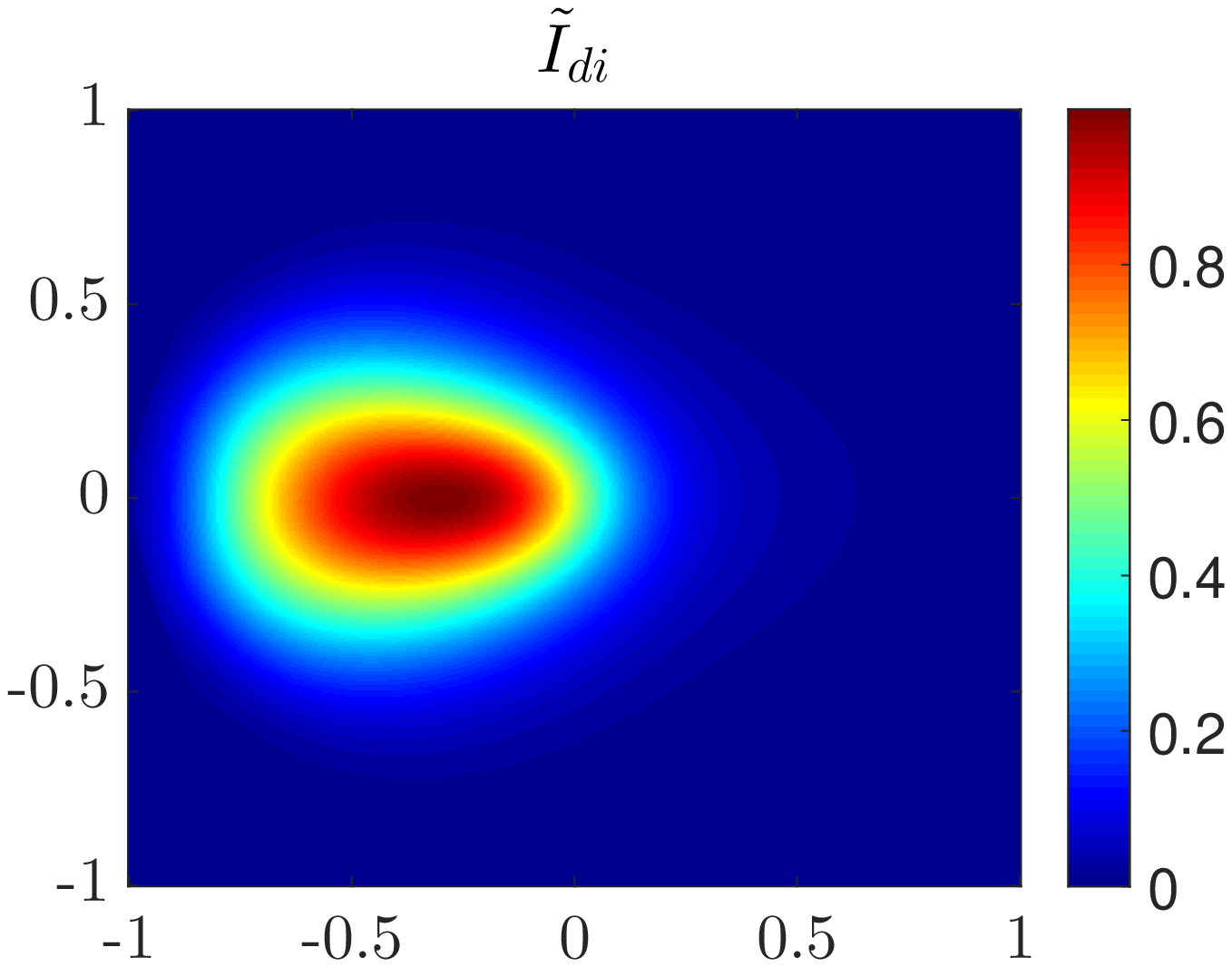}
    \end{subfigure}
    \caption{\textbf{Example 1}. 
    \mm{Top left (exact inclusions): conductivity inhomogeneity (orange),  
    potential inhomogeneity (blue); 
    Top right: monopole index $\tilde{I}_{\text{mo}}$ and dipole index $\tilde{I}_{\text{di}}$, with $f = \cos(\theta)$; 
    Bottom left: $\tilde{I}_{\text{mo}}$  and $\tilde{I}_{\text{di}}$, with $f = \cos(20\theta)$;
    Bottom right: $\gamma=0$. Left: $\tilde{I}_{\text{mo}}$, with $f = \cos(\theta)$; \,Right: 
    $\tilde{I}_{\text{di}} $, with $f = \cos(20\theta)$.}
     }
    \label{case2oneflux}
\end{figure}

\textbf{Example 1.} 
\mm{This example contains two different types of inhomogeneities: 
an inhomogeneity with $\sigma = 1.5$ located at the disk centered at $(-0.4,\,0)$ with radius $0.2$, and another inhomogeneity with $V = 15$ located at the disk centered at $(0.4,\,0)$ with radius $0.2$. 
We apply the boundary influxes of two different frequencies, 
$f_1 = \cos(\theta)$, $f_2 = \cos(20\theta)$, and show their index functions 
$\tilde{I}_{\text{mo}}$ and $\tilde{I}_{\text{di}}$ in Figs.\,\ref{case2oneflux}(b) and \ref{case2oneflux}(c). 
We can see, as the frequency of the boundary influx increases, the reconstruction by $\tilde{I}_{\text{di}} $ 
of the inhomogeneity with $\sigma$ located at left becomes more and more apparent, 
while the reconstruction by $\tilde{I}_{\text{mo}} $ of the inhomogeneity with $V$ located at right 
disappears eventually. 
Fig.\,\ref{case2oneflux} shows the reconstructions with Sobolev index $\gamma = 0$, 
from which we can see the reconstructions are much less sharp than the ones with $\gamma=1$. 
Therefore a non-zero $\gamma$ is essential for a sharper reconstruction. 

Similar numerical effects with the boundary influxes of different frequencies have been observed in many experiments. 
Therefore we will present in all subsequent examples only two measurement events. 
The first measurement is taken with a boundary influx of low frequency, i.e., $f = \cos(\theta)$, 
with which we calculate $\tilde I_{mo}$; the second measurement is taken with a boundary influx of 
high frequency, with which we compute $\tilde I_{di}$. 
}

\subsection{Decoupled reconstructions via the monopole and dipole index functions and 
appropriate choices of boundary influxes}

\mm{We are going to present three representative examples for reconstructing two types of inhomogeneities 
with appropriate choices of boundary influxes based on the strategy we proposed in section\,\ref{subsec_numerical_ver}.
In all our reconstructions for these examples, 
we do not assume any prior knowledge of the shapes, locations and ranges of values of the unknown 
inhomogeneous coefficients $\sigma$ and $V$. 
}

\begin{figure}
	\begin{subfigure}{2in}
    \includegraphics[scale = 0.35]{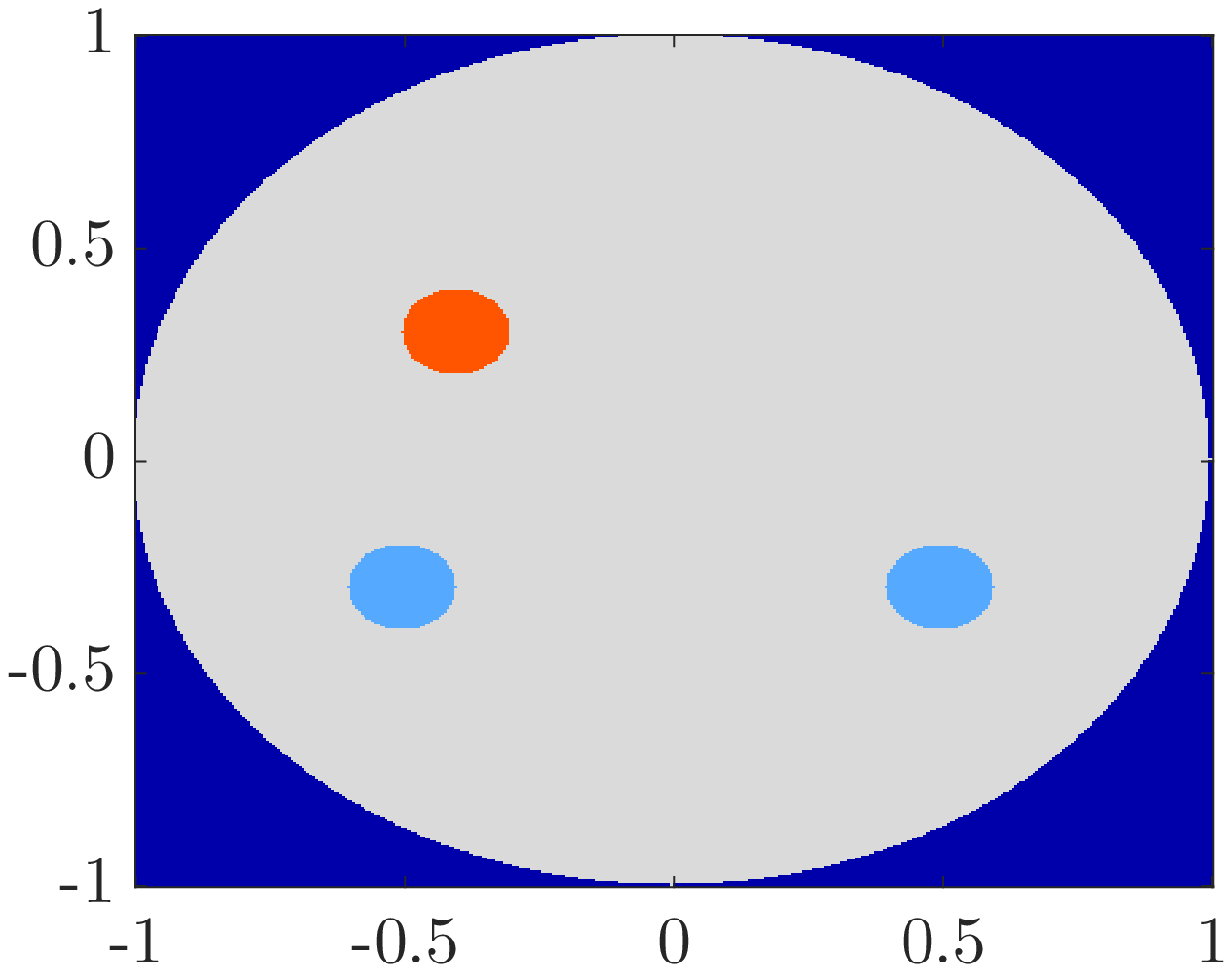}
    \centering
    \end{subfigure}
    \begin{subfigure}{2in}
    \includegraphics[scale = 0.35]{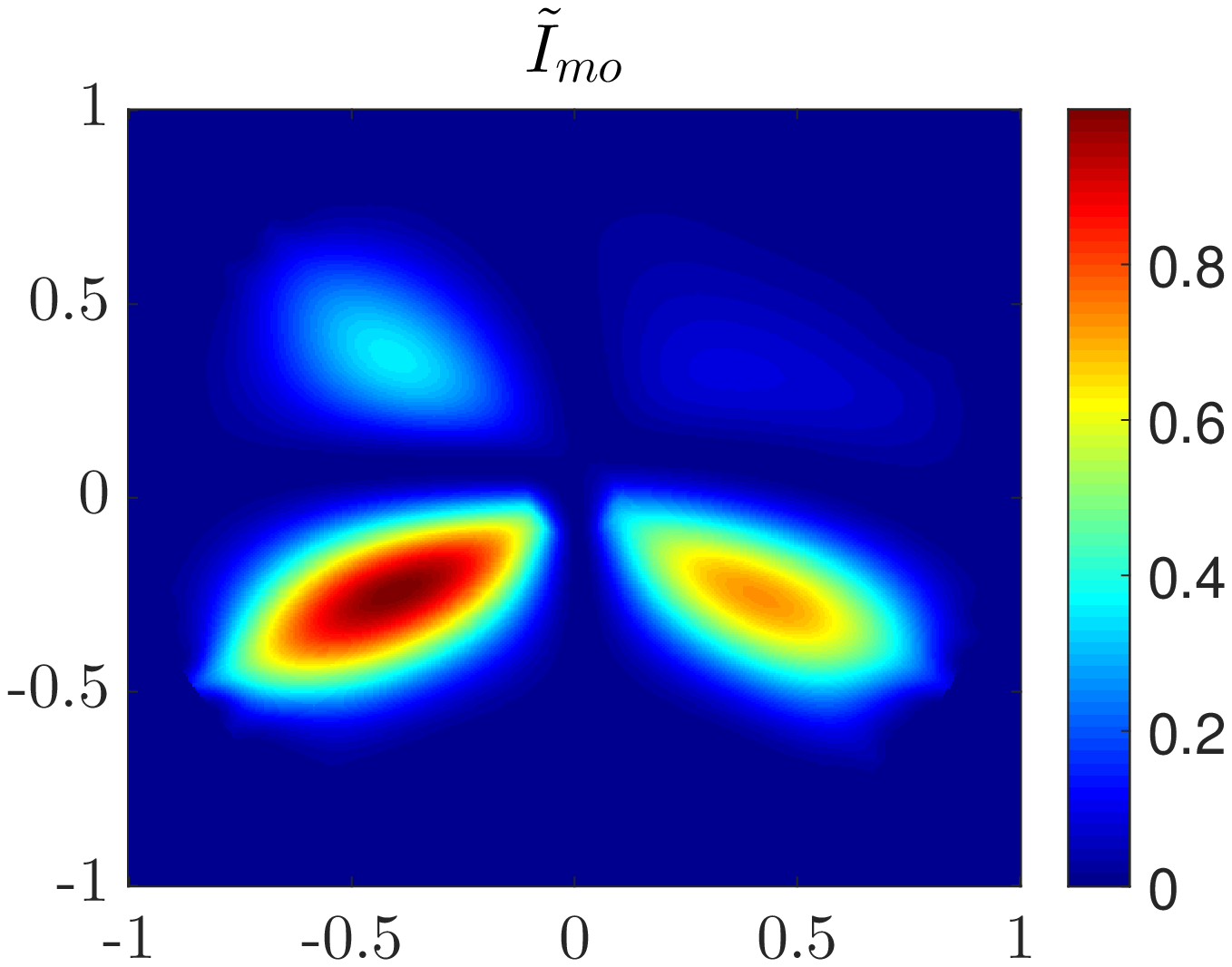}
                \centering
    \end{subfigure}
    \begin{subfigure}{2in}
        \includegraphics[scale = 0.35]{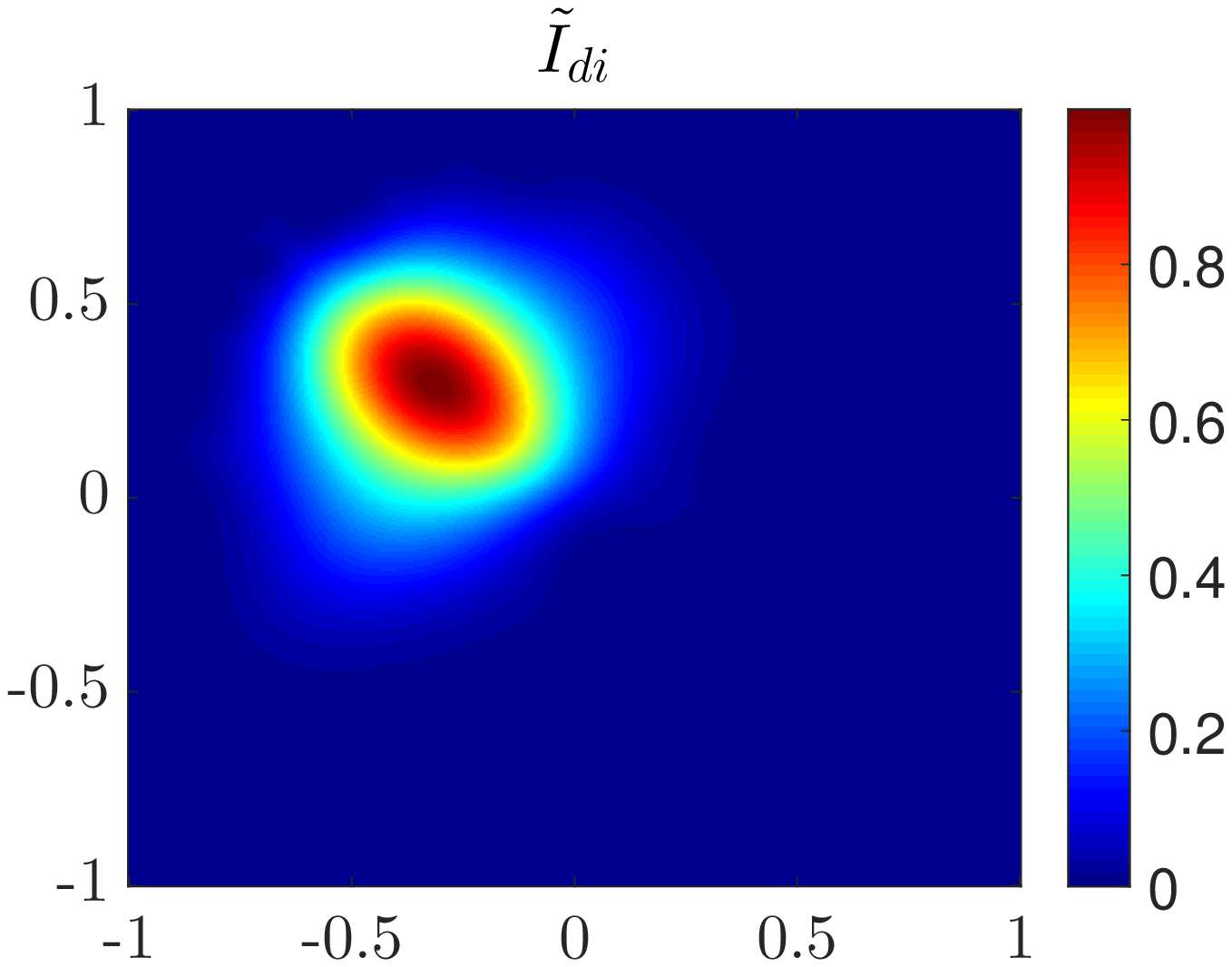}
                \centering
    \end{subfigure}
    \caption{\textbf{Example 2}. 
    \mm{Left (exact inclusions): conductivity inhomogeneity (orange),  
    potential inhomogeneities (blue); Middle: monopole index $\tilde{I}_{\text{mo}}$ with $f = \cos(\theta)$; 
    Right: dipole index $\tilde{I}_{\text{di}}$ with $f = \cos(20\theta)$.}
    }
    \label{case4twoflux}
\end{figure}

\mm{\textbf{Example 2.} In this example, we consider a medium with three inhomogeneities 
as indicated in Fig.\,\ref{case4twoflux}.
As we see, there are two inhomogeneities correspond to the potential $V = 15$, 
located at two disks centered at $(-0.5,-0.3)$ and $(0.5,-0.3)$ with radius $0.1$, respectively, 
and there is another inhomogeneity corresponding to the conductivity $\sigma = 1.5$, 
located at the disk centered at $(-0.4,0.4)$ with radius $0.1$.  
In Fig.\,\ref{case4twoflux}, we have plotted the monopole index $\tilde{I}_{\text{mo}}$ associated with the boundary influx 
$f = \cos(\theta)$, and the dipole index $\tilde{I}_{\text{di}}$ associated with the boundary influx $f = \cos(20\theta)$.  
As one can see from Fig.\,\ref{case4twoflux}, the two different types of inhomogeneities are decoupled: 
$\tilde{I}_{\text{mo}}$ shows the inhomogeneities with $V$, while $\tilde{I}_{\text{di}}$ shows the inhomogeneity 
with $\sigma$. It is surprising that even when the two types of inhomogeneities (both residing in the left part of $\Omega$) 
are very close to each other, the DSM could still separate them clearly. 
}

\begin{figure}
    \centering
        \begin{subfigure}{2in}    \includegraphics[scale = 0.35]{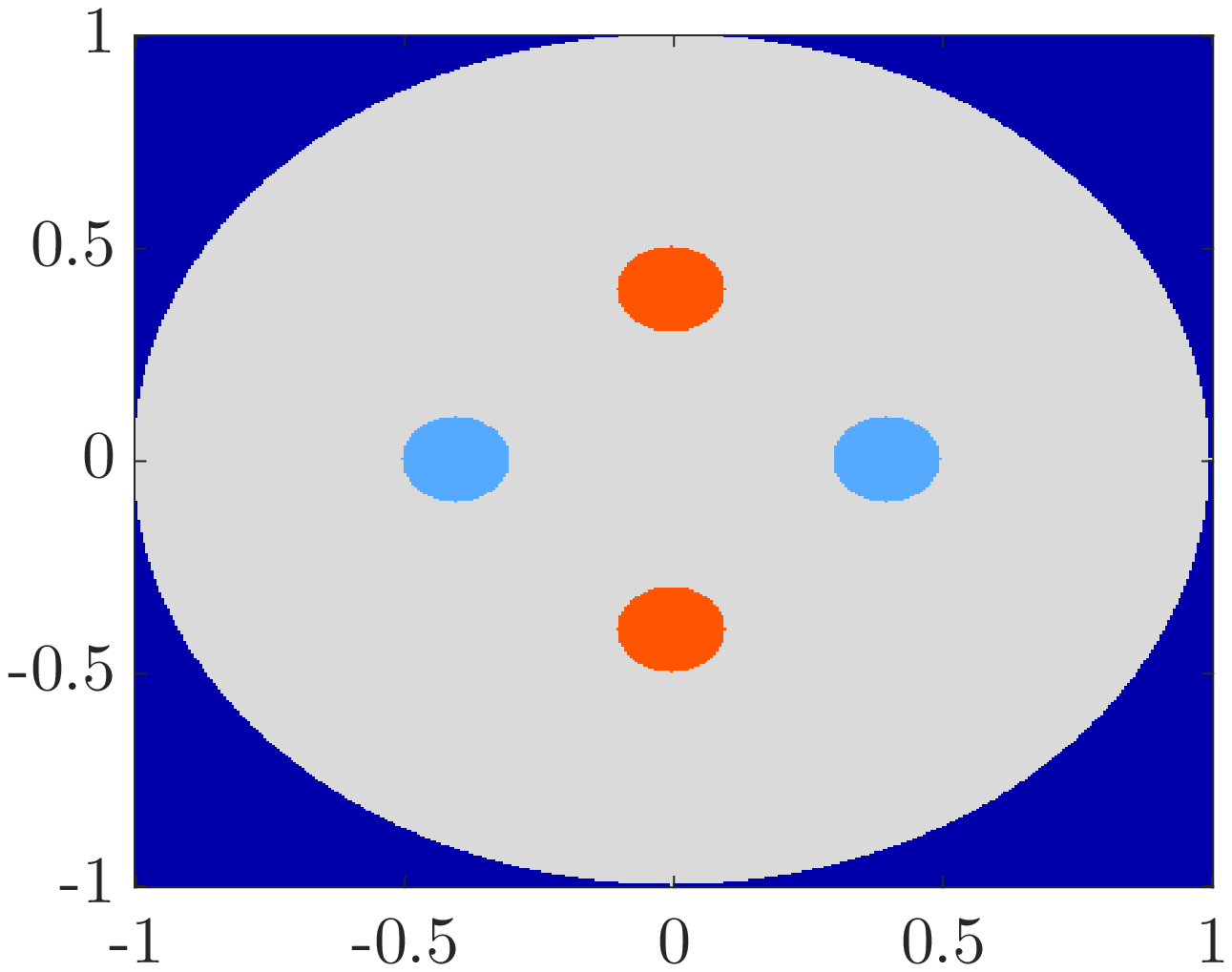}
    \end{subfigure}
    \begin{subfigure}{2in}
    \includegraphics[scale = 0.35]
        {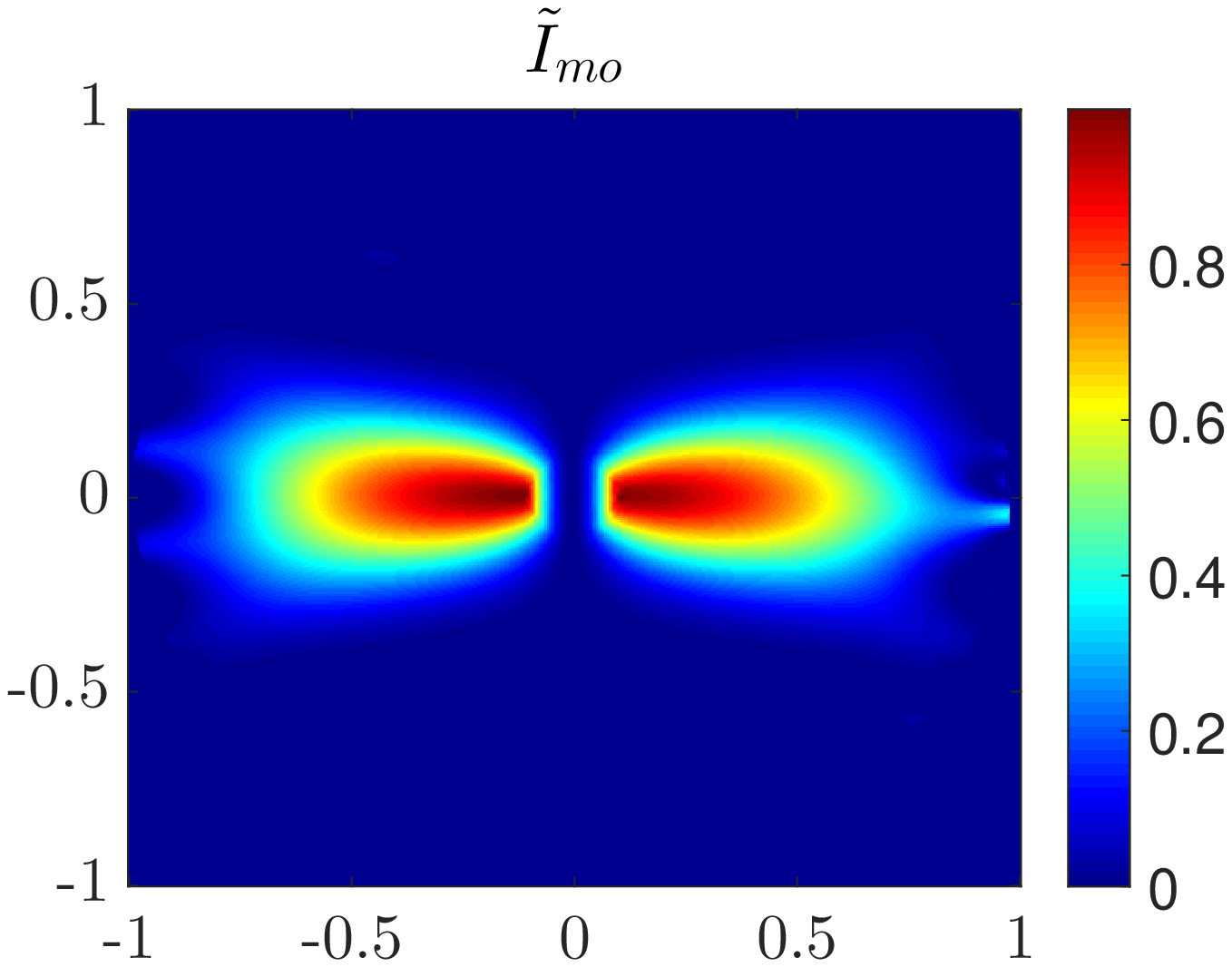}
        \end{subfigure}
  \begin{subfigure}{2in}
        \includegraphics[scale = 0.35]{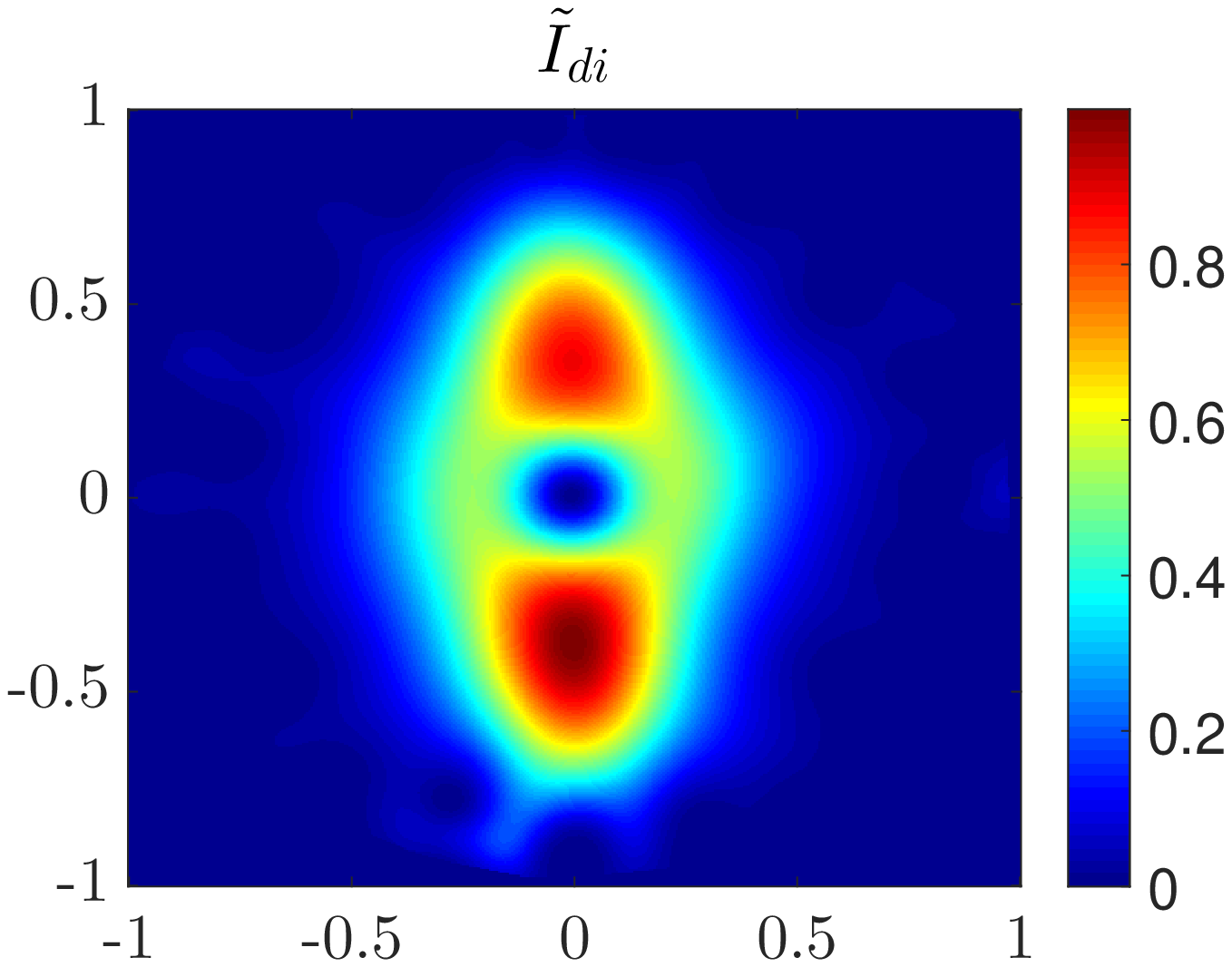}
    \end{subfigure}
    \caption{\textbf{Example 3}. 
    \mm{Left (exact inclusions): conductivity inhomogeneities (orange), potential inhomogeneities (blue); 
    Middle: monopole index $\tilde{I}_{\text{mo}}$ with $f = \cos(\theta)$; 
    Right: dipole index $\tilde{I}_{\text{di}}$ with $f = \cos(20\theta)$.}
    }
    \label{case1twoflux}
\end{figure}

\mm{\textbf{Example 3.} This is a more challenging example with four inhomogeneous inclusions as shown 
in Fig.\,\ref{case1twoflux}. As we see from the figure, 
there are two inhomogeneities corresponding to the conductivity $\sigma = 2.5$, 
located at two disks centered at $(0,0.4)$ and $(0,-0.4)$ with radius $0.1$, respectively; 
meanwhile there are two other inhomogeneities corresponding to the potential $V= 15$, 
located at two disks centered at $(0.4,0)$ and $(-0.4,0)$ with radius $0.1$, respectively.
Fig.\,\ref{case1twoflux} shows the monopole index $\tilde{I}_{\text{mo}}$ associated with the boundary influx 
$f = \cos(\theta)$ and the dipole index $\tilde{I}_{\text{di}}$ associated with the boundary influx $f = \cos(20\theta)$.  
The numerical reconstructions demonstrated the two different types of inhomogeneities are well separated: 
$\tilde{I}_{\text{mo}}$ recovers two inhomogeneities corresponding to $V$, while $\tilde{I}_{\text{di}}$ 
recovers two inhomogeneities corresponding to $\sigma$.  
This shows clearly the success of the DSM in decoupling the measurement data,
locate two different types of inhomogeneous inclusions and distinguish their types quite reasonably. 
}

\begin{figure}
    \begin{subfigure}{2in}
    \includegraphics[scale = 0.35]{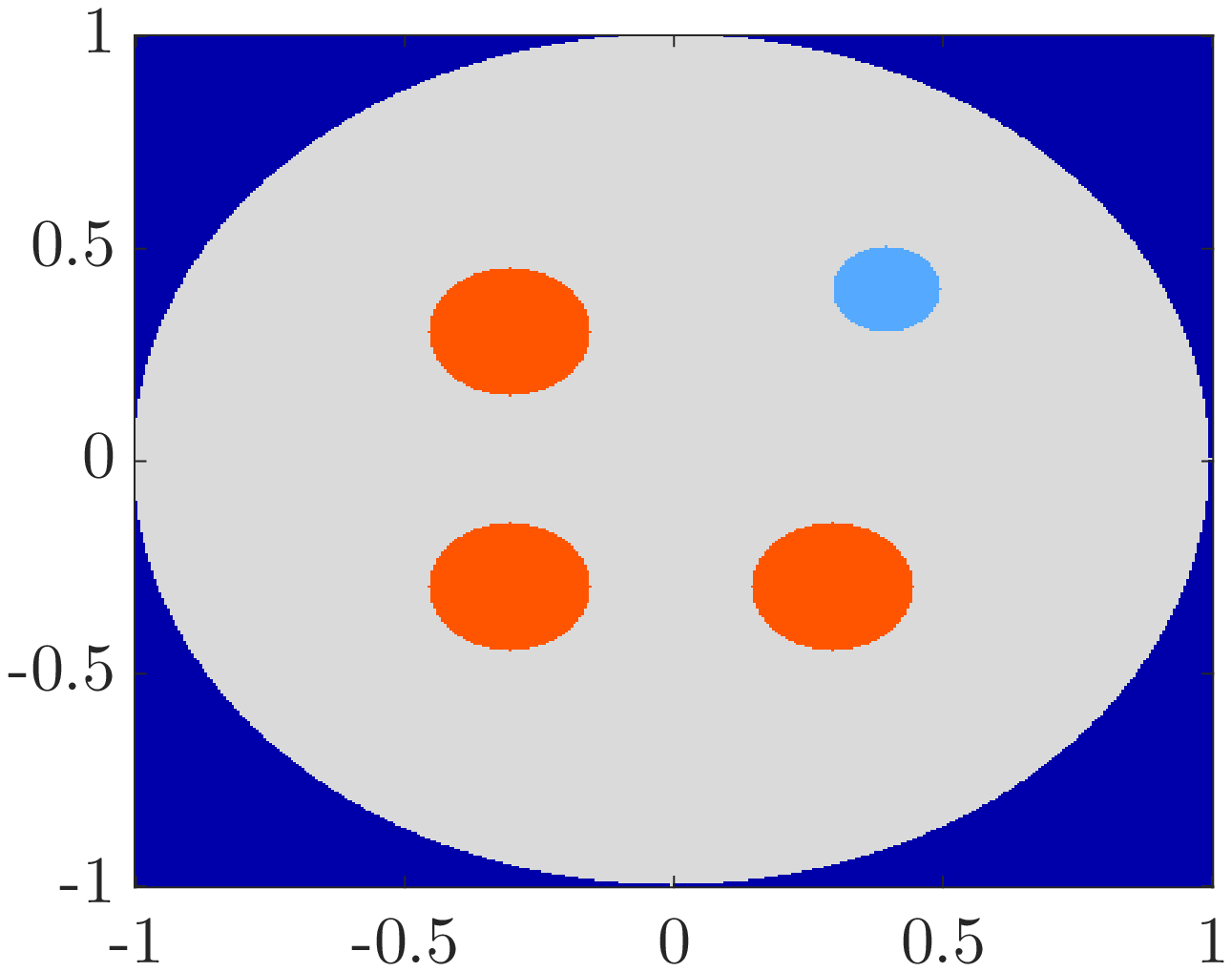}
    \centering
    \end{subfigure}
    \begin{subfigure}{2in}
    \includegraphics[scale = 0.35]{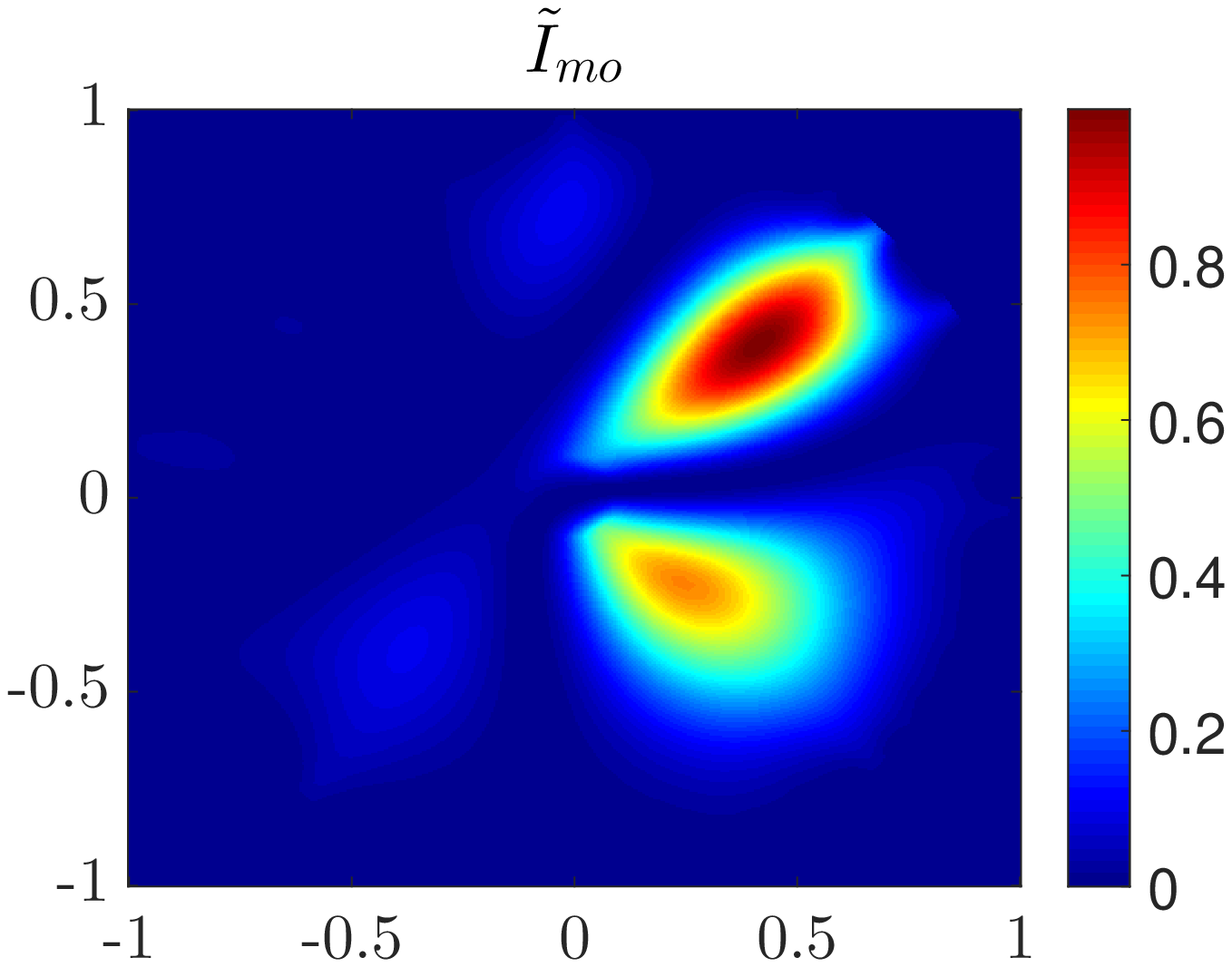}
                \centering
    \end{subfigure}
        \begin{subfigure}{2in}
        \includegraphics[scale = 0.35]{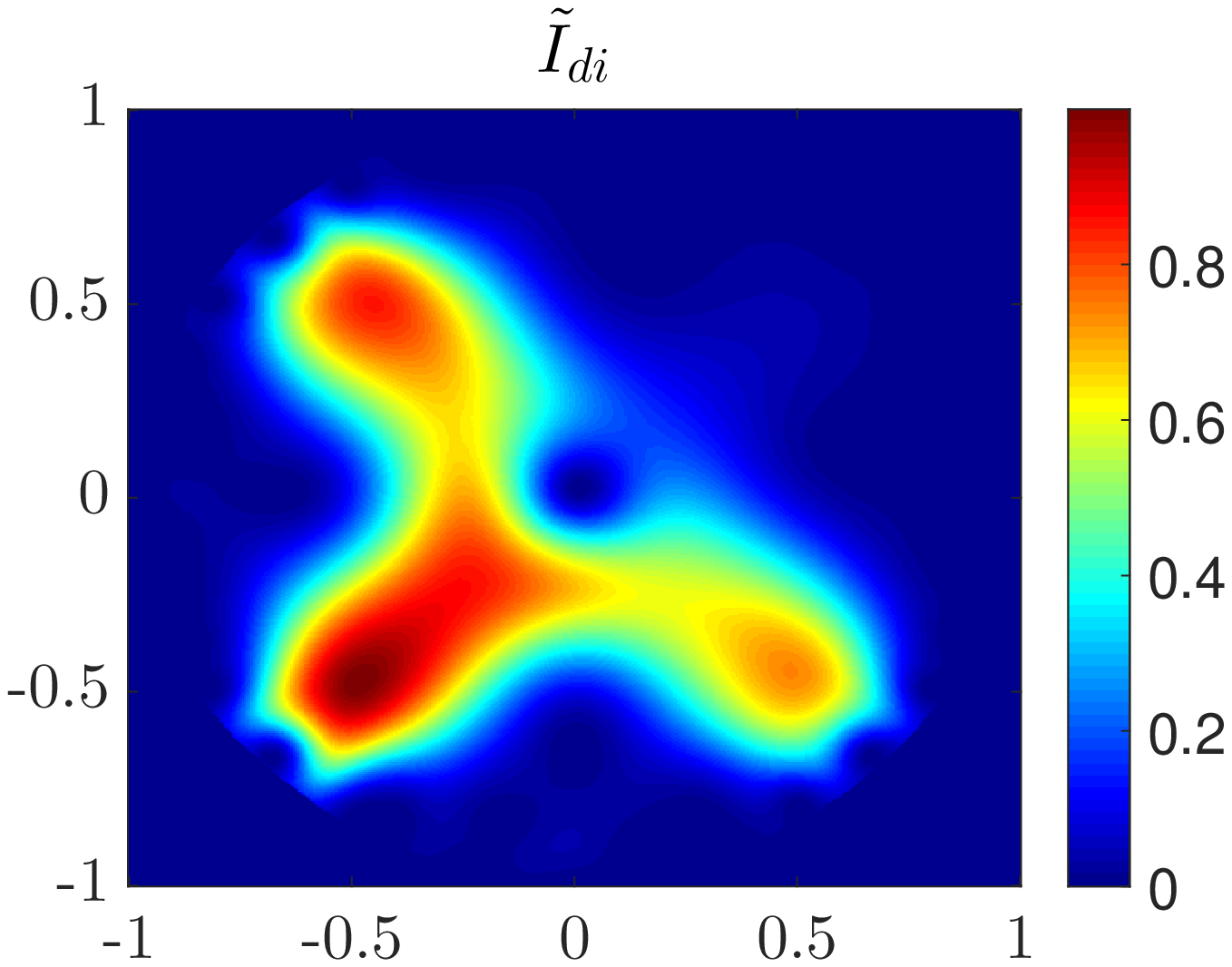}
                \centering
    \end{subfigure}
        \caption{\textbf{Example 4}. 
        \mm{Left (exact inclusions): conductivity inhomogeneities (orange), potential inhomogeneity (blue); 
        Middle: monopole index $\tilde{I}_{\text{mo}}$, with $f = \cos(\theta)$;
        Right: dipole index $\tilde{I}_{\text{di}}$, with $f = \cos(30\theta)$.
        }
}
    \label{case7}
\end{figure}

\textbf{Example 4.}
\mm{This example shows a medium with four inhomogeneous inclusions as in Fig.\,\ref{case7}.
We see three conductivity inhomogeneities with $\sigma = 2$ placed at three disks centered at $(-0.3,0.3)$, $(0.3,-0.3)$, 
and $(-0.3,-0.3)$ with radius $0.15$, and one potential inhomogeneity with $V = 22$ 
placed at the disk centered at $(0.4,0.4)$ with radius $0.1$.
Fig.\,\ref{case7} plots the monopole index $\tilde{I}_{\text{mo}}$ with the boundary influx $f = \cos(\theta)$ 
and the dipole index $\tilde{I}_{\text{di}}$ with the boundary influx $f = \cos( 30\theta)$. 
This example is quite surprising to see a satisfactory separation of the conductivity inhomogeneous inclusions 
from the potential inhomogeneities although the number of the former is three times of the latter.
We can further improve the sharpness of $\tilde{I}_{di}$ when the data is collected at more measurement points.
}

 \section{Concluding remarks}
We have proposed a novel direct sampling method for simultaneously reconstructing two different types of inhomogeneities 
inside a domain with boundary measurements collected from only one or two measurement events.  
This inverse problem is theoretically known to have no uniqueness in most cases, and is highly unstable and ill-posed.

A main feature of the new method is to design two distinct sets of probing functions, i.e., 
the monopole and dipole probing functions, which help decouple the respective signals coming from 
the monopole-type and dipole-type sources located in the sampling domain. Each type of sources carries the information of one distinctive type of inhomogeneity we aim to reconstruct. This enables us to decouple the boundary measurements and achieve reasonable simultaneous reconstructions. The direct sampling method relies on two index functions that can be computed 
in a fast, stable and highly parallel manner.
Numerical experiments have illustrated its stability in decomposing 
different signals coming from two types of inhomogeneities in measurement data, and 
its robustness against noise.

Our choice of the model inverse problem covers a general class of inverse coefficients problems that we encountered in applications, for instance, diffusion-based optical tomography, inverse electromagnetic scattering problem under transverse symmetry and ultrasound medical imaging. A very unique feature of the new method is its applications to the important 
scenarios when very limited data is available, e.g., only the data from one or two measurement event, to which most existing 
methods are not applicable.

Along this research topic, there are some interesting and important directions that deserve further exploration: 
extend the sampling method to a broader class of coefficients inverse problems with more complicated 
interaction terms, for instance, anisotropic electromagnetic scattering, fully anisotropic linear and nonlinear elasticity model, 
shallow water wave equation, Boltzmann transport equation, Klein-Gordon and Sine-Gordon equations, etc.; 
develop a unified framework of the direct sampling methods, with a concrete recipe for generating optimal probing functions 
and duality products for a given inverse problem.

%
%
%
%
%
%

\bibliographystyle{siamplain}
\nocite{*}
\bibliography{dsm_decoupling}
\citation

\end{document}